# Harmonic Analysis on the quantum Lorentz group

E. Buffenoir[*][†]
Laboratoire de Physique Mathématique et Théorique[†]
Université Montpellier 2, Place Eugène Bataillon
34000 Montpellier France

Ph. Roche[‡]
TH Division CERN,
CH-1211 Geneva 23 Switzerland
*on leave from*
CPT Ecole Polytechnique[§]
91128 Palaiseau Cedex France

January 18, 1999

**Abstract**

This work begins with a review of complexification and realification of Hopf algebras. We emphasize the notion of multiplier Hopf algebras for the description of different classes of functions (compact supported, bounded, unbounded) on complex quantum groups and the construction of the associated left and right Haar measure. Using a continuation of $6j$ symbols of $SU_q(2)$ with complex spins, we give a new description of the unitary representations of $SL_q(2,\mathbb{C})_{\mathbb{R}}$ and find explicit expressions for the characters of $SL_q(2,\mathbb{C})_{\mathbb{R}}$. The major theorem of this article is the Plancherel theorem for the Quantum Lorentz Group.

[*]e-mail:buffenoi@lpm.univ-montp2.fr
[†]Laboratoire du CNRS ESA 5032
[‡]e-mail:philippe.roche@CERN.CH
[§]Laboratoire Propre du CNRS UPR 14

# 1 Introduction

Our interest in the subject of harmonic analysis on $SL_q(2,\mathbb{C})_{\mathbb{R}}$ stems from our desire to apply the program of combinatorial quantization of Chern-Simons theory [18, 1, 9] to the case where the quantum group is a quantization of a non compact group. By $SL_q(2,\mathbb{C})_{\mathbb{R}}$ we mean a quantization, for $q$ real, of the group $SL(2,\mathbb{C})$ considered as a real Lie group. We have kept in mind that there is a close relationship between quantum gravity in $2+1$ dimensions and Chern-Simons theory with a non compact group of the type $SL(2,\mathbb{C})_{\mathbb{R}}, SL(2,\mathbb{R}) \times SL(2,\mathbb{R}), ISO(2,1)$ (depending on the sign of the cosmological constant) as first shown in [46] (which is a $2+1$ analogue of Ashtekar variables [4]). Any result in the project of quantization of Chern-Simons theory for these groups has spin-offs on the program of canonical quantization of $2+1$ quantum gravity (for a review on 2+1 quantum gravity, see [11]). Moreover the Chern-Simons topological invariants associated to $SL(2,\mathbb{C})_{\mathbb{R}}$ appear to be essential in the description of states in the Ashtekhar's program of canonical quantization of $3+1$ dimensional gravity. It will be clear in this paper that the "non-compacity" of $SL(2,\mathbb{C})_{\mathbb{R}}$ is at the origin of many properties which cannot be so easily translated from the compact case.
There are various ways to investigate the problem of quantization of Chern-Simons theory with a non compact group, two main contributions being the work of J.E. Nelson and T. Regge on the direct quantization of the algebra of observables [33] and the work of E. Witten on the geometric quantization of the space of flat connections for a complex group [47] (see also [26] in the case of $SL(2,\mathbb{R})$). If we could apply the program of combinatorial quantization, and particularly the results of [2], in the case where the group is non compact, we would be able to go much further than the approach of Nelson and Regge and obtain an Hilbert space of states, and a unitary representation of the algebra of observables on the space of states for any punctured surface using only representation theory of the associated quantum group. Unfortunately for the simplicity and fortunately for the richness of the theory, there are numerous problems to handle in order to use the combinatorial quantization program in these cases. All the difficulties are connected to the fact that harmonic analysis on quantization of non compact groups is completely different from the compact case (as in the classical case) and in the quantum case it is still not very much developed.

Harmonic analysis on compact quantum groups has been developed in the context of the compact matrix pseudo group by Woronowicz [48]. Its culminating result is the existence of a left and right integral and an analogue of Peter-Weyl theorem. We can say that in this case, functional analysis is restricted to very little and the theory is almost completely algebraic. The case of $SU_q(1,1)$ is now completely understood: classification of unitary representations of $SU_q(1,1)$ [32], proof of Plancherel Theorem [25, 28]. Harmonic analysis on quantization of non compact groups has made a lot of progress after the work of Woronowicz who recognized the central role of the theory of multipliers and affiliated elements of a non unital $C^\star$-algebra [50]. In their very important work [36], Podles and Woronowicz gave three main results: the construction of $\mathfrak{U}_q(sl(2,\mathbb{C})_{\mathbb{R}})$ as a quantum double of $\mathfrak{U}_q(su(2))$, an explicit Iwasawa decomposition that allows the construction of different classes of functions on $SL_q(2,\mathbb{C})_{\mathbb{R}}$ and the proof of the existence of a unimodular Haar measure. Their study is not restricted to the case of $SL_q(2,\mathbb{C})_{\mathbb{R}}$ and concerns all quantum doubles of any quantization of compact groups.

The aim of our work is twofold: the first goal is to obtain a complete description of the characters of $SL_q(2,\mathbb{C})_{\mathbb{R}}$ and prove a Plancherel formula. The second goal is to mix different points of view on quantum group theory: the $C^\star$ algebra approach developed by Woronowicz and his collaborators, and the $R-$matrix approach developed by the Russian school. We have tried to think of $\mathfrak{U}_q(sl(2,\mathbb{C})_{\mathbb{R}})$ as being the double of $\mathfrak{U}_q(su(2))$ and to forget about commutation relations; this has the advantage that all the results can be recast in terms of generalized $6j$ symbols and that the proofs are mainly graphical. We not only gain in clarity but obtain results like the description of characters and the Plancherel theorem, which would have been hard to obtain outside this formulation.
It will be clear to the reader that our constructions and proofs can be generalized, with some work, to the case of the quantization of any complex group. We have only analysed the case where $q$ is real. Although the study of $SL_q(2,\mathbb{C})_{\mathbb{R}}$, when $q$ is any complex number, is of great physical relevance (for quantum Liouville theory and for $2+1$ quantum gravity with a gravitational Chern-Simons term [47]), it appears that the general case requires quasi-Hopf algebras.
Note that we have abusively called in our work $SL_q(2,\mathbb{C})_{\mathbb{R}}$ the quantum Lorentz group although classically $SL(2,\mathbb{C})_{\mathbb{R}}$ is the simply connected covering of $SO(3,1)$. In [36] a quantization of $SO(3,1)$ is proposed, which is the quantum analogue of the relation $SO(3,1)=SL(2,\mathbb{C})_{\mathbb{R}}/\mathbb{Z}_2$.

# 2 Complexification and realification of Hopf algebras

In this chapter we have collected and given a symplifying treatment of results which are scattered in many places. We have tried to be very pedagogical and to clarify notions which appear quite confusing. In particular we have tried to distinguish the notion of Hopf algebra over $\mathbb{R}$ from the notion of star structure on the complexification. We have made an analysis for the case of star structure which are morphism of coalgebra and also for the other case of antimorphism of coalgebra, which is also of interest. We have distinguished the algebra $A$ from its conjugate $\overline{A}$ : these two algebras are often taken to be isomorphic which obscures the overall structure of realification. In particular we stress that the construction of a realification of a Hopf algebra $A$ is independent of a choice of star structure on $A$.

## 2.1 Star structures

In this section we give a self contained exposition of general properties of real forms of Hopf algebras. For a recent review of this subject see [31] and references therein.
We will first describe classical and standard notions such as realification, complexification.

Let $V$ be a $\mathbb{C}$ vector space, the realification of $V$ is the $\mathbb{R}$ vector space $V_{\mathbb{R}}$ obtained from $V$ by restricting the field $\mathbb{C}$ to $\mathbb{R}$. The vector space $\bar{V}$ is the vector space defined by $\bar{V}=V$ as an additive group and $\lambda\bar{x}=\overline{\bar{\lambda}x}$ where $x,\lambda\in V\times\mathbb{C}$ and $\bar{\lambda}$ is the complex conjugate of $\lambda$.

Let $W$ be a $\mathbb{R}$ vector space, the complexification of $W$ is the $\mathbb{C}$ vector space $W^{\mathbb{C}}=W\otimes_{\mathbb{R}}\mathbb{C}$. We trivially have $(W^{\mathbb{C}})_{\mathbb{R}}=W\oplus iW$ as $\mathbb{R}$ vector space and less trivially $(V_{\mathbb{R}})^{\mathbb{C}}=V\oplus\bar{V}$ as $\mathbb{C}$ vector space, where the latter isomorphism is given by $x+iy\mapsto(x+Jy,\overline{x-Jy})$ $x,y\in V$ where $J$ denotes the action of $\mathbb{C}$ on $V$ and $i$ denotes the formal variable used to describe complexification.

These notions and results generalize straightforwardly to the category of Lie algebras and associative algebras. For example, if $\mathfrak{G}$ is a Lie algebra over $\mathbb{R}$, there is a unique Lie algebra structure on the complexification $\mathfrak{G}^{\mathbb{C}}$ of $\mathfrak{G}$ given by $[x+iy,x'+iy']_{\mathfrak{G}^{\mathbb{C}}}=([x,x']_{\mathfrak{G}}-[y,y']_{\mathfrak{G}})+i([x,y']_{\mathfrak{G}}+[y,x']_{\mathfrak{G}})$. The complex conjugate $\overline{\mathfrak{G}}$ of a complex Lie algebra $\mathfrak{G}$ verifies $[\overline{x},\overline{y}]=\overline{[x,y]}$. If $A$ is an algebra, then $\overline{A}$ denotes the complex conjugate algebra of $A$ with algebra law given by $\overline{a}\overline{b}=\overline{ab}$. If $\pi$ is a representation of $A$ acting on a vector space $V$, then we can define a representation $\overline{\pi}$ of $\overline{A}$ acting on $\overline{V}$ by the formula: $\overline{\pi}(\overline{a})\overline{v}=\overline{\pi(a)(v)}$, $(a,v)\in A\times V$.

However these notions do not generalize so directly to the category of coalgebras. The process of complexification does not cause any problems: if $W$ is a $\mathbb{R}$ coalgebra and $\Delta:W\to W\otimes_{\mathbb{R}}W$ is the coproduct, then we can define on $W^{\mathbb{C}}$ a unique $\mathbb{C}$ coalgebra structure $\Delta^{\mathbb{C}}$ by extending $\Delta$ to $\Delta^{\mathbb{C}}$ by $\mathbb{C}$ linearity.
It is the process of realification which cannot be straightforwardly generalized. Indeed if $V$ is a $\mathbb{C}$ coalgebra, and $\Delta : V \to V \otimes_{\mathbb{C}} V$ is the coproduct, then the restriction of the field naturally defines a vector space $V_{\mathbb{R}}$ but the restriction of the coproduct that we will denote $\Delta_{\mathbb{R}}$ maps $V_{\mathbb{R}}$ to $(V \otimes_{\mathbb{C}} V)_{\mathbb{R}}$ and this space is not canonically embedded in $V_{\mathbb{R}} \otimes_{\mathbb{R}} V_{\mathbb{R}}$. So unless $\Delta$ satisfies very specific properties it is not possible to endow $V_{\mathbb{R}}$ with a structure of $\mathbb{R}$ coalgebra.
This simple remark indicates that there must be problems in constructing real coalgebras from complex ones, this problem will reappear recurrently in the Hopf algebra context.

Let us now introduce the important notion of star structure which selects a real form of a complex vector space. If $W$ is a $\mathbb{C}$ vector space we can ask the following question: what is the necessary and sufficient condition on $W$ in order that $W = V^{\mathbb{C}}$ where $V$ is a real vector space. The answer is that $W$ has to be endowed with an antilinear involutive map $\star : W \to W$, $V$ being defined as $V = \{x \in W, x^\star = -x\}$ (we could also have defined $V$ to be the eigenvectors of $\star$ of value 1). Of course, if $V$ is a vector space with a star structure then $\star$ is a $\mathbb{C}$ isomorphism between $V$ and $\overline{V}$. If $W, V$ are two vector spaces with star structures we will define a star structure on $V \otimes_{\mathbb{C}} W$ by $(a \otimes b)^\star = a^\star \otimes b^\star$. An important example is the case where $\mathfrak{G}$ is a complex vector space, and $V = \mathfrak{G}_{\mathbb{R}}$. Using the isomorphism $W = V^{\mathbb{C}} = \mathfrak{G} \oplus \overline{\mathfrak{G}}$, it is easy to show that the star structure is given by $(x \oplus \overline{y})^{\bigstar} = (-y \oplus -\overline{x})$.
This notion, being trivial in the vector space case is really important in the case of Lie algebras, and the following proposition plays a central role: let $W$ be a $\mathbb{C}$ Lie algebra, a real Lie algebra $V$ satisfying $W = V^{\mathbb{C}}$ is equivalently defined by an antilinear involutive antimorphism $\star : W \to W$, $V$ being defined by $V = \{x \in W, x^\star = -x\}$ (if we have chosen $\star$ to be an antimorphism, which is the standard choice, the sign has to be -).
If $A = \mathfrak{U}(W)$ is the universal envelopping algebra of $W$, a real Lie algebra $V$ is equivalently defined by an antilinear involutive antimorphism $\star : A \to A$, satisfying in addition the condition

$$(C) \ \ \forall a \in A, \ (\star \otimes \star)\Delta(a) = \Delta(a^\star), \tag{1}$$

$\mathfrak{U}(V)$ being defined by $\mathfrak{U}(V) = \{a \in A, a^\star = S(a)\}$ where $S = S^{-1}$ is the antipode of the universal envelopping algebra. $V$ is called a real form of $W$, and can also be defined by $V = \{a \in A, a^\star = -a, \Delta(a) = a \otimes 1 + 1 \otimes a\}$. In the important case where $A = \mathfrak{U}(\mathfrak{G})$ is the universal envelopping algebra of a real Lie algebra $\mathfrak{G}$, we have $A^{\mathbb{C}} = \mathfrak{U}((\mathfrak{G})^{\mathbb{C}}) = \mathfrak{U}(\mathfrak{G}) \otimes \overline{\mathfrak{U}(\mathfrak{G})}$ (the two copies are commuting) and the natural star structure $\bigstar$ is given by $(x \otimes \overline{y})^{\bigstar} = S(y) \otimes \overline{S(x)}$.
In the Hopf algebra context the standard philosophy is to define real form of a complex Hopf algebra as being a choice of star structure on the complex Hopf algebra and not to define them as being real Hopf algebras. This has the advantage of enlarging the notion of real form but constructions such as realification still have to be defined and will be studied in the following section.
We will choose the following definition: let $A$ be a complex Hopf algebra, a star structure on $A$ of type (C+) (resp. (C−)) is an antilinear involutive antimorphism $\star : A \to A$ and an invertible element $F_\pm \in A \otimes A$ such that :

$$\Delta(a^\star) = F_+(\star \otimes \star)(\Delta(a))F_+^{-1} \quad (C+), \qquad (\text{resp.} \qquad \Delta(a^\star) = F_-(\star \otimes \star)(\Delta'(a))F_-^{-1} \quad (C-) \ ).$$

When $A$ is quasitriangular the condition $C+$ and $C-$ are equivalent and we can choose the relation $F_+ = F_- R$. One motivation for this general definition involving $F_\pm$ is provided by a theorem of Drinfeld on quasi-Hopf algebra, stating that $\mathfrak{U}_\hbar(\mathfrak{G})$ is isomomorphic as quasi-Hopf algebra to $(\mathfrak{U}(\mathfrak{G}), \phi, R = exp(\hbar t))$ twisted by $F \in \mathfrak{U}(\mathfrak{G})^{\otimes 2}$. As a result if $\star$ is any involution selecting a real form on $\mathfrak{G}$ this $\star$ can be extended on $\mathfrak{U}(\mathfrak{G})$ and pulled back on $\mathfrak{U}_\hbar(\mathfrak{G})$. This $\star$ is trivially an involutive antimorphism of algebra, and satisfies condition $(C+)$ with $F_+ = F F^{\star \otimes \star}$.

Unfortunately no simple formula for $F$ is known up to now.
This definition appears restricted in the litterature and usually only $F_+ = 1$ and $F_- = 1$ are considered, this is also the choice we will take in this article.
This restriction has two drawbacks, first it will generate an assymmetry in the quasitriangular case between $(C+)$ and $(C-)$ and it will restrict severely the value of the deformation parameter.
It is easy to show that the relationship between the antipode and the $\star$ is as follows:

$$S \circ \star = \star \circ S, \text{ in the case (C}-) \tag{2}$$

$$S \circ \star = \star \circ S^{-1}, \text{ i.e } (S \circ \star)^2 = \text{id}, \text{ in the case (C+)} \tag{3}$$

As an example let us consider $A = \mathfrak{U}_q(sl(2,\mathbb{C}))$, where $q$ is in $\mathbb{C}^*$, it can be easily shown that in the case $q$ real there are three real forms $\mathfrak{U}_q(su(2))$ (C+), $\mathfrak{U}_q(su(1,1))$ (C+), $\mathfrak{U}_q(sl(2,\mathbb{R}))$ (C$-$) and in the case where $|q| = 1$ there are three other real forms $\mathfrak{U}_q(su(2))$ (C$-$), $\mathfrak{U}_q(su(1,1))$ (C$-$), $\mathfrak{U}_q(sl(2,\mathbf{R}))$ (C+). For other values of $q$ it is not possible to find star structure with $F = 1$. The study of real form of quantum groups has been first initiated in [19], and classified in the (C+) case for all $\mathfrak{U}_q(\mathfrak{G})$ in [42].
As already said before, star involutions on $\mathbb{C}$ Hopf algebras can perfectly exist without any $\mathbb{R}$ Hopf algebras associated to this star structure. Indeed if $A$ is a $\mathbb{C}$ Hopf algebra endowed with a star structure $\star$ then we can define $B = \{a \in A, a^\star = \theta(a)\}$ where $\theta$ is a linear involutive Hopf algebra automorphism of $A$, and study the case where $B$ is a $\mathbb{R}$ Hopf algebra such that $A = B^{\mathbb{C}}$. The automorphism $\theta$ has to satisfy $(\theta \circ \star)^2 = \text{id}$. If we choose the natural choice $\theta = S$, then it is easy to show that if $B$ is a real Hopf algebra such that $A = B^{\mathbb{C}}$ it then implies $S^2 = id$. This last condition severely restricts the structure of the Hopf algebra $A$ and excludes all the $\mathfrak{U}_q(\mathfrak{G})$. To our knowledge there is no general construction of $\theta$ which would ensure the existence of $B$ as a real Hopf algebra, so that in the rest of this article a real form of a $\mathbb{C}$ Hopf algebra is a definition of a star involution.

Let us end this paragraph with the following example, let $A = \mathfrak{U}_\hbar(sl(2))$, be the $\mathbb{C}[[\hbar]]$ Hopf algebra as defined by V.G.Drinfeld, it is the algebra defined by $H, J_+, J_-$ satisfying the standard commutation relations :

$$[H, J_\pm] = \pm 2 J_\pm, \qquad [J_+, J_-] = \frac{\sinh(\alpha\hbar H)}{\sinh(\alpha\hbar)} \tag{4}$$

and denote $q = \exp(\alpha\hbar)$. The coproduct is defined by:

$$\Delta(H) = H \otimes 1 + 1 \otimes H, \qquad \Delta(J_\pm) = q^{-\frac{H}{2}} \otimes J_\pm + J_\pm \otimes q^{\frac{H}{2}} \tag{5}$$

We can define a star structure: $\hbar^\star = \hbar, H^\star = H, J_+^\star = J_-, J_-^\star = J_+$ when $\alpha$ is real ($q$ real, (C+)) or $\alpha$ is imaginary ($|q| = 1$, (C$-$)). Let us denote $J_z = iH, J_x = (J_+ - J_-), J_y = i(J_+ + J_-)$, from the expression of the commutation relations, the $\mathbb{R}[[\hbar]]$ subalgebra $B$ of $A$ generated by $J_x, J_y, J_z$ is well defined and we have $A = B^{\mathbb{C}}$ as an algebra. From the expression of the coproduct it is easy to show that $B$ is a $\mathbb{R}[[\hbar]]$ Hopf algebra in the only case where $\alpha$ is pure imaginary. In both cases we can define a morphism of algebra $\theta$ on $A$ by $\theta(H) = -H, \theta(J_\pm) = -J_\pm$, and $B = \{a^\star = \theta(a)\}$, $\theta$ is always an anticomorphism of coalgebra, as a result $B$ is only a real Hopf algebra when $\star$ is of type (C$-$) i.e $\alpha$ pure imaginary.

Although $\mathfrak{U}_q(su(2))$ when $q$ real has been the major center of interest in the mathematical litterature, the other real form $\mathfrak{U}_q(su(2))$ with $|q| = 1$ has shown to be of central importance in connection with Conformal Field Theory and Chern-Simons theory [30][1].
Let $(A, m, \Delta)$ be a Hopf algebra, in the rest of this article we will use the following convention: $A_{op}$ denotes the Hopf algebra equal to $(A, m_{op}, \Delta)$ with $m_{op} = m \circ \sigma$ and $A^{op}$ the Hopf algebra equal to $(A, m, \Delta^{op})$ with $\Delta^{op} = \sigma \circ \Delta$, where $\sigma$ is the permutation operator. Let $A$ be a Hopf algebra, and $\overline{A}$ the Hopf algebra with coproduct $\overline{\Delta}$ defined by $\overline{\Delta}(\overline{a}) = \sum_{(a)} \overline{a_1} \otimes \overline{a_2}$. As soon as $A$ is endowed with a $\star$ structure of type $(C+)$ we can identify $A$ with $\bar{A}^{op}$ using for example the isomorphisms of Hopf algebra: $\rho^{\pm} : A \rightarrow \overline{A}^{op}, x \mapsto \overline{S^{\pm 1}(x^\star)}$.
If $\pi$ is a representation of $A$ acting on a vector space $V$, we can define the representation $\overline{\pi}$ of $A$ acting on $\overline{V}$ given by $\overline{\pi}(x) = \overline{(\pi \circ \rho^{\pm})(x)}$. If $V$ is an hermitian vector space, a representation of $A$ is said to be unitary if $\pi(a^\star) = \pi(a)^\dagger$. This notion is easily shown to generalize the standard notion of unitary representation of real Lie algebra.

Let $A^* = End_{\mathbb{C}}(A, \mathbb{C})$, if $A$ has a star structure of type (C$-$) (resp. (C+)) then $A^*$ can be endowed with star structures of type (C$-$) (resp. (C+)) using the following definitions:

$$\alpha^\star(a) = \overline{\alpha(a^\star)} \text{ if } \alpha, a \in A^\star \times A \text{ of type (C}-\text{)} \tag{6}$$
$$(resp. \qquad \alpha^\star(a) = \overline{\alpha(S^{\pm 1} a^\star)} \text{ if } \alpha, a \in A^\star \times A \text{ of type (C+)).} \tag{7}$$

In the rest of this article, we will always choose the definition with $S^{-1}$ in the case (C+) ( this in order to allow later a construction of star structure on the quantum double of $A$).
If $(A, R)$ is a quasitriangular Hopf algebra we will always denote $R^{(+)} = R$, $R^{(-)} = \sigma(R)^{-1}$ and $R' = R_{21} = \sigma(R)$. Let us consider a quasitriangular Hopf algebra with a real form of type (C+) or (C$-$), we automatically get a constraint on $R$ which is as follows:

$$[RR^{\star\otimes\star}, \Delta^{op}(a)] = 0 \qquad \forall a \in A \text{ in the case (C+)} \tag{8}$$
$$[R^{-1}R^{\star\otimes\star}, \Delta(a)] = 0 \qquad \forall a \in A \text{ in the case (C}-\text{).} \tag{9}$$

As a result it is consistent to study a smaller class of real form of quasitriangular Hopf algebras satisfying:

$$R^{\star\otimes\star} = R^{-1} \quad \text{or} \quad R^{\star\otimes\star} = R' \qquad \text{in the case (C+)} \tag{10}$$
$$R^{\star\otimes\star} = R'^{-1} \quad \text{or} \quad R^{\star\otimes\star} = R \qquad \text{in the case (C}-\text{).} \tag{11}$$

By convention the case where $R^{\star\otimes\star} \in \{R, R'\}$ is called the real case, and $R^{\star\otimes\star} \in \{R^{-1}, R'^{-1}\}$ the imaginary case ( the reason for this definition follows from the fact that in the $\mathfrak{U}_q(su(2))$ we have $R^{\star\otimes\star} = R'$ if $q$ is real and $R^{\star\otimes\star} = R^{-1}$ if $q$ is unimodular).
Let us recall the notion of twisted (or braided) tensor product of Hopf Algebras which will be extensively used in our work. Let $A$ and $B$ be two Hopf algebras and $F$ an element of $B \otimes A$, satisfying

$$(\Delta_B \otimes id)(F) = F_{23}F_{13}, \qquad (id \otimes \Delta_A)(F) = F_{12}F_{13} \tag{12}$$

the twisted Hopf algebra denoted $A \otimes_F B$ is the Hopf algebra equal to $A \otimes B$ as an algebra and with coproduct $\Delta$ defined by:

$$\Delta(a \otimes b) = F_{23}\Delta^A_{13}(a)\Delta^B_{24}(b)F_{23}^{-1}. \tag{13}$$

It then follows that $F$ satisfies

$$(id \otimes S_A)(F) = (S_B^{-1} \otimes id)(F) = F^{-1} \tag{14}$$

the antipode on $A \otimes_F B$ being given by the formula:

$$S^{\pm 1}(a \otimes b) = F_{21}^{-1}(S_A^{\pm 1}(a) \otimes S_B^{\pm 1}(b))F_{21}. \tag{15}$$

The quantum double of a Hopf algebra is a particular case of this construction. Let $A$ be a Hopf algebra of finite dimension, let $e_i$ be any choice of basis of $A$ and let us denote by $e^i$ the dual basis, we define $C$ to be the Hopf algebra $C = A^* \otimes_F A_{op}$ where $F^{-1} = \sum_i e_i \otimes e^i$. It is easy to show that $F = \sum_i e_i \otimes S_{A^*}^{-1} e^i$ satisfies (12). The quantum double is by definition the quasitriangular Hopf algebra $D(A) = C^*$ with $R$ matrix $R = \sum_i e_i \otimes 1 \otimes 1 \otimes e^i$. As a coalgebra it is given by $D(A) = A \otimes A^{*op}$, the algebra law being given by:

$$(a \otimes \xi).(b \otimes \eta) = \sum_{(b),(\xi)} ab_2 \otimes \xi_2 \eta < \xi_1, S^{-1} b_1 >< \xi_3, b_3 >, \tag{16}$$

where $a, b \in A\ ,\quad \xi, \eta \in A^*$, and we have used the Sweedler notation $\Delta(x) = \sum_{(x)} x_1 \otimes x_2$. We will use the notation $D(A) = A \hat{\otimes} A^{*op}$ to remind the reader that inside $D(A)$, $(A \otimes 1)$ and $(1 \otimes A^{*op})$ do not commute but satisfy the braided equalities (16).

The construction of the quantum double is due to V.G. Drinfeld [17], and its construction using a twisted tensor product appears first in [38], note however that the definition of Drinfeld gives $D(A) = A \otimes A^{*op}$ whereas Reshetikhin and Semenov one gives $D(A) = A^{*op} \otimes A$. In the infinite dimensional case, which is the case of interest for the study of deformation of Lie algebras, it is still possible to use this construction after having precisely defined the convergence of the infinite sum $F$. It is easy to show that if $A$ is endowed with a star structure $\star$ of type $C\pm$ this star structure can be extended uniquely to a star structure of type $C\pm$ denoted $\bigstar$ on $D(A)$ by the definitions:

$$\begin{aligned} \alpha^{\bigstar}(a) &= \overline{\alpha(a^\star)} \quad \text{if } \alpha, a \in A^{*op} \times A \ \text{ in the case (C}-) \\ \alpha^{\bigstar}(a) &= \overline{\alpha(S^{-1} a^\star))} \quad \text{if } \alpha, a \in A^{*op} \times A \ \text{ in the case (C+)}. \end{aligned} \tag{17}$$

Let $A$ be a quasitriangular Hopf algebra, let us write $R = \sum_i a_i \otimes b_i$, then the element $u = \sum_i S(b_i) a_i$ is invertible, its inverse is $u^{-1} = \sum_i S^{-2}(b_i) a_i$ and satisfies:

$$S^2(x) = uxu^{-1} \tag{18}$$

$$\Delta(u) = (u \otimes u)(R'R)^{-1} = (R'R)^{-1}(u \otimes u). \tag{19}$$

This quasitriangular Hopf algebra is a ribbon Hopf algebra, if and only if there exists an invertible element $v$ such that:

$$v \quad \text{is a central element} \tag{20}$$

$$v^2 = uS(u), \epsilon(v) = 1, S(v) = v \tag{21}$$

$$\Delta(uv^{-1}) = uv^{-1} \otimes uv^{-1} \tag{22}$$

$$\Delta(v) = (v \otimes v)(R'R)^{-1} = (R'R)^{-1}(v \otimes v). \tag{23}$$

We will denote $\mu$ the group like element $\mu = uv^{-1}$.

Let $A$ be a Hopf algebra, we define an action (resp. an antiaction) of $A$ on $A^*$ called $\overset{R}{\triangleright}$ (resp. $\overset{L}{\triangleright}$) by :

$$a \overset{L}{\triangleright} b = (id \otimes a)\Delta(b) \qquad a \overset{R}{\triangleright} b = (a \otimes id)\Delta(b) \quad \forall a \in A\ , b \in A^*. \tag{24}$$

It is easy to see that these operations satisfy:

$$aa' \overset{R}{\triangleright} b = a \overset{R}{\triangleright} (a' \overset{R}{\triangleright} b); \qquad aa' \overset{L}{\triangleright} b = a' \overset{L}{\triangleright} (a \overset{L}{\triangleright} b) \quad \forall a, a' \in A, \forall b \in A^* \tag{25}$$

$$a \overset{R}{\triangleright} (bc) = \sum_{(a)} a_1 \overset{R}{\triangleright} b \ \ a_2 \overset{R}{\triangleright} c\ ; \qquad a \overset{L}{\triangleright} (bc) = \sum_{(a)} a_1 \overset{L}{\triangleright} b \ \ a_2 \overset{L}{\triangleright} c \quad \forall a \in A, \forall b, c \in A^* \tag{26}$$

An action (or antiaction) satisfying the last property will be called derivation.

We can define two other actions, called adjoint actions, of $A$ on $A^*$ defined as follows:

$$a \overset{\mathrm{ad}^{\pm}}{\triangleright} b = \sum_{(a),(b)} < S^{\pm 1}(a_1), b_1 > b_2 < a_2, b_3 >, \quad \forall a \in A, \forall b \in A^*. \tag{27}$$

Note that, unless the coproduct on $A$ is cocommutative, these adjoint actions are not derivations. We will make a constant use of the following property: let $b$ be an element of $A^*$,

$$(a \overset{\mathrm{ad}^{-}}{\triangleright} b = \epsilon(a) b, \forall a \in A) \quad \Leftrightarrow \quad (\forall a \in A, a \overset{R}{\triangleright} b = a \overset{L}{\triangleright} b) \tag{28}$$

which is straightforward to verify. In this case, $b$ is said $ad^{\pm}$ invariant.

If $A$ is a ribbon Hopf algebra and $\pi$ is a finite dimensional representation of $A$, acting on $V$, the matrix elements of $\pi$ are elements of $A^*$ and the linear forms $tr_V(\pi(.))$ (resp. $tr_V(\pi(\mu^{-1})\pi(.))$ ) are respectively invariant elements of $A^*$ under the action $\overset{\mathrm{ad}^{-}}{\triangleright}$ (resp. $\overset{\mathrm{ad}^{+}}{\triangleright}$). These are very standard results in quantum group theory. Let $A$ be a Hopf algebra, a right (resp. left) invariant integral $h$ is an element $h_R$ (resp. $h_L$) of $A^*$ satisfying the following equation:

$$\forall a \in A, \sum_{(a)} a_1 h_R(a_2) = h_R(a)\, 1 \quad (\text{resp. } \sum_{(a)} h_L(a_1) a_2 = h_L(a)\, 1) \tag{29}$$

We will assume that $A$ is endowed with a star structure and non zero right and left invariant integrals, satisfying $\forall a \in A$, $h_L(a^\star) = \overline{h_L(a)}$, and $h_R(a^\star) = \overline{h_R(a)}$. We can define hermitian scalar product on $A^*$ as follows

$$(b|c)_R = h_R(b^\star c), \;\; (b|c)_L = h_L(b^\star c). \tag{30}$$

It is easy to show that the actions $\overset{R,L}{\triangleright}$ of $A$ on $A^*$ satisfy:

$$(a \overset{R}{\triangleright} b|c)_R = (b|a^\star \overset{R}{\triangleright} c)_R \qquad (S^{-2}(a) \overset{L}{\triangleright} b|c)_L = (b|a^\star \overset{L}{\triangleright} c)_L. \tag{31}$$

In particular the action $\overset{R}{\triangleright}$ is unitary.

In the case where there exists a non zero element $h$ of $A^*$ which is simultaneously a right integral and a left integral, we can define a hermitian product: $\forall b, c \in A^*, (b|c) = h(b^\star c)$. It is easy to show that the adjoint action $\overset{\mathrm{ad}^{-}}{\triangleright}$ of $A$ on $A^\star$ is unitary for this hermitian product. In the case where moreover $A$ is a ribbon Hopf algebra, we can define another hermitian product $((.|.))$ defined by $\forall b, c \in A^*, ((b|c)) = \sum_{(b),(c)} < \mu, b_1 >< \mu, c_1 > (b_2|c_2)$ for which the action $\overset{\mathrm{ad}^{+}}{\triangleright}$ is unitary.

## 2.2 Realification of Hopf Algebras

Let $A$ be a complex Hopf algebra, and let us denote $B$ the Hopf algebra $B = A \otimes_{\mathbb{C}} \overline{A}^{op}$ with the star structure defined by $(x \otimes \overline{y})^{\bigstar} = Sy \otimes \overline{S^{-1}x}$.

It is easy to show that $(B, \bigstar)$ is a star Hopf algebra of type $(C+)$ which will be called realification of $A$. If $A$ is quasitriangular and if we denote by $R$ its R-matrix, $\overline{A}^{op}$ is also quasitriangular with the natural R-matrix $\overline{R}'$, then $B$ has at least two quasitriangular structure:

$$\mathcal{R}^{(\pm)} = R_{13}^{(\pm)} \overline{R}_{42}^{(\pm)} \text{ of real type} \tag{32}$$

$$\mathcal{R}^{(\pm)} = R_{13}^{(\mp)} \overline{R}_{42}^{(\pm)} \text{ of imaginary type.} \tag{33}$$

It is still possible to twist this construction and we will use for this purpose the construction of twisted product. Let $Q$ be an invertible element of $\overline{A}^{op} \otimes A$ such that

$$(\Delta_{\overline{A}^{op}} \otimes id)(Q) = Q_{23}Q_{13}, \quad (id \otimes \Delta_A)(Q) = Q_{12}Q_{13}, \quad Q_{12} = \overline{Q_{21}}, \tag{34}$$

we will denote as usual $A \otimes_Q \overline{A}^{op}$ the Hopf algebra obtained from $A \otimes_{\mathbb{C}} \overline{A}^{op}$ by twisting. It is easy to show that this Hopf algebra can be endowed with a star structure of type $C+$ by the following formula :

$$(x \otimes \bar{y})^\star = Q_{21}^{-1}(S_A y \otimes \overline{S_A^{-1} x})Q_{21}. \tag{35}$$

In the case where $A$ is quasitriangular, the equations imply that

$$R_{23}Q_{12}Q_{13} = Q_{13}Q_{12}R_{23} \tag{36}$$

and similarly $A \otimes_Q \overline{A}^{op}$ has two quasitriangular structure (obtained by twisting):

$$\mathcal{R}^{(\pm)} = Q_{14}R_{13}^{(\pm)}\overline{R}_{42}^{(\pm)}Q_{23}^{-1} \text{ of real type} \tag{37}$$

$$\mathcal{R}^{(\pm)} = Q_{14}R_{13}^{(\mp)}\overline{R}_{42}^{(\pm)}Q_{23}^{-1} \text{ of imaginary type.} \tag{38}$$

This construction is general and completely independent of the existence of a star structure on $A$. In the particular case where $A$ is endowed with a star structure it has been developped by [13].

## 2.3 Quantum Double and Realification of Factorizable Hopf Algebras

Recall that if $A$ is a quasitriangular Hopf algebra, it is said factorizable if the linear map:

$$A^* \to A$$
$$\xi \mapsto (\xi \otimes id)(RR')$$

is an isomorphism of vector space.

When $A$ is a factorizable Hopf algebra, then the theorem of factorization of [38] applies and we have an isomorphism of Hopf algebra between the quantum double $D(A)$ and $A \otimes_{R^{-1}} A$:

$$\begin{aligned} \phi : \quad & D(A) \to A \otimes_{R^{-1}} A \\ & a \otimes \xi \mapsto \sum_{(a),(\xi)} a_1 \hat{R}^{(+)}(\xi_1) \otimes a_2 \hat{R}^{(-)}(\xi_2) \end{aligned} \tag{39}$$

where $a, \xi \in A \times A^{*op}$ and $\hat{R}^{(\pm)} : A^{*op} \to A$ are the algebra homomorphisms defined by $R^{(\pm)}(\xi) = (\xi \otimes id)(R^{\pm})$. It is then easy to compute the image of the universal $\mathcal{R}_D$ matrix of $D(A)$ under this isomorphism:

$$(\phi \otimes \phi)(\mathcal{R}_D^{(\pm)}) = R_{14}^{(-)}R_{24}^{(\mp)}R_{13}^{(\pm)}R_{23}^{(+)}. \tag{40}$$

In the case where $A$ is factorizable, endowed with a star structure of type $C+$, and if $R$ satisfies $R^{\star\otimes\star} = R'$ then $(id \otimes \rho^\pm)\phi$ is an isomorphism of star Hopf algebras between $D(A)$ and $A \otimes_{Q^\pm} \overline{A}^{op}$ with $Q^\pm = (\rho^\pm \otimes id)(R^{-1})$.
From this isomorphism we obtain that $D(A)$ is isomorphic to $A \otimes_{R^{-1}} A$ with star structure $(x \otimes y)^\star = R'(y^\star \otimes x^\star)R'^{-1}$.
This last theorem is very important because it shows that one realification of a star Hopf algebra

$A$ is the quantum double of $A$ [36][13].
On the otherhand, if $A$ is factorizable, endowed with a star structure of type $C+$, and if $R$ satisfies the other condition $R^{\star\otimes\star} = R^{-1}$ then $D(A)$ is isomorphic to $A \otimes_{\mathbb{C}} A$ with star structure defined by $(x \otimes y)^{\bigstar} = x^{\star} \otimes y^{\star}$.

Up to this point the reader can be puzzled by the fact that there exist more than one realification of a Hopf algebra depending on the choice of $Q$. In order to classify the possible realification we face two problems: what is the structure of the elements $Q$ satisfying equations (34) ? What are the possible structures of a realification, eventually generalizing or not the structure we gave by twisted product? In the case of $sl(2,\mathbb{C})_{\mathbb{R}}$ this problem has been studied in great details by Zakrewski and Woronowicz [52][51].

# 3 The Quantum Lorentz Group

## 3.1 Quantum Deformation of the Envelopping Algebra of the Lorentz Group

If $A$ is a star algebra, $\mathrm{Mat}_n(\mathbb{C}) \otimes A$ will always be endowed with the star structure defined by $(M^{\dagger})^i_j = (M^j_i)^{\star}$ where $M \in \mathrm{Mat}_n(\mathbb{C}) \otimes A$.

Let $W$ be a hermitian vector space, a representation $\pi : A \to End(W)$ is unitary if and only if $\pi(a^{\star}) = \pi(a)^{\dagger}, \forall a \in A$.

In this section we review the construction of the quantum deformation of the Lie algebra of the Lorentz group and precise our conventions.

We first define $A = \mathfrak{U}_q(su(2))$ for $q \in ]0,1[$, as being the star Hopf algebra defined by the generators $J_{\pm}, q^{\pm J_z}$, and the relations:

$$q^{\pm J_z} q^{\mp J_z} = 1 \qquad q^{J_z} J_{\pm} q^{-J_z} = q^{\pm 1} J_{\pm} \qquad [J_+, J_-] = \frac{q^{2J_z} - q^{-2J_z}}{q - q^{-1}}. \tag{41}$$

The coproduct is defined by

$$\Delta(q^{\pm J_z}) = q^{\pm J_z} \otimes q^{\pm J_z} \qquad \Delta(J_{\pm}) = q^{-J_z} \otimes J_{\pm} + J_{\pm} \otimes q^{J_z} \tag{42}$$

and the star structure is given by:

$$(q^{J_z})^{\star} = q^{J_z} \ , \qquad J_{\pm}^{\star} = q^{\mp 1} J_{\mp}. \tag{43}$$

This Hopf algebra is a ribbon quasi-triangular Hopf algebra. A description of the universal $R$ matrix, of the element $\mu, v$ is given in the Appendix.

Finite dimensional representations of $\mathfrak{U}_q(su(2))$ are completely reducible, and irreducible representations are completely classified by a couple $(\omega, K) \in \{1,-1,i,-i\} \times \frac{1}{2}\mathbb{Z}^+$, the dimension of the representation is $2K+1$ and the spectrum of $q^{J_z}$ is included in $\{\omega q^n, n \in \frac{1}{2}\mathbb{Z}\}$. The appearance of $\omega \in \{-1,i,-i\}$ comes from the existence of non trivial outer automorphisms $\tau_{\omega}$ of $\mathfrak{U}_q(su(2))$, defined by $\tau_{\omega}(q^{J_z}) = \omega q^{J_z}, \tau_{\omega}(J_+) = \omega^2 J_+, \tau_{\omega}(J_-) = J_-$. These automorphisms have no classical counterpart, and are a source of annoying features. In all the quantization schemes with specialization of the parameter to a complex number these automorphisms appear.

It is easy to show that an irreducible representation is unitarizable if and only if $\omega \in \{+1,-1\}$.

We will denote in the rest of this article $Irr(A)$ the set of all equivalency classes of finite dimensional irreducible and unitary representations with $\omega = 1$. The tensor product of elements of $Irr(A)$ is completely reducible in elements of $Irr(A)$.

Each of these classes $\overset{K}{\pi}$ is labelled by its spin $K$ (we will use a capital letter to denote it), and let us define $\overset{K}{V}$ as being the vector space associated to the representation of spin $K$. The representation of spin 0 is associated to the counit. The representation of dimension $2K+1$ and associated to $\omega$, $(\omega^4=1)$ is $\overset{K}{\pi}\circ\tau_\omega$.

Let $M$ be an element of $End(\overset{K}{V})$, we will define $\mathrm{tr}_q(M) = tr_{\overset{K}{V}}(\overset{K}{\mu}{}^{-1}M)$. In particular the q-dimension of the representation $\overset{K}{\pi}$ is by definition $[d_K]=\mathrm{tr}_q(1)=[2K+1]$.

The tensor product of elements of $Irr(A)$ behaves as in the classical case:

$$\overset{I}{\pi}\otimes\overset{J}{\pi} = \bigoplus_{K=|I-J|}^{I+J}\overset{K}{\pi}. \tag{44}$$

It will be convenient to introduce at this point the following notation: if $I,J,K$ are elements of $\frac{1}{2}\mathbb{Z}^+$, we define $Y(I,J,K)$ as follows:

$$Y(I,J,K)=\begin{cases} 1 & \text{if} \quad I+J-K, J+K-I, K+I-J\in\mathbb{Z}^+ \\ 0 & \text{otherwise}\end{cases}$$

We also define

$$\tilde{Y}(I,J)=\begin{cases} 1 & \text{if} \quad I+J, I-J\in\mathbb{Z}^+ \\ 0 & \text{otherwise}\end{cases}$$

Let us denote by $(\overset{K}{e}_i | i=1...dim\overset{K}{V})$ a particular unitary basis of $\overset{K}{V}$, and $(\overset{K}{e}{}^i | i=1...dim\overset{K}{V})$ its dual basis.

For any representations $\overset{I}{\pi},\overset{J}{\pi},\overset{K}{\pi}$ of $Irr(A)$ we define the Clebsch-Gordan maps $\psi^K_{IJ}$ (resp. $\phi^{IJ}_K$) as elements of $\mathrm{Hom}_A(\overset{I}{V}\otimes\overset{J}{V},\overset{K}{V})$ (resp. $\mathrm{Hom}_A(\overset{K}{V},\overset{I}{V}\otimes\overset{J}{V})$). These Clebsch-Gordan maps will be chosen to satisfy the following properties:

when $Y(I,J,K)=0$ we have $\psi^K_{IJ}=0, \phi^{IJ}_K=0$,

when $Y(I,J,K)=1$ $\psi^K_{IJ},\phi^{IJ}_K$ are non zero and defined by Clebsch-Gordan coefficients:

$$\phi^{IJ}_K(\overset{K}{e}_c)=\sum_{a,b}\left(\begin{matrix} a & b \\ I & J\end{matrix}\middle|\begin{matrix} K \\ c\end{matrix}\right)\overset{I}{e}_a\otimes\overset{J}{e}_b \qquad \psi^K_{IJ}(\overset{I}{e}_a\otimes\overset{J}{e}_b)=\sum_c\left(\begin{matrix} c \\ K\end{matrix}\middle|\begin{matrix} I & J \\ a & b\end{matrix}\right)\overset{K}{e}_c. \tag{45}$$

These Clebsch-Gordan maps are at this point defined up to phases. We will impose first that $\psi^K_{IJ}=(\phi^{IJ}_K)^\dagger$. These maps are still defined up to a phase which is fixed up to a sign by the requirement that the Clebsch-Gordan coefficients are real (this is possible because $q$ is real). These last signs are completely determined by the convention of Wigner : $\left(\begin{matrix} I & -J \\ I & J\end{matrix}\middle|\begin{matrix} K \\ I-J\end{matrix}\right)\in\mathbb{R}^{+*}$ .

With these conventions, the Clebsch-Gordan coefficients can be explicitely computed and satisfy numerous properties, some of them being stated in the Appendix.

To a representation $\pi$ of $\mathfrak{U}_q(su(2))$, we can naturally associate three representations, the conjugate representation $\overline{\pi}$, defined by $\overline{\pi}(x)=\overline{\pi(S^{-1}(x^\star))}$ and two contragredient representations $\hat{\pi}$ (resp. $\check{\pi}$) defined by $\hat{\pi}(x)={}^t\pi(S(x))$ (resp. $\check{\pi}={}^t\pi(S^{-1}(x))$ ) for any $x\in\mathfrak{U}_q(su(2))$.

In the special case of $\mathfrak{U}_q(su(2))$ and if we take $\pi=\overset{I}{\pi}$ then these three representations are equivalent to $\overset{I}{\pi}$, in particular we have $\overline{\pi}(x)=\overset{I}{W}\pi(x)\overset{I}{W}{}^{-1}$ with $\overset{I}{W}$ matrices which components are easily shown to be defined up to a scalar : $\overset{I}{W}{}^{\bar{m}}_{m'}=c_I\delta^m_{-m'}q^m(-1)^{I-m}$.

We will define the linear forms $\overset{I}{u}{}_i^j = < \overset{I}{e}{}^j | \overset{I}{\pi}(.) | \overset{I}{e}_i >$ and $Pol(SU_q(2))$ the star Hopf algebra generated as a vector space by $(\overset{I}{u}{}_i^j)_{I\in\frac{1}{2}\mathbb{Z}^+, i,j=1\dots 2I+1}$. By a direct application of the definitions, we have

$$\overset{I}{u}_1 \overset{J}{u}_2 = \sum_M \phi^{IJ}_M \, \overset{M}{u} \, \psi^M_{IJ}, \tag{46}$$

$$\Delta(\overset{I}{u}{}^a_b) = \sum_c \overset{I}{u}{}^a_c \otimes \overset{I}{u}{}^c_b. \tag{47}$$

The first equation implies the exchange relations

$$\overset{IJ}{R}_{12} \, \overset{I}{u}_1 \overset{J}{u}_2 = \overset{J}{u}_2 \overset{I}{u}_1 \, \overset{IJ}{R}_{12}, \tag{48}$$

where $\overset{IJ}{R} = (\overset{I}{\pi} \otimes \overset{J}{\pi})(R)$.

Due to the unitarity of the elements of $Irr(A)$, we have:

$$(\overset{I}{u}{}^a_b)^\star = S(\overset{I}{u}{}^b_a) \quad i.e. \quad \overset{I}{u}\overset{I}{u}{}^\dagger = \overset{I}{u}{}^\dagger \overset{I}{u} = 1. \tag{49}$$

The elements $\overset{I}{u}$ satisfy also the relation: $\overset{\bar{I}}{u} = \overset{I}{W} \overset{I}{u} \overset{I}{W}{}^{-1}$. If we define $\overset{I}{w}_{im} = \overset{I}{W}{}^{\bar{i}}_m$ then unitarity of the representation implies that:

$$S(\overset{I}{u}{}_i^j) = \sum_{m,n} \overset{I}{w}_{im} \overset{I}{u}{}^m_n \overset{I}{w}{}^{-1\,nj}. \tag{50}$$

Up to this point we will often use Einstein convention of summation over the indices labelling the vectors of a representation (the magnetic moments denoted by a small letter) and not on the indices labelling the representations (the spin of the representation denoted by a capital letter).

$Pol(SU_q(2))$ can equivalently be viewed as the star Hopf algebra generated by the matrix elements of $\overset{\frac{1}{2}}{u}$ with the defining relation (48) (with $I, J = \frac{1}{2}$) and the relation (49) (with $I = \frac{1}{2}$). Explicit relations are recalled in the Appendix.

$\mathfrak{U}_q(su(2))$ being a factorizable Hopf algebra it is possible to give a different but equivalent form of the defining relations of $\mathfrak{U}_q(su(2))$. Let us introduce, for each $I \in \frac{1}{2}\mathbb{Z}^+$ the elements $\overset{I}{L}{}^{(\pm)} \in \mathrm{End}(\overset{I}{V}) \otimes A$ defined by $\overset{I}{L}{}^{(\pm)} = (\overset{I}{\pi} \otimes id)(R^{(\pm)})$.

The duality bracket is given by

$$< \overset{I}{L}{}^{(\pm)}_1, \overset{J}{u}_2 >= \overset{IJ}{R}{}^{(\pm)}_{12} \tag{51}$$

These matrices satisfy the relations:

$$\overset{I}{L}{}^{(\pm)}_1 \overset{J}{L}{}^{(\pm)}_2 = \sum_K \phi^{IJ}_K \, \overset{K}{L}{}^{(\pm)} \, \psi^K_{IJ} \tag{52}$$

$$\overset{IJ}{R}_{12} \, \overset{I}{L}{}^{(+)}_1 \overset{J}{L}{}^{(-)}_2 = \overset{J}{L}{}^{(-)}_2 \overset{I}{L}{}^{(+)}_1 \, \overset{IJ}{R}_{12} \tag{53}$$

$$\Delta(\overset{I}{L}{}^{(\pm)a}{}_b) = \sum_c \overset{I}{L}{}^{(\pm)c}{}_b \otimes \overset{I}{L}{}^{(\pm)a}{}_c \tag{54}$$

$$(\overset{I}{L}{}^{(\pm)a}{}_b)^\star = S^{-1}(\overset{I}{L}{}^{(\mp)b}{}_a) \quad \text{i.e} \quad \overset{I}{L}{}^{(\pm)\dagger} \overset{I}{L}{}^{(\mp)} = \overset{I}{L}{}^{(\mp)} \overset{I}{L}{}^{(\pm)\dagger} = 1 \tag{55}$$

The first equation implies the exchange relations

$$\overset{IJ}{R}_{12}\overset{I}{L}{}_1^{(\pm)}\overset{J}{L}{}_2^{(\pm)} = \overset{J}{L}{}_2^{(\pm)}\overset{I}{L}{}_1^{(\pm)}\overset{IJ}{R}_{12} \tag{56}$$

Explicit relations between generators of $\mathfrak{U}_q(su(2))$ and matrix elements of $\overset{I}{L}{}^{(\pm)}$ are recalled in the Appendix.

The quantum double of $\mathfrak{U}_q(su(2))$ is the Hopf algebra $\mathcal{D} = \mathfrak{U}_q(su(2))\hat{\otimes}Pol(SU_q(2))^{op}$. We will denote by $\overset{I}{g}{}^a_b$ the elements $\overset{I}{u}{}^a_b$ embedded in $Pol(SU_q(2))^{op} \subset \mathcal{D}$. This is a star algebra with the following definition of star: $(\overset{I}{g}{}^a_b)^\star = S^{-1}(\overset{I}{g}{}^b_a)$ and also (55).

The algebra law satisfies the relations (46) (with $\overset{I}{u}$ replaced by $\overset{I}{g}$) and also the exchange relations:

$$\overset{IJ}{R}{}^{(\pm)}_{12}\overset{I}{L}{}^{(\pm)}_1\overset{J}{g}_2 = \overset{J}{g}_2\overset{I}{L}{}^{(\pm)}_1\overset{IJ}{R}{}^{(\pm)}_{12}. \tag{57}$$

The coproduct on $\mathcal{D}$ satisfies

$$\Delta_{\mathcal{D}}(\overset{I}{L}{}^{(\pm)a}{}_b) = \sum_c \overset{I}{L}{}^{(\pm)c}{}_b \otimes \overset{I}{L}{}^{(\pm)a}{}_c \quad , \quad \Delta_{\mathcal{D}}(\overset{I}{g}{}^a_b) = \sum_c \overset{I}{g}{}^c_b \otimes \overset{I}{g}{}^a_c. \tag{58}$$

$\mathcal{D}$ is a quantum deformation of the envelopping algebra of $sl(2,\mathbb{C})_{\mathbb{R}}$ [36][34][16]. In order to understand this point it is important to recall the principle of quantum duality which appears first in [17] and developed in [39]. Let $\star$ be the involution on $sl(2,\mathbb{C})$ selecting the compact form $su(2)$. $\mathfrak{U}_q(su(2))$ being a star Hopf algebra, $(sl(2,\mathbb{C}),\star)$ inherits a star Lie bialgebra structure. By duality it gives on $sl(2,\mathbb{C})^*$ a structure of star Lie bialgebra which real form is isomorphic as a $\mathbb{R}$ Lie algebra to the Lie algebra $an(2)$, the real Lie algebra of $2\times 2$ lower triangular complex matrices with real diagonal of zero trace. At the quantum level, an easy application of the quantum duality principle, shows that $Pol(SU_q(2)) = \mathfrak{U}_q(su(2))^* = \mathfrak{U}_q(su(2)^*) = \mathfrak{U}_q(an(2))$ (at this point we are quite unprecise on the exact signification of the dual, this will be corrected later on). As a result $\mathcal{D} = \mathfrak{U}_q(su(2)) \otimes \mathfrak{U}_q(an(2))$ as a vector space, containing $\mathfrak{U}_q(su(2))$ and $\mathfrak{U}_q(an(2))$ as subalgebras. It can be shown that the classical limit of $\mathcal{D}$ is the Lie algebra $sl(2,\mathbb{C})_{\mathbb{R}}$, and the classical limit of the previous decomposition is the Iwasawa decomposition. In the appendix the reader will find the complete defining set of commutation relations of $\mathcal{D}$. We will use the notation $\mathcal{D} = \mathfrak{U}_q(sl(2,\mathbb{C})_{\mathbb{R}})$.

Another natural decomposition, in the classical case, of $\mathfrak{U}(sl(2,\mathbb{C})_{\mathbb{R}})\otimes\mathbb{C}$ as a complex algebra is the decomposition in two commuting copies $\mathfrak{U}(sl(2,\mathbb{C}))\otimes\mathfrak{U}(sl(2,\mathbb{C}))$. In order to obtain such a description in the quantum case, let us now describe the factorization theorem applied to $\mathcal{D}$. Using the explicit form of the isomorphism $\phi$ and the expressions of the coproduct (58), we can write:

$$\phi(\overset{I}{L}{}^{(\pm)}) = \overset{I}{M}{}^{(r\pm)}\overset{I}{M}{}^{(l\pm)} \tag{59}$$

$$\phi(\overset{I}{g}) = \overset{I}{M}{}^{(r-)}\overset{I}{M}{}^{(l+)} \tag{60}$$

where the components of $\overset{I}{M}{}^{(l\pm)}$ (resp. $\overset{I}{M}{}^{(r\pm)}$) are located in the subalgebra $A\otimes 1$ (resp. $1\otimes A$) of $A\otimes A$.

From a direct application of the factorization theorem, we have the following relations:

$$\overset{IJ}{R}_{12}\overset{I}{M}{}_1^{(i\pm)}\overset{J}{M}{}_2^{(i\pm)} = \overset{J}{M}{}_2^{(i\pm)}\overset{I}{M}{}_1^{(i\pm)}\overset{IJ}{R}_{12} \tag{61}$$

$$\overset{IJ}{R}{}_{12}^{(\pm)}\overset{I}{M}{}_1^{(i\pm)}\overset{J}{M}{}_2^{(i\mp)} = \overset{J}{M}{}_2^{(i\mp)}\overset{I}{M}{}_1^{(i\pm)}\overset{IJ}{R}{}_{12}^{(\pm)} \tag{62}$$

$$\overset{I}{M}{}_1^{(l\sigma)}\overset{J}{M}{}_2^{(r\sigma')} = \overset{J}{M}{}_2^{(r\sigma')}\overset{I}{M}{}_1^{(l\sigma)} \tag{63}$$

with $i \in \{r, l\}, \sigma, \sigma' \in \{+, -\}$.

This factorization theorem gives an immediate description of the finite dimensional representations of $\mathcal{D}$. Indeed from the isomorphism of complex algebra between $\mathcal{D}$ and $\mathfrak{U}_q(sl(2,\mathbb{C})) \otimes \mathfrak{U}_q(sl(2,\mathbb{C}))$ we obtain that each finite dimensional representation of $\mathcal{D}$ is completely reducible and that the finite dimensional irreducible representations of $\mathcal{D}$ are labelled by a pair $(\pi, \pi')$ of irreducible finite dimensional representations of $\mathfrak{U}_q(sl(2,\mathbb{C}))$. If $\pi, \pi'$ are two irreducible finite dimensional representations acting on $V$ and $V'$, $\pi \otimes \pi'$ is an irreducible representation of $\mathfrak{U}_q(sl(2,\mathbb{C})) \otimes \mathfrak{U}_q(sl(2,\mathbb{C}))$, from which we deduce the associated representation $\rho(\pi, \pi')$ of $\mathcal{D}$ acting on $V \otimes V'$ by the expressions:

$$\rho(\pi,\pi')(\overset{I}{g}_1) = \overset{I\pi'}{R}{}_{13}^{(-)}\overset{I\pi}{R}{}_{12}^{(+)}, \qquad \rho(\pi,\pi')(\overset{I}{L}{}_1^{(\pm)}) = \overset{I\pi'}{R}{}_{13}^{(\pm)}\overset{I\pi}{R}{}_{12}^{(\pm)}. \tag{64}$$

At this point the reader is invited to look at the appendix where the commutation relations of generators of $\mathcal{D}$ are explicitely computed. From these relations it is easy to show that there exists a non trivial outer automorphism of $\mathcal{D}$ defined by $\tau_s(a) = -a, \tau_s(b) = -b, \tau_s(c) = -c, \tau_s(d) = -d, \tau_s(q^{J_z}) = q^{J_z}, \tau_s(J_\pm) = J_\pm$. This implies that $\epsilon_s = \epsilon \circ \tau_s$ is a one dimensional representation, which can be written as $\epsilon_s = \rho(\overset{0}{\pi}_i, \overset{0}{\pi}_{-i})$ and has therefore no classical analogue. This representation is called strange [36] (is not a smooth representation in the sense of these authors) and will reappear when we will study infinite dimensional representations of $\mathcal{D}$. It has been shown in [41] that finite dimensional representations of $\mathfrak{U}_q(sl(2,\mathbb{C})_{\mathbb{R}})$ are completely reducible and that the only finite dimensional irreducible representations, up to equivalence, of $\mathcal{D}$ are $\rho(\overset{I}{\pi}, \overset{J}{\pi})$ and $\rho(\overset{I}{\pi}, \overset{J}{\pi}) \circ \tau_s$.

The center of the algebra $\mathcal{D}$ is a polynomial algebra with two generators $\Omega^\pm$ defined by:

$$\Omega^\pm = \mathrm{tr}_q(\overset{\frac{1}{2}}{L}{}^{(\mp)\,-1}\overset{\frac{1}{2}}{g}). \tag{65}$$

An easy computation shows that

$$\phi(\Omega^+) = \mathrm{tr}_q(\overset{\frac{1}{2}}{M}{}^{(l-)\,-1}\overset{\frac{1}{2}}{M}{}^{(l+)}) \tag{66}$$

$$\phi(\Omega^-) = \mathrm{tr}_q(\overset{\frac{1}{2}}{M}{}^{(r+)\,-1}\overset{\frac{1}{2}}{M}{}^{(r-)}) \tag{67}$$

which are the quadratic Casimirs of each copies in the intrinsic form found by [19].

The action of the star involution is easy to compute and we have $(\Omega^+)^\star = \Omega^-$. The action of $\tau_s$ on $\Omega^\pm$ is given by: $\tau_s(\Omega^\pm) = -\Omega^\pm$. In the appendix we have worked out the precise connection betwen our conventions and those of [36][35].

## 3.2 Algebras of functions on $SL_q(2,\mathbb{C})$

### 3.2.1 Introduction

This chapter is intended to be a pedagogical introduction to the theory of multiplier Hopf algebras as first described and studied by S.L.Woronowicz [50]. It is not the aim of this chapter to give a detailed account of this theory but we will give a simple exposition of the main ideas and of the tools used in the rest of this work. The reader is urged to read the papers of Woronowicz [50], Podles and Woronowicz[36] and Van Daele [43, 44] where he will find a complete exposition of multiplier Hopf algebras as well as detailed proofs of the theorems. The aim of the theory of multipliers, applied to the quantum group case, is to develop an appropriate mathematical framework where quantization of non compact group naturally fits in. As we will see the formalism is not completely obvious, and it is the aim of this work to show that it is possible to develop a complete harmonic analysis of $SL_q(2,\mathbb{C})$ in this framework.

Let us first explain the philosophy of multipliers. Let $K$ be a compact group, and let $Pol(K)$ be the polynomial algebra of matrix coefficients of finite dimensional representations of $K$. $K$ being compact, we can endow $Pol(K)$ with a structure of normed algebra with a star structure and unit element:

$$\forall f \in Pol(K), \|f\| = \sup_{x\in K} |f(x)| \quad f^\star(x) = \overline{f(x^{-1})}. \tag{68}$$

The completion of this normed algebra gives the commutative $C^\star$ algebra $(C(K), \| \ \|)$ with $C(K)$ the space of complex continuous functions on $K$.

In the quantum case we can proceed along the same lines using the theory of compact matrix pseudogroup of Wonorowicz [48]. We will explain it for example in the case of $SU_q(2)$, see also for this example[45]. Let $Pol(SU_q(2))$ be the star Hopf algebra generated by the matrix elements of $Irr(\mathfrak{U}_q(su(2)))$. We can define on $Pol(SU_q(2))$ a norm defined by:

$$\|a\| = \sup_\pi \|\pi(a)\| \tag{69}$$

where the supremum is taken over all unitary representations of $Pol(SU_q(2))$ acting on a Hilbert space of finite or infinite dimension (in the case of $Pol(SU_q(2))$ these unitary representations have been classified in [45]). Of course we have to show that this definition is a norm. This is clearly the case because $Pol(SU_q(2))$ is generated by the coefficients of $\overset{\frac{1}{2}}{g} = \begin{bmatrix} a & b \\ -qb^\star & a^\star \end{bmatrix}$ satisfying, in particular, the relation $aa^\star + bb^\star = 1$. For each unitary representation $\pi$ of $Pol(SU_q(2))$, the last relation implies $\|\pi(a)\| \leq 1, \|\pi(b)\| \leq 1$. It can be shown that $\|x\| = 0 \Rightarrow x = 0$, see [45, 12]. After completion of $(Pol(SU_q(2)), \| \ \|)$ we obtain a $C^\star$ algebra denoted $(Fun(SU_q(2)), \| \ \|)$ which can be thought as being the quantum analog of the space of continuous functions on $SU_q(2)$.

The harmonic analysis on $SU_q(2)$, and more generally on each compact matrix pseudo group has been studied in the work [48]. The axioms of matrix pseudo groups are sufficiently general to include the category of compact group as well as all known examples of functions on compact form of quantum groups. It can be shown from these axioms that there exists a unique left and right integral and an analog of the Peter-Weyl theorem can be proved from first principles.

The coproduct $\Delta$ maps $Pol(SU_q(2))$ to $Pol(SU_q(2)) \otimes_{alg} Pol(SU_q(2))$, where we have denoted by $\otimes_{alg}$ the algebraic tensor product. After completion $\Delta$ is uniquely extended to $\underline{\Delta} : Fun(SU_q(2)) \to Fun(SU_q(2)) \otimes Fun(SU_q(2))$ where $\otimes$ denotes a suitable tensor product of $C^\star$ algebra (note that there is no ambiguity in the definition of the tensor product $Fun(SU_q(2)) \otimes Fun(SU_q(2))$ because $Fun(SU_q(2))$ is nuclear [49]. )

Let now $G$ be a locally compact but not compact group, we can still define $Pol(G)$ the polynomial algebra of matrix coefficients of finite dimensional representation of $G$. The only but crucial difference with the compact case is that $Pol(G)$ is no more a normed algebra with the definition of the norm (68), the reason being that the supremum is generically equal to infinity. The notion of convergence on $Pol(G)$ using the norm convergence is therefore useless. It is still possible to define weaker convergence using semi norms such as the convergence on all compact subspaces or the pointwise convergence. Unfortunately, because the notion of point looses its meaning in the quantum case, the generalization of weaker convergences cannot be directly generalized. So, contrary to the compact case, we cannot construct the algebra of continuous functions from a completion of the polynomial algebra. This is the major obstacle to define noncompact quantum group, the algebra of polynomials on the quantum group is not sufficient to reconstruct the space of continuous functions.

There are situations where we can go round this major hurdle, which is handled using the concept of multiplier. The philosophy is to revert the basic objects: the given space is a $C^\star$ algebra which can be thought of as being a quantization of the space of continuous functions vanishing at infinity. From this algebra, we can derive from pure algebraic constructions the set of multipliers of $A$ which plays the role of the algebra of bounded functions and the set of affiliated elements of $A$ which plays the role of continuous functions. The polynomial algebra is a derived set consisting of particular affiliated elements. At this point we have to stress that if we are given a quantization of the polynomial algebra of some non compact group there is no theorem which garanty that there exists a $C^\star$ algebra $A$ such that this polynomial algebra would be a certain set of affiliated elements of $A$. Woronowicz has even shown that $SU_q(1,1)$ (with $q \in \mathbb{R}$ can be defined as a star Hopf algebra but there exists no $C^\star$ algebra for which the elements of $SU_q(1,1)$ are affiliated [50].

### 3.2.2 Multipliers: the algebraic theory

This section is taken from the article of Van Daele [43], its aim is to recall notions that are not so well known and to recall definitions and notations that will be used all throughout our work.

Let $A$ be an associative algebra with or without a unit element (we will see that the construction of multipliers is really interesting in the case where $A$ admits no unit element). An element $\rho \in \mathrm{End(A)}$ is a left (resp. right) multiplier if

$$\forall a, b \in A \quad \rho(ab) = \rho(a)b \quad (\text{resp.} \quad \rho(ab) = a\rho(b)). \tag{70}$$

A multiplier $\rho$ is a couple $\rho = (\rho_1, \rho_2)$ with $\rho_1$ left multiplier, $\rho_2$ right multiplier satisfying:

$$\forall a, b \in A \quad a\rho_1(b) = \rho_2(a)b. \tag{71}$$

We will denote by $R(A)$ (resp. $L(A)$) the set of right (resp. left) multipliers, and by $M(A)$ the set of multipliers. $R(A), L(A)$ are subalgebras of End(A), $M(A)$ is also a subalgebra of $End(A) \times End(A)$.

To each element of $A$ we can associate a left multiplier $\rho_l(a)$, a right multiplier $\rho_r(a)$ and a multiplier $\rho(a)$ defined by:

$$\rho_l(a)(b) = ab \quad , \quad \rho_r(a)(b) = ba \quad \text{and} \quad \rho(a) = (\rho_l(a), \rho_r(a)) \tag{72}$$

It is easy to show that if the algebra law on $A$ is non degenerate, i.e $\forall b \in A \;\; ab = 0 \Rightarrow a = 0$ and $\forall b \in A \;\; ba = 0 \Rightarrow a = 0$ then the maps $A \xrightarrow{\rho_l} L(A), A \xrightarrow{\rho_r} R(A), A \xrightarrow{\rho} M(A)$ are monomorphisms of algebras.
If $A$ admits a unit element then all these maps are trivially isomorphisms. The notion of multiplier is therefore interesting only in the case where $A$ does not have a unit element.

In the rest of this chapter we will always assume that the multiplication law on $A$ is non degenerate. From the last proposition it will be convenient to identify $A$ as a subspace of $L(A)$ and of $R(A)$.
We will use a nice convention which has to be carefully used: if $\rho_1 \in L(A)$ and $\rho_2 \in R(A)$ we will use the notation $\rho_1 b = \rho_1(b)$ and $b\rho_2 = \rho_2(b)$.

When the algebra law is non degenerate the mappings $M(A) \to L(A), (\rho_1, \rho_2) \mapsto \rho_1$ and $M(A) \to R(A), (\rho_1, \rho_2) \mapsto \rho_2$ are easily seen to be injections, so it will be convenient to identify $M(A)$ with a subspace of $L(A)$ or with a subspace of $R(A)$ (this should create no confusion). The equation (71) is rewritten $(a\rho)b = a(\rho b)$ using all these identifications.

Let us give two examples which should illuminate the construction.

Example 1: A commutative example.
Let $X$ be a set and $A = \mathcal{F}_f(X)$ the set of functions on $X$ which are equal to zero except on a finite subset of $X$. If $X$ is infinite, $A$ does not admit a unit element.
It is easy to show that $L(A) = R(A) = M(A) = \mathcal{F}(X)$ where $\mathcal{F}(X)$ is the algebra of functions on $X$.

Example 2: A non commutative example.
Let $X$ be a set and $n$ a map from $X$ to $\mathbb{N}$. We define $A_x = \mathrm{Mat}_{n(x)}(\mathbb{C})$ and $A = \bigoplus_{x\in X} A_x$. $A$ is an algebra with unity if and only if $X$ is finite. When $n$ is the constant map $n(x) = 1$ we recover the example 1.

Let $(a_x)_{x\in X} \in \prod_{x\in X} A_x$, we can define a left multiplier $\rho'_l((a_x)_{x\in X}) \in L(A)$ by

$$\rho'_l((a_x)_{x\in X})(\oplus_{x\in X} b_x) = \oplus_{x\in X} a_x b_x. \tag{73}$$

It is easy to show that $\rho'_l$ realizes an isomorphism of algebra between $L(A)$ and $\prod_{x\in X} A_x$. Similarly we have $R(A) = L(A) = M(A) = \prod_{x\in X} A_x$.

This example although quite trivial is the core of the construction of multipliers applied to the functions on $SL_q(2, \mathbb{C})$.

If $A$ is a star algebra, we can define antilinear antimorphisms of algebra $\dagger : L(A) \to R(A)$ (resp $\dagger : R(A) \to L(A)$) using the definitions:

$$a\rho^\dagger = (\rho a^\star)^\star, \forall a \in A\ , \forall \rho \in L(A). \tag{74}$$

$$\rho^\dagger a = (a^\star \rho)^\star, \forall a \in A\ , \forall \rho \in R(A) \tag{75}$$

This definition is chosen such that $\dagger$ is equal to $\star$ on elements of $A$.

On $M(A)$ we can define a star structure by $(\rho_1, \rho_2)^\dagger = (\rho_2^\dagger, \rho_1^\dagger)$.

### 3.2.3 Multiplier and affiliated elements: the $C^\star$ algebra case

In this section we recall the definitions and the notations of the work of Wonorowicz [50]. Let $A$ be a $C^\star$ algebra with or without a unit element. We define $M_b(A)$ the space of bounded multipliers as the subspace of $M(A)$ consisting of elements $\rho = (\rho_1, \rho_2)$ such that $\rho_1$ and $\rho_2$ are bounded. It will be convenient to identify $M(A)$ as a closed subalgebra of the bounded linear operator acting on $A$. $M(A)$ is a $C^\star$ algebra with unity. Let us recall the heuristic definition of an affiliated element to $A$. Let $T$ be a linear operator $T : A \to A$ with domain $D(T)$ dense in $A$, $T$ is an element affiliated to $A$, and we write $T\eta A$ if bounded fonctions of $T$ belongs to $M_b(A)$. Woronowicz chooses the function $z \mapsto z(1 + z\bar{z})^{-\frac{1}{2}}$ which leads him to the precise definition: let $T$ be a linear operator $T : A \to A$ with domain $D(T)$ dense in $A$, $T$ is affiliated to $A$ if and only if there exists $Z \in M_b(A)$ such that $\|Z\| \leq 1$ and

$$(x \in D(T), y = Tx) \Leftrightarrow (\exists a \in A, x = (1 - Z^\dagger Z)^{\frac{1}{2}} a, y = Za). \tag{76}$$

Let us give two important examples which are the exact counterparts of the examples 1 and 2 given in the previous section.

Example 1.

Let $X$ be a locally compact space and let $A = Fun_0(X)$ be the $C^\star$ algebra of continuous functions which vanishes at infinity with norm $\|f\| = \sup_{x\in X}|f(x)|$. The set of bounded multipliers of $A$ is identified to the set of continuous bounded functions on $X$, whereas the affiliated elements are the set of complex continuous functions on $X$.
In this example, we can define the star algebra $Fun_c(X)$ of compact supported functions, this is not a Banach algebra when $X$ is non compact, and its completion with the previous norm gives the $C^\star$ algebra $A$.

Example 2.

Let $\{A_n, n \in \mathbb{N}\}$ be a family of $C^\star$ algebras with units, we define $A = \sum^{\oplus} A_n$ to be the $C^\star$ algebra of sequences $(a_n)_{n\in\mathbb{N}}$ with $\lim_{n\to+\infty}\|a_n\| = 0$, with the norm $\|(a_n)_{n\in\mathbb{N}}\| = \sup_n \|a_n\|$. The bounded multiplier are identified with the set of sequences $(a_n)_{n\in\mathbb{N}}$ such that $\sup_n\|a_n\| < +\infty$ and the set of affiliated elements of $A$ is identified with the star algebra $\prod_{n\in\mathbb{N}} A_n$. Let $A_c$ be the star algebra $A_c = \bigoplus_{n\in\mathbb{N}} A_n$, this is not a Banach algebra, but its completion using the previous norm gives the $C^\star$ algebra $A$.

Up to this point we have not introduced any structure of Hopf algebra on $A$, this is handled using the concept of multiplier Hopf algebra.

### 3.2.4 Multiplier Hopf algebras: the algebraic theory

In this section we give the main motivation of the concept of multiplier Hopf algebra and we recall the basic ingredients of this theory.

The little following example [43] should illuminate the necessity of working with multiplier Hopf algebras.

Let $G$ be a group, and let $A = \mathcal{F}_f(G)$ be the space of functions which vanish on $G$ except on a finite subset of $G$. $G$ being a group, we would expect that by duality $\mathcal{F}_f(G)$ is endowed with a structure of Hopf algebra. This is not the case: we can define $\Delta : \mathcal{F}(G) \to \mathcal{F}(G\times G)$, therefore $\Delta$ maps $\mathcal{F}_f(G)$ to $\mathcal{F}(G\times G)$ which is not, when $G$ is infinite, included in $\mathcal{F}_f(G)\otimes\mathcal{F}_f(G)$. However we have $\mathcal{F}(G\times G) = M(\mathcal{F}_f(G)\otimes\mathcal{F}_f(G))$, so that we are naturally led to the situation where $\Delta : A \to M(A\otimes A)$.

However from the expression of $\Delta(f)$ it is easy to show that we still have that $\Delta(f)(1\otimes g)$ and $\Delta(f)(g\otimes 1)$ are elements of $A\otimes A$ with $f, g\in A$.

The notion of multiplier Hopf algebra is the abstraction of this situation. We will not give an exposition of this theory which is very well exposed in [43]. All the axioms of Hopf algebra can be extended to this framework.

Let us give a simple example which has been coined discrete quantum group [44], and which will be of upmost importance in the case of $SL_q(2,\mathbb{C})_{\mathbb{R}}$. A discrete quantum group is a multiplier star Hopf algebra $(A,\Delta)$ where $A$ is a direct sum of matrix algebras. The notations are those of Example 2 of the previous subsections. Let $A = \bigoplus_{x\in X} A_x$ where $A_x = \mathrm{Mat}_{n(x)}(\mathbb{C})$, and let us define $\Delta : A \to M(A\otimes A) = \prod_{x,y} A_x\otimes A_y$, we will assume that the axioms of a multiplier Hopf algebra are satisfied. Let $e_x$ be the unit of $A_x$, it is a central projection in $A$ and we have $Ae_x = A_x$. From the fact that $\Delta(a)(1\otimes b)$ and $\Delta(a)(b\otimes 1)$ are in $A\otimes A$, we easily obtain the proposition 2.2 of [44]: consider the property $\Delta(e_z)(e_x\otimes e_y)\neq 0, x, y, z\in X$ ; given $x, y$ there are only finitely many $z$ such that this property holds, given $y, z$ there are only finitely many $x$ such that this property holds.
Van Daele defines a notion of right invariant integral, a left invariant integral appropriate to the framework of multiplier Hopf algebra, and show that for discrete quantum groups there exists a unique left integral and a unique right integral.

### 3.2.5 Multiplier Hopf algebras: the $C^\star$ algebra theory

This subject will not be developped here, it is a rapidly evolving subject which basic axioms are not completely fixed up to now, see the references [50], [36], [29]. The motivations of these works is to give axioms of non compact quantum groups, sufficiently general to include all the known examples and to prove in this theory an existence theorem for a left integral, and to begin a large program of harmonic analysis.

In the rest of our work we will only use the results obtained in [36].

In our work we will make an extensive use of the notion of multiplier Hopf algebra. A very good reference on this subject is [43, 44] in the pure algebraic case. We will also use the notions of $C^*$ multiplier Hopf algebra, and affiliated elements which is well exposed in [50, 29] and also [36] in the case of $SL_q(2,\mathbb{C})$.

The aim of the theory of multipliers, applied to the quantum group case, is to develop an appropriate mathematical framework where quantization of non compact group naturally fits in. As we will see the formalism is not completely obvious, and it is the aim of our work to show that it is possible to develop a complete harmonic analysis of $SL_q(2,\mathbb{C})$ in this framework.

Let us review briefly the basic ideas. $Pol(SU_q(2))$ is a star Hopf algebra and we can define on $Pol(SU_q(2))$ a norm defined by:

$$||a|| = \sup_\pi ||\pi(a)||, \ \forall a \in Pol(SU_q(2)) \tag{77}$$

where the supremum is taken over all unitary representations of $Pol(SU_q(2))$ acting on a Hilbert space (note that $\pi(a)$ is a bounded operator and that the supremum always exists and is finite [45]). After completion of $(Pol(SU_q(2))$ with respect to this norm we obtain a $C^\star$ algebra denoted $(Fun(SU_q(2)), ||\ \ ||)$ which can be thought as being the quantum analog of the space of continuous functions on $SU(2)$. The coproduct $\Delta : Pol(SU_q(2)) \rightarrow Pol(SU_q(2)) \otimes Pol(SU_q(2))$ extends by density to a morphism of $C^*$ algebras: $\Delta : Fun(SU_q(2)) \rightarrow Fun(SU_q(2)) \otimes Fun(SU_q(2))$.

From the quantum duality principle, it is natural to view the linear forms on $Pol(SU_q(2))$ as quantum analogues of classes of functions on $AN$.

More precisely Podles and Woronowicz were led to introduce the following construction: let $(\overset{I}{E}{}^i_j)$ be the elements of $(Pol(SU_q(2)))^*_{op}$ defined by the duality bracket: $< \overset{I}{E}{}^i_j, \overset{J}{g}{}^m_n >= \delta_{I,J}\delta^i_n\delta^m_j$.

It is easy to see that these elements are independent and that they are naturally endowed with the following structure of star algebra, dual to the coproduct on the elements $\overset{I}{g}{}^i_j$:

$$\overset{I}{E}{}^k_l \overset{J}{E}{}^r_s = \delta_{I,J}\overset{I}{E}{}^k_s\delta^r_l, \quad (\overset{I}{E}{}^i_j)^\star = \overset{I}{E}{}^j_i \tag{78}$$

As a result the structure of the algebra generated by $\overset{I}{E}{}^k_l$ is isomorphic to the star algebra $\bigoplus_{I\in\frac{1}{2}\mathbb{Z}^+} \mathrm{Mat}_{2\mathrm{I}+1}(\mathbb{C})$.

We define $Fun_c(AN_q(2))$ (the quantum analog of compact supported functions) to be the star algebra $\bigoplus_I \mathrm{Mat}_{2\mathrm{I}+1}(\mathbb{C})$ which basis is $\overset{I}{E}{}^k_l$. Note that this algebra has no unit element. This algebra is endowed with a structure of normed algebra defined as follows: $||\oplus_I a_I|| = sup_I ||a_I||$ where $||a_I||$ is the usual sup norm of finite dimensional matrix.

$Fun_0(AN_q(2))$ (the quantum analog of functions vanishing at infinity) is the $C^\star$ algebra $Fun_0(AN_q(2)) = \sum_I{}^\oplus \mathrm{Mat}_{2\mathrm{I}+1}(\mathbb{C})$ completion of $Fun_c(AN_q(2))$ with respect to this norm.

By duality $Fun_c(AN_q(2))$ is endowed with a structure of multiplier Hopf algebra, called discrete quantum group [44]. The coproduct $\Delta$ maps $Fun_c(AN_q(2))$ to $M(Fun_c(AN_q(2) \otimes Fun_c(AN_q(2))$ where $M(B)$ denotes the algebra of multiplier of $B$. In this particular case, $M(Fun_c(AN_q(2)\otimes Fun_c(AN_q(2))) = \prod_{I,J\in\frac{1}{2}\mathbb{Z}^+}\mathrm{Mat}_{2\mathrm{I}+1}(\mathbb{C})\otimes\mathrm{Mat}_{2\mathrm{J}+1}(\mathbb{C})$. It is easy to compute the coproduct by duality:

$$\Delta(\overset{I}{E}{}^i_j) = \sum_{J,K}{}^{\oplus}\begin{pmatrix} n & s & \Big| & I \\ J & K & \Big| & j\end{pmatrix}\begin{pmatrix} i & \Big| & J & K \\ I & \Big| & m & r\end{pmatrix}\overset{J}{E}{}^m_n\otimes\overset{K}{E}{}^r_s \tag{79}$$

where this infinite sum as to be understood as the multiplier

$$(\begin{pmatrix} n & s & \Big| & I \\ J & K & \Big| & j\end{pmatrix}\begin{pmatrix} i & \Big| & J & K \\ I & \Big| & m & r\end{pmatrix}\overset{J}{E}{}^m_n\otimes\overset{K}{E}{}^r_s\;)_{J,K\in\frac{1}{2}\mathbb{Z}^+}. \tag{80}$$

In order not to obscure formulas we will stick to the notation with the sum. It is therefore an easy task to check the axioms of a star multiplier Hopf algebra on $Fun_c(AN_q(2))$.

Up to this point we can endow $Pol(SU_q(2)$ with a different star structure, let us denote by $Pol(SU'_q(2))$ this star algebra and denote by $\overset{I}{k}{}^m_n$ the element $\overset{I}{u}{}^m_n$ when considered as element of $Pol(SU'_q(2))$ : the involution is then defined as $(\overset{I}{k}{}^m_n)^\star = S^{-1}(\overset{I}{k}{}^n_m)$.

As a result it is natural to define the normed star algebra $Fun_{cc}(SL_q(2,\mathbb{C})_\mathbb{R}) = Pol(SU'_q(2)\otimes Fun_c(AN_q(2))$, the normed star algebra $Fun_c(SL_q(2,\mathbb{C})_\mathbb{R}) = Fun(SU'_q(2))\otimes Fun_c(AN_q(2))$ where the tensor product $\otimes$ always denotes the algebraic tensor product.

Naturally we can define the $C^\star$ algebra $Fun_0(SL_q(2,\mathbb{C})_\mathbb{R}) = Fun(SU'_q(2))\otimes Fun_0(AN_q(2))$ and denote by $Fun(SL_q(2,\mathbb{C})_\mathbb{R})$ the set of affiliated elements to $Fun_0(SL_q(2,\mathbb{C})_\mathbb{R})$. We have $Fun(SL_q(2,\mathbb{C})_\mathbb{R}) = \prod_{I\in\frac{1}{2}\mathbb{Z}^+} Fun(SU'_q(2))\otimes\mathrm{Mat}_{2\mathrm{I}+1}(\mathbb{C})$. Note that this last algebra is not endowed with a norm.

$Fun_{cc}(SL_q(2,\mathbb{C})_\mathbb{R})$ is endowed with a structure of multiplier Hopf algebra whereas $Fun_0(SL_q(2,\mathbb{C})_\mathbb{R})$ is endowed with a structure of multiplier $C^\star$ Hopf algebra. The coalgebra structure on these last two spaces has to be twisted, according to the structure of the dual of a quantum double. More precisely we have $Fun_{cc}(SL_q(2,\mathbb{C})_\mathbb{R}) = Pol(SU'_q(2)\otimes_\mathcal{F} Fun_c(AN_q(2))$ and $Fun_0(SL_q(2,\mathbb{C})_\mathbb{R}) = Fun(SU'_q(2))\otimes_\mathcal{F} Fun_0(AN_q(2))$ as multiplier coalgebras, with $\mathcal{F}$ the invertible multiplier of $M(\;Fun_c(AN_q(2))\otimes Pol(SU_q(2))\;)$, defined by:

$$\mathcal{F} = \sum_{J,x,y}{}^{\oplus}\overset{J}{E}{}^x_y\otimes\overset{J}{k}{}^y_x, \tag{81}$$

(this element satisfies the equations (12) and is the bicharacter of [36]). All this is described in [36].

Let us now describe explicitely the full structure of multiplier Hopf algebra of $Fun_{cc}(SL_q(2,\mathbb{C})_\mathbb{R})$.

The explicit structure of multiplier algebra of $Fun_{cc}(SL_q(2,\mathbb{C})_\mathbb{R})$ is defined by the following equations:

$$\begin{aligned}
&\overset{A}{k}{}^i_j\otimes\overset{B}{E}{}^k_l\cdot\overset{C}{k}{}^m_n\otimes\overset{D}{E}{}^p_q = \sum_F\begin{pmatrix} i & m & \Big| & F \\ A & C & \Big| & r\end{pmatrix}\begin{pmatrix} s & \Big| & A & C \\ F & \Big| & j & n\end{pmatrix}\delta^p_l\delta_{B,D}\overset{F}{k}{}^r_s\otimes\overset{B}{E}{}^k_q\\
&\Delta(\overset{A}{k}{}^i_j\otimes\overset{B}{E}{}^k_l) = \mathcal{F}^{-1}_{23}\;(\sum_{C,D,m,p,q,r,s}\begin{pmatrix} q & s & \Big| & B \\ C & D & \Big| & l\end{pmatrix}\begin{pmatrix} k & \Big| & C & D \\ B & \Big| & p & r\end{pmatrix}\overset{A}{k}{}^i_m\otimes\overset{C}{E}{}^p_q\otimes\overset{A}{k}{}^m_j\otimes\overset{D}{E}{}^r_s\;)\;\mathcal{F}_{23}\\
&\epsilon(\overset{A}{k}{}^i_j\otimes\overset{B}{E}{}^k_l) = \delta^i_j\delta_{B,0}\\
&(\overset{A}{k}{}^i_j\otimes\overset{B}{E}{}^k_l)^\star = S^{-1}(\overset{A}{k}{}^j_i)\otimes\overset{B}{E}{}^l_k\\
&\text{with }\;\mathcal{F}^{-1}_{12} = \sum_{J,x,y}\overset{J}{E}{}^x_y\otimes S^{-1}(\overset{J}{k}{}^y_x)
\end{aligned}\tag{82}$$

It is easy to check the axioms of a multiplier Hopf algebra.

It is also easy to check that, although the coproduct $\Delta$ maps $Fun_{cc}(SL_q(2,\mathbb{C})_{\mathbb{R}})$ to $M(Fun_{cc}(SL_q(2,\mathbb{C})_{\mathbb{R}})\otimes Fun_{cc}(SL_q(2,\mathbb{C})_{\mathbb{R}}))$, one can still define actions $\overset{R}{\triangleright}$ (resp. $\overset{L}{\triangleright}$) of $\mathfrak{U}_q(sl(2,\mathbb{C})_{\mathbb{R}})$ on $Fun_{cc}(SL_q(2,\mathbb{C})_{\mathbb{R}})$ by the formulas (26).

As an algebra we can define $Pol(SL_q(2,\mathbb{C})_{\mathbb{R}}) = Pol(SU'_q(2)) \otimes \mathfrak{U}_q(su(2))_{op}$, which elements are affiliated elements to $Fun_0(SL_q(2,\mathbb{C})_{\mathbb{R}})$.

## 3.3 Iwasawa decomposition of $SL_q(2,\mathbb{C})_{\mathbb{R}}$

Let $\overset{(I,J)}{\pi} = \rho(\overset{I}{\pi},\overset{J}{\pi})$ be an irreducible representation of $\mathfrak{U}_q(sl(2,\mathbb{C})_{\mathbb{R}})$, its matrix elements $\overset{IJ}{G}$ are affiliated elements of $Fun_0(SL_q(2,\mathbb{C})_{\mathbb{R}})$. We have the decomposition $Pol(SL_q(2,\mathbb{C})_{\mathbb{R}}) = Pol(SU'_q(2)) \otimes \mathfrak{U}_q(su(2))_{op}$ which is in duality with $\mathfrak{U}_q(sl(2,\mathbb{C})_{\mathbb{R}})$. Let us denote by $\overset{I}{T}{}^{(\pm)} \in \mathfrak{U}_q(su(2))_{op}$ the elements $\overset{I}{L}{}^{(\pm)}$ of $\mathfrak{U}_q(su(2))$ when considered as elements of $Pol(SL_q(2,\mathbb{C})_{\mathbb{R}})$.

The duality bracket is defined by:

$$\mathfrak{U}_q(su(2))\hat{\otimes}Pol(SU_q(2))^{op} \times Pol(SU'_q(2)) \otimes \mathfrak{U}_q(su(2))_{op} \to \mathbb{C}$$
$$< \overset{I}{L}{}^{(\pm)}_1 \otimes \overset{A}{g}_{1'}, \overset{J}{k}_2 \otimes \overset{B}{T}{}^{(\pm)}_{2'} >= \overset{IJ}{R}{}^{(\pm)}_{12} \overset{AB}{R'}{}^{(\pm)}_{1'2'} \tag{83}$$

From these definitions it is easy to show that the algebra law of $Pol(SL_q(2,\mathbb{C})_{\mathbb{R}})$ on the elements $\overset{I}{k}$ and $\overset{I}{T}{}^{(\pm)}$ satisfies:

$$\overset{I}{k}_1\overset{J}{k}_2 = \sum_K \phi^{IJ}_K \overset{K}{k}\, \psi^K_{IJ}, \tag{84}$$
$$\overset{I}{T}{}^{(\pm)a}{}_b \overset{J}{T}{}^{(\pm)c}{}_d = \sum_K \left(\begin{array}{cc|c} c & a & K \\ J & I & r \end{array}\right) \overset{K}{T}{}^{(\pm)r}{}_s \left(\begin{array}{c|cc} s & J & I \\ K & d & b \end{array}\right) \tag{85}$$
$$\overset{IJ}{R}{}^{ax}_{by}\, \overset{J}{T}{}^{(-)y}{}_z \overset{I}{T}{}^{(+)b}{}_c = \overset{I}{T}{}^{(+)a}{}_b \overset{J}{T}{}^{(-)x}{}_y\, \overset{IJ}{R}{}^{yb}_{cz} \tag{86}$$
$$\overset{I}{k}_1 \overset{J}{T}{}^{(\pm)}_2 = \overset{J}{T}{}^{(\pm)}_2 \overset{I}{k}_1 \tag{87}$$

(the reader should notice that the formulas for the product of $\overset{I}{T}{}^{(\pm)}$ are inverted with respect to those of $\overset{I}{L}{}^{(\pm)}$, this is because $\overset{I}{T}{}^{(\pm)} \in \mathfrak{U}_q(su(2))_{op}$).

The action of the untwisted coproduct is given by:

$$\Delta^0(\overset{I}{k}{}^a_b) = \sum_c \overset{I}{k}{}^a_c \otimes \overset{I}{k}{}^c_b \quad \Delta^0(\overset{I}{T}{}^{(\pm)a}{}_b) = \sum_c \overset{I}{T}{}^{(\pm)c}{}_b \otimes \overset{I}{T}{}^{(\pm)a}{}_c \tag{88}$$

From the duality bracket, it is easy to show that

$$< \overset{I}{k}_2 \overset{I}{T}{}^{(-)-1}_2 \overset{J}{k}_3 \overset{J}{T}{}^{(+)-1}_3, \overset{K}{g}_1 >= \overset{KJ}{R}{}^{(-)}_{13} \overset{KI}{R}{}^{(+)}_{12} \quad < \overset{I}{k}_2 \overset{I}{T}{}^{(-)-1}_2 \overset{J}{k}_3 \overset{J}{T}{}^{(+)-1}_3, \overset{K}{L}{}^{\pm}_1 >= \overset{KJ}{R}{}^{(\pm)}_{13} \overset{KI}{R}{}^{(\pm)}_{12}, \tag{89}$$

comparing these expressions with (64) we deduce that $\overset{IJ}{G}_{12} = \overset{I}{k}_1 \overset{I}{T}{}^{(-)-1}_1 \overset{J}{k}_2 \overset{J}{T}{}^{(+)-1}_2$.

If we take $(I,J) = (\frac{1}{2},0)$ this is exactly the quantum analog of the Iwasawa decomposition, first found, by direct computation in [36]. The fact that the last expressions are containing $\overset{I}{T}{}^{(\pm)-1}$ comes from the fact that the coproduct on $\overset{I}{T}{}^{(\pm)}$ is inverted with respect to the usual coproduct on $\overset{I}{k}$. With our conventions the matrix $\overset{I}{T}{}^{(-)}$ is lower triangular whereas $\overset{I}{T}{}^{(+)}$ is upper triangular, and their expressions as affiliated elements of $Fun_0(AN_q(2))$ is the following:

$$\overset{I}{T}{}^{(\pm)-1}{}^{a}_{b} = \sum_{J\in\frac{1}{2}\mathbb{Z}^+}{}^{\oplus} \overset{IJ}{R}{}^{(\pm)-1}{}^{ax}_{by} \overset{J}{E}{}^{y}_{x} \tag{90}$$

## 3.4 Characters of finite dimensional representations of $\mathfrak{U}_q(sl(2,\mathbb{C})_\mathbb{R})$

We will be particularly interested in the trace of the representation $\overset{(A,B)}{\pi}$ which can be written as $\overset{(A,B)}{\chi} = tr_{11'} \overset{AB}{G}{}_{11'}$ : it is an affiliated element of $Fun_0(SL_q(2,\mathbb{C})_\mathbb{R})$ which can therefore be represented as an element of $\prod_{I\in\frac{1}{2}\mathbb{Z}^+} Fun(SU_q(2))\otimes Mat_{2I+1}(\mathbb{C})$. We will now use the generalized Iwasawa decomposition proved above to explicitely find the expression of this affiliated element.

Theorem 1: Expression of the characters of finite dimensional representations

*$\overset{(A,B)}{\chi}$ is an affiliated element of $Fun_0(SL_q(2,\mathbb{C})_\mathbb{R})$ which can be explicitely computed in term of $6j$ symbols of $SU_q(2)$ :*

$$\overset{(A,B)}{\chi} = \sum_K{}^{\oplus} \sum_{M,N} {}_f\Lambda^{KN}_M(A,B) \begin{pmatrix} n & z \\ M & K \end{pmatrix}\!\!\left|\begin{matrix} N \\ u \end{matrix}\right) \left(\begin{matrix} u \\ N \end{matrix}\right|\!\! \begin{pmatrix} K & M \\ y & m \end{pmatrix} \overset{M}{k}{}^{m}_{n} \otimes \overset{K}{E}{}^{y}_{z} \tag{91}$$

*where ${}_f\Lambda^{KN}_M(A,B) \in \mathbb{C}$ is given by*

$$ {}_f\Lambda^{KN}_M(A,B) = \sum_D \left\{\begin{matrix} A & B \\ K & N \end{matrix}\middle| \begin{matrix} M \\ D \end{matrix}\right\}^2 \frac{v_D}{v_B}\left(\frac{v_M}{v_N v_K}\right)^{\frac{1}{2}} = \sum_D \left\{\begin{matrix} B & A \\ K & N \end{matrix}\middle| \begin{matrix} M \\ D \end{matrix}\right\}^2 \frac{v_A}{v_D}\left(\frac{v_N v_K}{v_M}\right)^{\frac{1}{2}}, \tag{92}$$

*(the index f stands for finite dimensional). Note that in the previous equation (92) the sum on $D$ is finite, and that if $K$ is fixed then the sum in (91) involves just a finite number of $M$.*

Proof: We have:

$$\overset{(A,B)}{\chi} = tr_A(\overset{A}{k}\overset{A}{T}{}^{(-)-1}) tr_B(\overset{B}{k}\overset{B}{T}{}^{(+)-1}) = \overset{A}{k}{}^{a}_{b}\overset{A}{k}{}^{c}_{d} \otimes \sum_{K,L}{}^{\oplus} \overset{AK}{R}{}^{(-)-1}{}^{bx}_{ay} \overset{BL}{R}{}^{(+)-1}{}^{dz}_{ct} \overset{K}{E}{}^{y}_{x} \overset{L}{E}{}^{t}_{z}$$

$$= \sum_K{}^{\oplus} \sum_M \begin{pmatrix} a & c \\ A & B \end{pmatrix}\!\!\left|\begin{matrix} M \\ m \end{matrix}\right) \overset{AK}{R}{}^{(-)-1}{}^{bx}_{ay} \overset{BK}{R}{}^{(+)-1}{}^{dz}_{cx} \left(\begin{matrix} n \\ M \end{matrix}\right|\!\! \begin{pmatrix} A & B \\ b & d \end{pmatrix} \overset{M}{k}{}^{m}_{n} \otimes \overset{K}{E}{}^{y}_{z}$$

$$= \sum_{KM} \left[ \vphantom{\Big|} \right]_v \overset{M}{k}{}^{d}_{e} \otimes \overset{K}{E}{}^{c}_{b}.$$

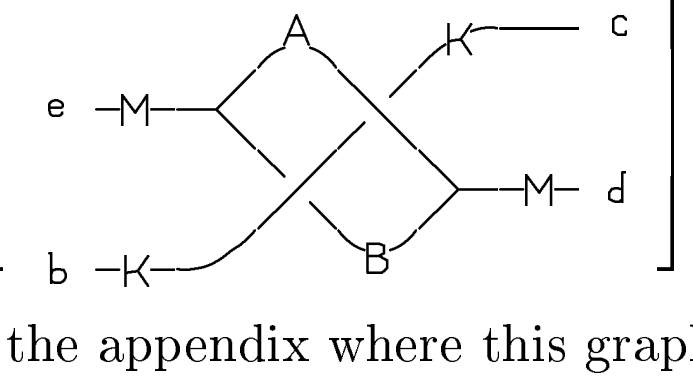

The reader is invited to look at the appendix where this graphical representation is explained. This last expression can be written, using twice the unitarity on Clebsch-Gordan coefficients as:

$$\sum_{K,M,N} {}_f\Lambda^{KN}_M(A,B) \begin{pmatrix} n & z \\ M & K \end{pmatrix}\!\!\left|\begin{matrix} N \\ u \end{matrix}\right) \left(\begin{matrix} u \\ N \end{matrix}\right|\!\! \begin{pmatrix} K & M \\ y & m \end{pmatrix} \overset{M}{k}{}^{m}_{n} \otimes \overset{K}{E}{}^{y}_{z}$$

where ${}_f\Lambda_M^{KN}(A,B)$ is defined by:

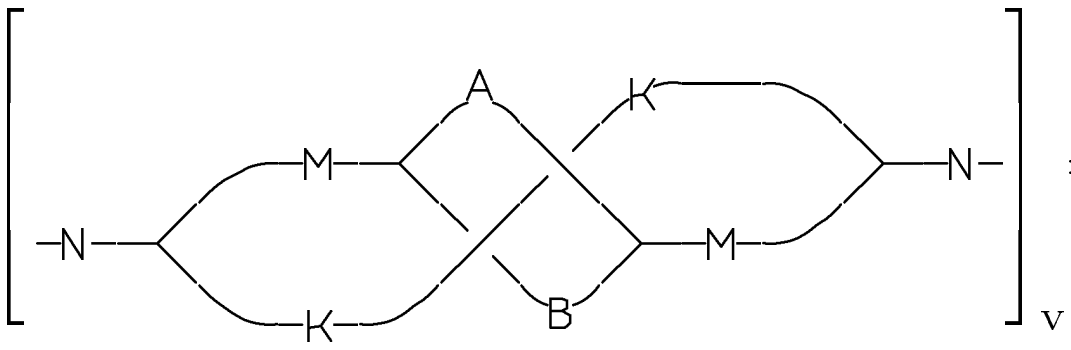

$$= {}_f\Lambda_M^{KN}(A,B)\mathrm{id}_{\overset{N}{V}} \qquad (93)$$

The left handside is equal to:

$$\sum_D$$

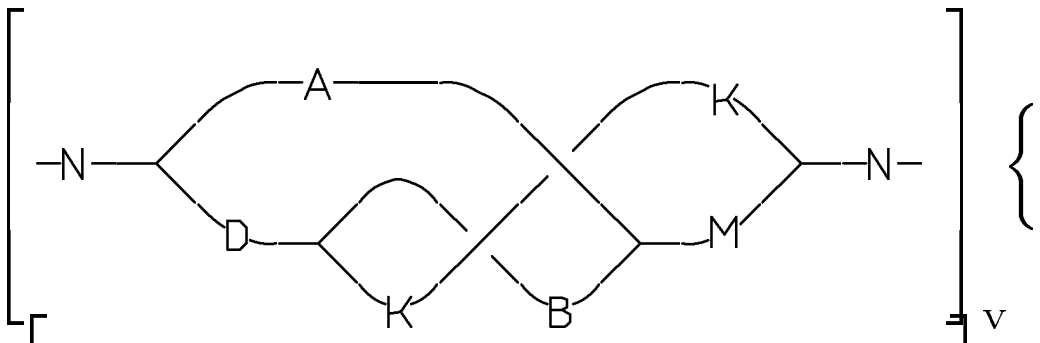

$$\left\{\begin{matrix} A & B \\ K & N \end{matrix}\middle| \begin{matrix} M \\ D \end{matrix}\right\} =$$

$$= \sum_D$$

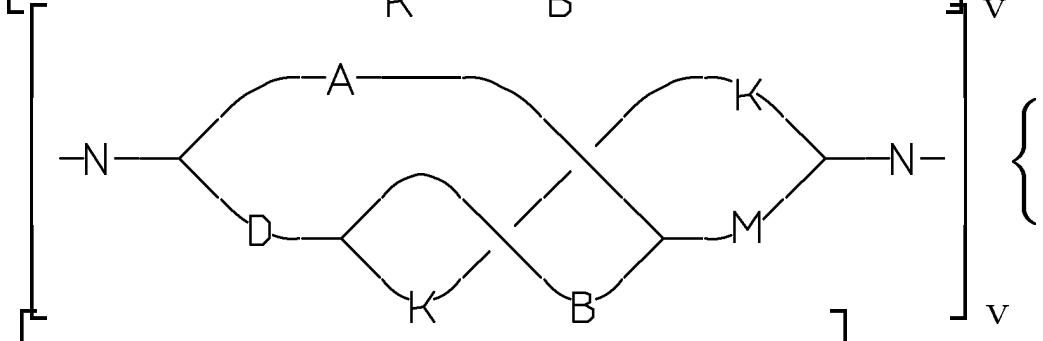

$$\left\{\begin{matrix} A & B \\ K & N \end{matrix}\middle| \begin{matrix} M \\ D \end{matrix}\right\} \frac{v_D}{v_B v_K} =$$

$$= \sum_D$$

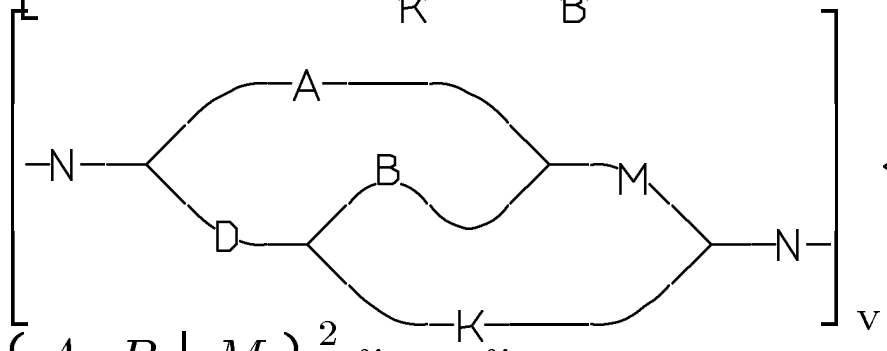

$$\left\{\begin{matrix} A & B \\ K & N \end{matrix}\middle| \begin{matrix} M \\ D \end{matrix}\right\} \frac{v_D}{v_B}\left(\frac{v_M}{v_N v_K}\right)^{\frac{1}{2}}$$

$$= \sum_D \left\{\begin{matrix} A & B \\ K & N \end{matrix}\middle| \begin{matrix} M \\ D \end{matrix}\right\}^2 \frac{v_D}{v_B}\left(\frac{v_M}{v_N v_K}\right)^{\frac{1}{2}}.$$

From the definition of the 6j symbols it is easy to show that if *Y(A,B,M).Y(K,N,D).Y(K,B,D).Y(A,B,D)=* 0 then $\left\{\begin{matrix} A & B \\ K & N \end{matrix}\middle| \begin{matrix} M \\ D \end{matrix}\right\} = 0$. As a result these constraints impose that if $K$ is fixed then the sum in (91) involves just a finite number of $M$. This can be written as a picture in terms of coefficients of $R$ matrix in IRF representation. Indeed, it is trivial to show using the formula of ${}_f\Lambda_M^{KN}(A,B)$ and the formulas for $R$ matrix in IRF picture that:

$${}_f\Lambda_M^{KN}(A,B) = \sum_D$$

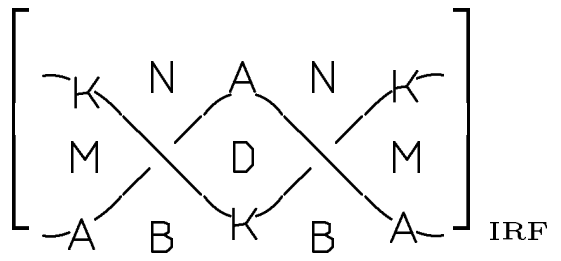

$$\left(\frac{v_N}{v_M v_K}\right)^{\frac{1}{2}} \qquad (94)$$

We can also express ${}_f\Lambda_M^{KN}(A,B)$ with the following formula:

$${}_f\Lambda_M^{KN}(A,B) = \sum_D$$

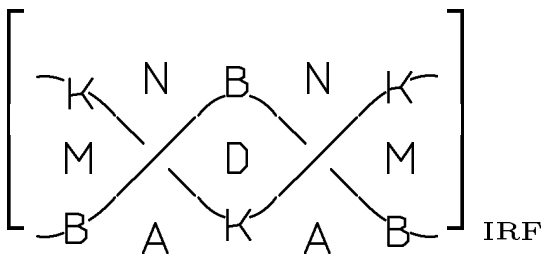

$$\left(\frac{v_M v_K}{v_N}\right)^{\frac{1}{2}}. \qquad (95)$$

In order to show this last equation we have to modify a little bit the previous proof. From equation (93), we have:

$${}_f\Lambda_M^{KN}(A,B)\mathrm{id}_{\overset{N}{V}} =$$

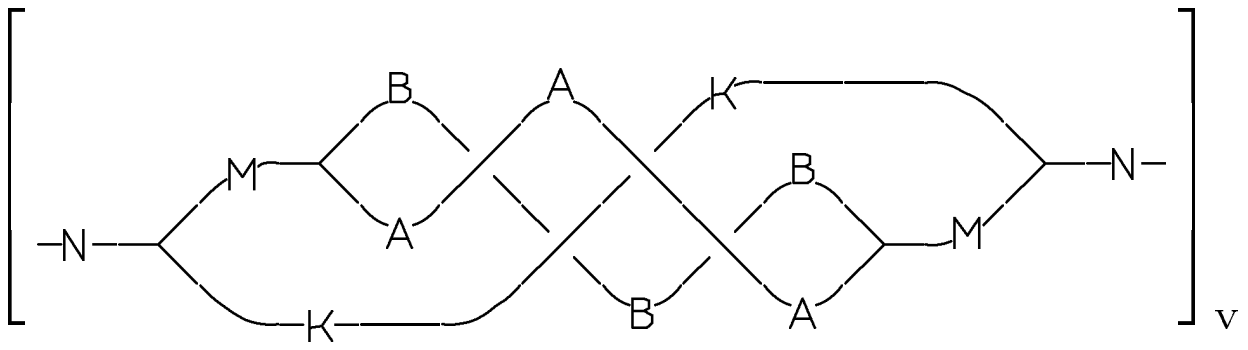

$=$ 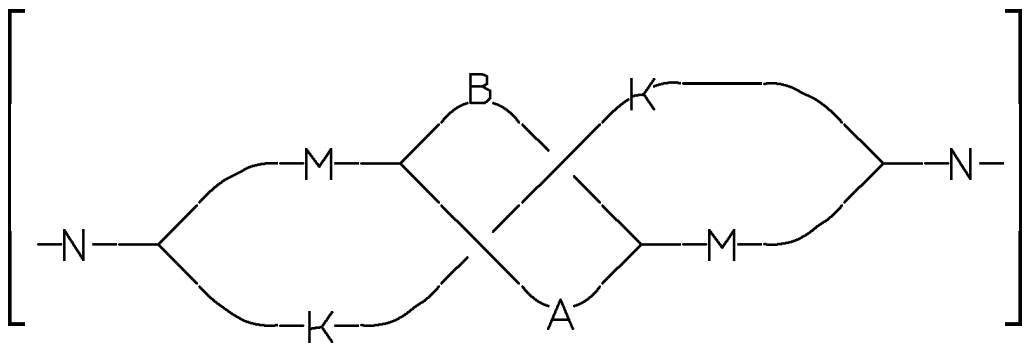 $_{\mathrm{V}}$ .

From this diagram we can apply the succession of moves of the previous proof, and we finally obtain a new expression for ${}_f\Lambda^{KN}_M(A,B)$ :

$${}_f\Lambda^{KN}_M(A,B)=\sum_D\left\{\begin{matrix} B & A \\ K & N\end{matrix}\middle|\begin{matrix} M \\ D\end{matrix}\right\}^2\frac{v_A}{v_D}\left(\frac{v_N v_K}{v_M}\right)^{\frac{1}{2}}=\sum_D$$

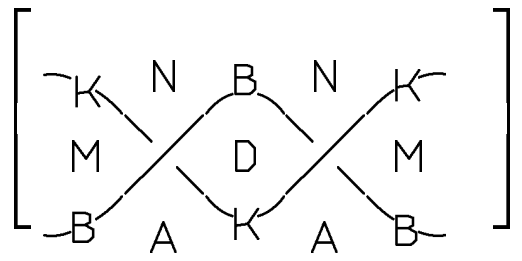

$$_{\mathrm{IRF}}\left(\frac{v_M v_K}{v_N}\right)^{\frac{1}{2}}.$$

□

Let us now study the case of the strange representations $\overset{(A,B)}{\pi}\circ\tau_s$. It is easy to show from the definition of $\tau_s$ that the character of the one dimensional representation $\epsilon_s$ is the affiliated element $\sum^{\oplus}_{J\in\frac{1}{2}\mathbb{Z}^+}(-1)^{2J}\overset{J}{E}{}^x_x$ (this is the way it appears in [36]). As a result from the fact that $\overset{(A,B)}{\pi}\circ\tau_s=\overset{(A,B)}{\pi}\otimes\epsilon_s$, we obtain that the character of the representation $\overset{(A,B)}{\pi}\circ\tau_s$ can be written as:

$$\sum_K^{\oplus}\sum_{M,N}(-1)^{2K}{}_f\Lambda^{KN}_M(A,B)\left(\begin{matrix} n & z \\ M & K\end{matrix}\middle|\begin{matrix} N \\ u\end{matrix}\right)\left(\begin{matrix} u \\ N\end{matrix}\middle|\begin{matrix} K & M \\ y & m\end{matrix}\right)\overset{M}{k}{}^m_n\otimes\overset{K}{E}{}^y_z. \quad (96)$$

## 3.5 Left and right integral on $SL_q(2,\mathbb{C})_{\mathbb{R}}$ and its classical limit

One of the main result of [36] is the proof that there exists a right and left integral on $Fun_c(SL_q(2,\mathbb{C})_{\mathbb{R}})$.

Let $h_{ANr}$ and $h_{ANl}$ be the linear forms defined on $Fun_c(AN_q(2))$ by the formula

$$h_{ANr}(\overset{I}{E}{}^r_s)=[d_I]\overset{I}{\mu}{}^{-1r}{}_s,\quad h_{ANl}(\overset{I}{E}{}^r_s)=[d_I]\overset{I}{\mu}{}^r_s \quad (97)$$

then these linear forms are respectively right and left integral of the discrete quantum group $Fun_c(AN_q(2))$.

A complete study of existence and uniqueness of right and left integral has been done in the case of discrete quantum group in [44].

Let $h_K$ be the unique normalized right and left integral on $Fun(SU'_q(2))$ . It is defined on $Fun_c(SU_q(2))$ by $h_K(\overset{I}{k}{}^r_s)=\delta_{I,0}$. Then the linear form $h=h_K\otimes h_{ANr}$, satisfies:

$$h(\overset{I}{k}{}^a_b\otimes\overset{J}{E}{}^r_s)=\delta_{I,0}[d_J]\overset{J}{\mu}{}^{-1r}{}_s, \quad (98)$$

This linear form is a left and right integral on $Fun_c(SL_q(2,\mathbb{C})_{\mathbb{R}})$. This result as already been proved in [36], but we think that the following computation which is straighforward gives an alternative proof of this important theorem.

Proof: Let $h$ be the linear form defined by (98). We have:

$$(h\otimes \mathrm{id})\Delta(\overset{I}{k}{}^m_p \otimes \overset{C}{E}{}^r_s) =$$

$$= (h\otimes \mathrm{id}) \sum_{A,B,M,N} \left(\begin{array}{cc|c} a' & b' & C \\ A & B & s \end{array}\right) \left(\begin{array}{c|cc} r & A & B \\ C & a & b \end{array}\right) \overset{I}{k}{}^m_n \otimes \overset{M}{E}{}^x_y \overset{A}{E}{}^a_{a'} \overset{N}{E}{}^z_t \otimes \overset{M}{k}{}^y_x \overset{I}{k}{}^n_p S^{-1}(\overset{N}{k}{}^t_z) \otimes \overset{B}{E}{}^b_{b'}$$

$$= \sum_{A,B} \left(\begin{array}{cc|c} a' & b' & C \\ A & B & s \end{array}\right) \left(\begin{array}{c|cc} r & A & B \\ C & a & b \end{array}\right) \delta_{I,0} \overset{A}{\mu}{}^{-1}{}^x_t \overset{A}{k}{}^a_x S^{-1}(\overset{A}{k}{}^t_{a'}) \otimes \overset{B}{E}{}^b_{b'} [d_A]$$

$$= \sum_{A,B} \left(\begin{array}{cc|c} a' & b' & C \\ A & B & s \end{array}\right) \left(\begin{array}{c|cc} r & A & B \\ C & a & b \end{array}\right) \delta_{I,0} \overset{A}{\mu}{}^{-1}{}^a_{a'} 1 \otimes \overset{B}{E}{}^b_{b'} [d_A]$$

$$= \sum_A$$

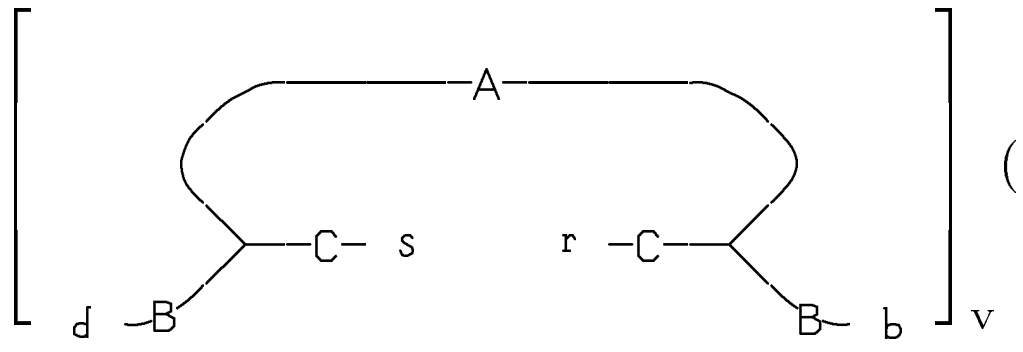

$$(-1)^{-2A}[d_A]\delta_{I,0} 1 \otimes \overset{B}{E}{}^b_d$$

$$= \sum_A$$

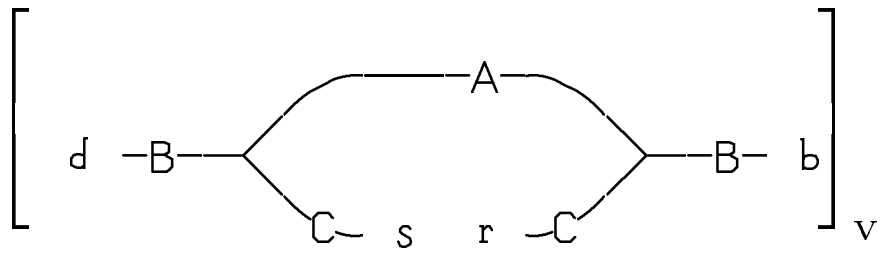

$$(-1)^{-2A+2C-2B}\delta_{I,0}\frac{[d_A][d_C]}{[d_B]} 1 \otimes \overset{B}{E}{}^b_d$$

$$= \sum_A$$

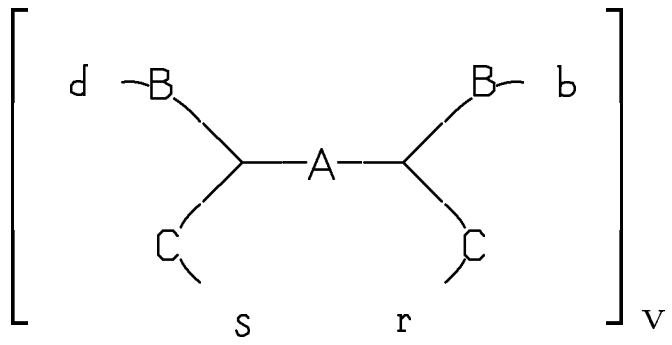

$$(-1)^{2A-2B}\delta_{I,0}[d_C] 1 \otimes \overset{B}{E}{}^b_d$$

$$=$$

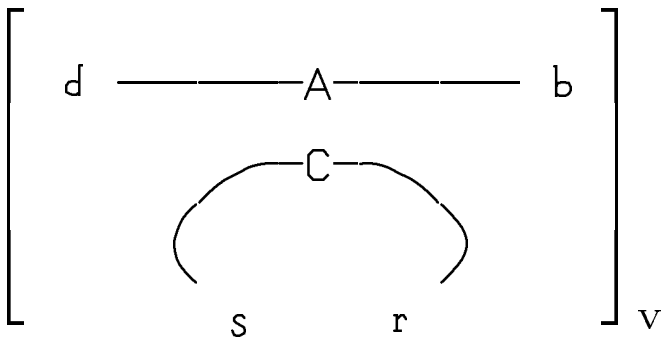

$$(-1)^{2C}\delta_{I,0}[d_C] 1 \otimes \overset{B}{E}{}^b_d$$

$$= \sum_{B,b,d} \overset{B}{E}{}^b_d \delta_{I,0}[d_C] 1 \otimes \overset{C}{\mu}{}^{-1}{}^r_s = h(\overset{I}{k}{}^m_p \otimes \overset{C}{E}{}^r_s) 1$$

□

Using this integral we can define on $Fun_{cc}(SL_q(2,\mathbb{C})_{\mathbb{R}})$ an hermitian form defined by $(a,b) = h(a^\star b), \forall a,b \in Fun_{cc}(SL_q(2,\mathbb{C})_{\mathbb{R}})$. As a result $Fun_{cc}(SL_q(2,\mathbb{C})_{\mathbb{R}})$ is equiped with a norm $|| \, ||_{L^2}$ defined by $||a||_{L^2} = (a,a)^{\frac{1}{2}}$. We will denote by $L^2(SL_q(2,\mathbb{C})_{\mathbb{R}})$ the completion of $Fun_{cc}(SL_q(2,\mathbb{C})_{\mathbb{R}})$ with respect to this norm.

One puzzling aspect of the construction of $Fun_c(SL_q(2,\mathbb{C})_{\mathbb{R}})$ is that we know from general arguments that this multiplier Hopf algebra should correspond to the quantum analog of the compact supported functions on $SL(2,\mathbb{C})_{\mathbb{R}}$, but the precise limit $q \to 1$ is not so easy to understand. Although we precisely know what are the classical limit of the elements $\overset{I}{k}$, it is not at all obvious what is the limit, if any, of the elements $\overset{I}{E}$. In particular we would like to understand in which sense the formula of $h$ gives, in the classical limit, the well known left and right invariant Haar measure of $SL(2,\mathbb{C})_{\mathbb{R}}$. In the rest of this part we will explain how to understand the classical limit of the expressions (97).

An element $g \in AN(2)$ can be written as $g = \begin{pmatrix} a & 0 \\ n & a^{-1} \end{pmatrix}$ with $a \in \mathbb{R}^{+*}, n \in \mathbb{C}$. As usual $a$ and $n$ can be thought as being functions of $g$.

The expressions of the right and left Haar measure $dm_r$ and $dm_l$ on $Fun_c(AN(2))$ are given by

$$\int \phi(g) dm_l(g) = \int \phi(g) a \; da \; d^2n \;\; , \;\; \int \phi(g) dm_r(g) = \int \phi(g) a^{-3} \; da \; d^2n. \tag{99}$$

The aim of what follows is to prove that the expressions (97) are certain types of Jackson integral discretizing the previous classical Haar measures.

A straightforward computation, shows that $\overset{\frac{1}{2}}{T}{}^{(-)-1} = \begin{pmatrix} \hat{a} & 0 \\ (\frac{q^2+1}{2})^{\frac{1}{2}}\hat{n} & \hat{a}^{-1} \end{pmatrix}$ with the commutation relations:

$$\hat{a}\hat{a}^\star = \hat{a}^\star\hat{a}, \quad \hat{a}\hat{n} = q^{-1}\hat{n}\hat{a}, \quad \hat{a}\hat{n}^\star = q\hat{n}^\star\hat{a}, \quad \hat{n}\hat{n}^\star - \hat{n}^\star\hat{n} = \frac{2}{q+q^{-1}}(q^{-1}-q)(\hat{a}^2 - \hat{a}^{-2}). \tag{100}$$

The algebra generated by $\hat{a}, \hat{n}, \hat{n}^\star$ is isomorphic to $\mathfrak{U}_q(su(2))$ when $q \neq 1$ by quantum duality principle.

In the limit where $q \to 1$ the elements $\hat{a}, \hat{n}, \hat{n}^\star$ commute as they should. The element $\hat{C} = \frac{1}{2}(\hat{n}\hat{n}^\star + \hat{n}^\star\hat{n}) + (\hat{a}^2 + \hat{a}^{-2})$ is in the center of this algebra. The elements $\hat{C}$ and $\hat{a}$ are both self-adjoint positive affiliated elements and commute. The projections of the affiliated elements $\hat{a}, \hat{n}, \hat{n}^\star$ on each component $Mat_{2I+1}(\mathbb{C}) \subset Fun_c(AN(2))$ are given by the following expression:

$$\hat{a}_I = \overset{I}{\pi}(q^{J_z}) \qquad \hat{n}_I = q^{\frac{1}{2}}(q^2-1)(\frac{2}{q+q^{-1}})^{\frac{1}{2}}\overset{I}{\pi}(J_-) \qquad \hat{n}_I^\star = q^{-\frac{1}{2}}(q^2-1)(\frac{2}{q+q^{-1}})^{\frac{1}{2}}\overset{I}{\pi}(J_+). \tag{101}$$

Let $\phi : \mathbb{R}^+ \times \mathbb{R}^+ \to \mathbb{R}$ be continuous with compact support, we can define $\phi(\hat{a}, \hat{C})$ to be the affiliated element which components are $(\phi(\hat{a}, \hat{C}))_I = \phi(\hat{a}_I, \hat{C}_I)$. It will also be convenient to define $\phi(\underline{a, C})$ to be the function on $AN(2)$ defined by $\phi(a, C)(g) = \phi(a(g), C(g))$ where $C(g) = n(g)\underline{n(g)} + a^2(g) + a^{-2}(g)$.

We have the following theorem:

Theorem 2:

*The classical left and right measure of $AN(2)$ are related to the quantum left and right measure Haar measure by the following formula:*

$$\lim_{q \to 1} (q - q^{-1})^3 \;\; \pi \;\; h_{ANr}(\phi(\hat{a}, \hat{C})) = \int \phi(a, C) dm_r(g) \tag{102}$$

$$\lim_{q \to 1} (q - q^{-1})^3 \;\; \pi \;\; h_{ANl}(\phi(\hat{a}, \hat{C})) = \int \phi(a, C) dm_l(g) \tag{103}$$

*where $\phi$ is any continuous function $\mathbb{R}^+ \times \mathbb{R}^+ \to \mathbb{R}$ with compact support*

Proof: The right handside of the first equality is equal to

$$2\pi \int_0^{+\infty} \int_0^{+\infty} \phi(a, \rho^2 + a^2 + a^{-2}) \rho a^{2\epsilon_{r,l}-1} \; d\rho \; da$$

where $\epsilon_r = -1, \epsilon_l = 1$. After the change of variable $\rho^2 + a^2 + a^{-2} = 2\cosh(x), a = \exp(y)$ this integral is rewritten:

$$4\pi \int_{-\infty}^{+\infty} dx \int_{-x/2}^{x/2} \phi(e^y, 2\cosh x) e^{2\epsilon_{r,l} y} \sinh x \; dy$$

This integral is transformed into a Riemann sum by discretizing $x = \hbar(2J+1), y = m\hbar$ where $q = e^{\hbar}$, $J \in \frac{1}{2}\mathbb{Z}^+$ and $m \in \frac{1}{2}\mathbb{Z}$. As a result we have:

$$\int \phi(a, C) dm_{r,l}(g) = \lim_{\hbar \to 0} \hbar^2 8\pi \sum_{J \in \frac{1}{2}\mathbb{Z}^+} \sum_{m=-J}^{J} \phi(q^m, q^{2J+1} + q^{-2J-1}) q^{\epsilon_{r,l}} (q^{2J+1} - q^{-2J-1}).$$

This last expression can also be written as:

$$\lim_{\hbar \to 0} \hbar^2 8\pi \sum_{J \in \frac{1}{2}\mathbb{Z}^+} \sum_{m=-J}^{J} \phi(q^m, \frac{2}{q+q^{-1}}(q^{2J+1} + q^{-2J-1})) q^{\epsilon_{r,l}} (q^{2J+1} - q^{-2J-1}). \tag{104}$$

because of $\frac{2}{q+q^{-1}} = 1 + O(\hbar)$. An easy computation, using (101) shows that $C_I = \frac{2}{q+q^{-1}}(q^{2I+1} + q^{-2I-1})1$, as a result we obtain that (104) is exactly

$$\lim_{q \to 1} \{ (q - q^{-1})^3 \; \pi \sum_{J \in \frac{1}{2}\mathbb{Z}^+} tr( \; \phi(\hat{a}_J, \hat{C}_J) \; ) \; [d_J] \}. \tag{105}$$

□

It would be very interesting to extend this result to a larger class of functions. The following heuristic argument shows that this is certainly possible but will require a deeper analysis. Assume that $\phi$ is now a function of three variables, we would like to define $\phi(\hat{a}, \hat{n}, \hat{n}^\star)$, this is not possible in general because $\hat{n}, \hat{n}^\star$ do not commute, but the difference of two choices of ordering will be of order $\hbar$ so that when one compute the left handside of (102, 103) in the limit where $\hbar$ goes to zero the choice of ordering is irrelevant. It is easy to see from the structure of the representation $\overset{I}{\pi}$ that the only monomials $\hat{n}^j \hat{n}^{\star k}$ such that $tr_J(\overset{J}{\mu} \hat{n}^j \hat{n}^{\star k})$ is non zero are those for which $j = k$. From the previous argument we can always choose the ordering $(\frac{1}{2}\hat{n}\hat{n}^\star + \frac{1}{2}\hat{n}^\star\hat{n})^j$ and we are back to the case of a function $\phi(a, C)$. In particular this argument explains the presence of the trace in the formula $h_{ANr,l}$ : taking the trace corresponds in the classical limit to taking the mean over the angular variable $\theta$ defined by $n = \rho e^{i\theta}, \rho \in \mathbb{R}^{+*}$.

# 4 Representation theory and characters

## 4.1 Irreducible unitary representations

Irreducible representations of $SL(2,\mathbb{C})_{\mathbb{R}}$ have been first classified by I.M.Gelfand and M.A.Naimark [21], see also the very good monographs [20][37] devoted to $SL(2,\mathbb{C})_{\mathbb{R}}$. Irreducible unitary representations of $\mathfrak{U}_q(sl(2,\mathbb{C})_{\mathbb{R}})$ have been classified in [35]. W.Pusz shows that if $(V, \pi)$ is an irreducible unitary representation of $\mathfrak{U}_q(sl(2,\mathbb{C})_{\mathbb{R}})$ then $\pi(\Omega^\pm) = \omega^\pm id_V$ and $\omega^\pm$ are complex numbers such that $\omega^+ = \bar{\omega}^-$. Moreover the irreducible unitary representations are entirely classified, up to equivalence, by the value of $\omega^+$. More precisely $\omega^+, \omega^-$ can take only the following allowed value $\omega^+ = \frac{[2(2\aleph_0+1)]}{[2\aleph_0+1]}, \omega^- = \frac{[2(2\aleph_1+1)]}{[2\aleph_1+1]}$ with the following constraint on $\aleph_0, \aleph_1$:

1. $\aleph_0 = \aleph_1 = 0$. In this case $\pi$ is the trivial one dimensional representation $\pi = \epsilon$.

2. $\aleph_0 = \aleph_1 = \frac{i\pi}{2\hbar}$. In this case $\pi$ is the one dimensional representation $\pi = \epsilon_s$.

3. $2\aleph_0 + 1 = \frac{1}{2}(m + i\rho) \quad 2\aleph_1 + 1 = \frac{1}{2}(-m + i\rho), \quad m \in \mathbb{N}^*, \rho \in [0, \frac{4\pi}{\hbar}[$ or $m = 0, \rho \in [0, \frac{2\pi}{\hbar}]$. $\pi$ is an infinite dimensional representation belonging to the so called principal series. We will denote by $\hat{\mathcal{D}}_p$ the set of representations belonging to the principal series.

4. $2\aleph_0 + 1 = 2\aleph_1 + 1 = \rho \in ]0, 1[$. $\pi$ is an infinite dimensional representation belonging to the so called complementary series.

5. $2\aleph_0 + 1 = 2\aleph_1 + 1 = \rho + \frac{i\pi}{\hbar}, \rho \in ]0, 1[$. $\pi$ is an infinite dimensional representation belonging to the so called strange complementary series.

One can map one complementary series to the other one by tensoring it with the representation $\epsilon_s$. The principal series and the standard complementary series, after proper rescaling of the generators tend, when $q$ goes to one, to the principal and complementary series of $SL(2,\mathbb{C})_{\mathbb{R}}$, whereas the limit of the strange complementary representations, when $q \to 1$, is singular and disappear in this limit.

The proof follows the following path: when $V$ is infinite dimensional its is shown that $V = \bigoplus_{I=\frac{m}{2}}^{+\infty} \overset{I}{V}$, with $\pi(a) = \oplus_I \overset{I}{\pi}(a)$ for $a \in \{q^{J_z}, J_+, J_-\}$ (note that the multiplicity are always 1) and by a Wigner-Eckart theorem we have that $\pi(\overset{\frac{1}{2}}{g}{}^i_j)$ maps $\overset{K}{V}$ to $\overset{K-1}{V} \oplus \overset{K}{V} \oplus \overset{K+1}{V}$. Explicit formulas for the action of $\overset{\frac{1}{2}}{g}{}^i_j$ are computed. We will show in this section that one can find expressions of $\pi(\overset{I}{g}{}^i_j)$ for any $I \in \frac{1}{2}\mathbb{Z}^+$ and this is nicely expressed in term of continuation of $6j$ symbols. At this point the reader is invited to read the appendix where we give an introduction to complex continuation of $6j$ symbols and their graphical representations, and to read the articles [5, 6, 14, 23] where this theory is developped.

Theorem 3: Unitary Irreducible representations of $\mathfrak{U}_q(sl(2,\mathbb{C})_\mathbb{R})$

*Let $\pi$ be an infinite dimensional unitary irreducible representation of $\mathfrak{U}_q(sl(2,\mathbb{C})_\mathbb{R})$, and assume that $\aleph_0, \aleph_1$ take one of the value described above and associated to the principal series or to one of the two complementary series. $\pi$ is equivalent to the following representation of $\mathfrak{U}_q(sl(2,\mathbb{C})_\mathbb{R})$ : the representation space is given by $V = \bigoplus_{A\in\frac{1}{2}\mathbb{Z}^+, A\pm\frac{m}{2}\in\mathbb{N}} \overset{A}{V}$ and the action of generators in an orthonormal basis $\{\overset{C}{e}_r, C\in\frac{1}{2}\mathbb{Z}^+, r=-C,\dots,C\}$ is given by*

$$\overset{B}{L}{}^{(\pm)i}{}_j \cdot \overset{C}{e}_r = \overset{C}{e}_n \overset{BC}{R}{}^{(\pm)in}_{jr} \tag{106}$$

$$\overset{B}{g}{}^i_j \cdot \overset{C}{e}_r = \sum_{DE} \overset{E}{e}_p \begin{pmatrix} p & i \\ E & B \end{pmatrix}\!\!\left|\begin{matrix} D \\ x\end{matrix}\right. \begin{pmatrix} x \\ D \end{pmatrix}\!\!\left|\begin{matrix} B & C \\ j & r\end{matrix}\right. \Lambda^{BD}_{EC}(\aleph_0,\aleph_1) \tag{107}$$

*where the complex numbers $\Lambda^{BD}_{EC}(\aleph_0,\aleph_1)$ are defined by*

$$\Lambda^{BC}_{AD}(\aleph_0,\aleph_1) = \sum_{\aleph_2 | Y(B,\aleph_2,\aleph_1)=1} \left\{\begin{matrix} B & C \\ \aleph_0 & \aleph_1 \end{matrix} \middle| \begin{matrix} A \\ \aleph_2 \end{matrix}\right\} \left\{\begin{matrix} B & C \\ \aleph_0 & \aleph_1 \end{matrix} \middle| \begin{matrix} D \\ \aleph_2 \end{matrix}\right\} \frac{v_{\aleph_2}}{v_{\aleph_1}} \frac{v_A^{\frac14} v_D^{\frac14}}{v_B^{\frac12} v_C^{\frac12}} \tag{108}$$

$$= \sum_{\aleph_2 | Y(B,\aleph_2,\aleph_1)=1} \left\{\begin{matrix} B & C \\ \aleph_1 & \aleph_0 \end{matrix} \middle| \begin{matrix} A \\ \aleph_2 \end{matrix}\right\} \left\{\begin{matrix} B & C \\ \aleph_1 & \aleph_0 \end{matrix} \middle| \begin{matrix} D \\ \aleph_2 \end{matrix}\right\} \frac{v_{\aleph_0}}{v_{\aleph_2}} \frac{v_B^{\frac12} v_C^{\frac12}}{v_A^{\frac14} v_D^{\frac14}}. \tag{109}$$

*These coefficients satisfy the following relations:*

$$\sum_E \left\{\begin{matrix} A & E \\ B & U \end{matrix} \middle| \begin{matrix} C \\ D \end{matrix}\right\} \Lambda^{AC}_{FE} \Lambda^{BD}_{EP} = \sum_K \left\{\begin{matrix} F & A \\ B & U \end{matrix} \middle| \begin{matrix} C \\ K \end{matrix}\right\} \left\{\begin{matrix} A & B \\ P & U \end{matrix} \middle| \begin{matrix} K \\ D \end{matrix}\right\} \Lambda^{KU}_{FP} \tag{110}$$

$$\sum_D \Big(\frac{v_A v_{\frac12}}{v_D}\Big)^{\pm\frac12} \frac{[2D+1]}{[2A+1]} \Lambda^{\frac12 D}_{AA} = \omega_\pm Y(A,\aleph_0,\aleph_1). \tag{111}$$

<u>Proof:</u> From Pusz classification, we know that $V = \bigoplus_{A\in\frac{1}{2}\mathbb{Z}^+, A\pm\frac{m}{2}\in\mathbb{N}} \overset{A}{V}$. The elements of $\mathfrak{U}_q(su(2))$ act on each $\overset{I}{V}$ by $\overset{I}{\pi}$ and we can always assume that the elements of $\mathfrak{U}_q(an(2))$ act as follows:

$$\overset{C}{g}{}^j_l \cdot \overset{P}{e}_z = \sum_E \xi \begin{bmatrix} EC \\ CP \end{bmatrix}^{uj}_{lp} \overset{E}{e}_u.$$

The relations (106) is the action of the representation $\overset{C}{\pi}$ expressed in term of $\overset{B}{L}{}^{(\pm)}$. We first impose on $\xi$ the constraint obtained from the commutation relations between $\mathfrak{U}_q(su(2))$ and $\mathfrak{U}_q(an(2))$. We have

$$(\overset{BC}{R}{}^{(\pm)ij}_{kl} \overset{B}{L}{}^{(\pm)k}{}_m \overset{C}{g}{}^l_n) \cdot \overset{P}{e}_z = \sum_E \overset{E}{e}_u (\overset{BC}{R}{}^{(\pm)ij}_{kl} \overset{BE}{R}{}^{(\pm)ku}_{mp} \xi \begin{bmatrix} EC \\ CP \end{bmatrix}^{pl}_{nz})$$

$$(\overset{C}{g}{}^j_l \overset{B}{L}{}^{(\pm)i}{}_k \overset{BC}{R}{}^{(\pm)kl}_{mn}) \cdot \overset{P}{e}_z = \sum_E \overset{E}{e}_u (\xi \begin{bmatrix} EC \\ CP \end{bmatrix}^{uj}_{lp} \overset{BP}{R}{}^{(\pm)ip}_{kz} \overset{BE}{R}{}^{(\pm)kl}_{mn})$$

Therefore there exist complex numbers $\Lambda^{CD}_{EP}$ such that $\xi \begin{bmatrix} EC \\ CP \end{bmatrix}^{uj}_{lp} = \sum_D \begin{pmatrix} u & j \\ E & C \end{pmatrix}\!\!\left|\begin{matrix} D \\ x\end{matrix}\right. \begin{pmatrix} x \\ D \end{pmatrix}\!\!\left|\begin{matrix} C & P \\ l & p\end{matrix}\right. \Lambda^{CD}_{EP}$. We now study the constraint on the numbers $\Lambda^{CD}_{EP}$ coming from the algebra law on $\mathfrak{U}_q(an(2))$.

We first have:

$$(\overset{Ai}{g}_j \overset{Bk}{g}_l)\cdot \overset{P}{e}_z = \sum_{FEDC} \left(\begin{array}{cc|c} v & i & C \\ F & A & y \end{array}\right) \left(\begin{array}{c|cc} y & A & E \\ C & j & p \end{array}\right) \left(\begin{array}{cc|c} p & k & D \\ E & B & x \end{array}\right) \left(\begin{array}{c|cc} x & B & P \\ D & l & z \end{array}\right) \Lambda^{AC}_{FE}\Lambda^{BD}_{EP}\, \overset{F}{e}_v$$

$$= (\left(\begin{array}{cc|c} i & k & K \\ A & B & m \end{array}\right) \overset{Km}{g}_n \left(\begin{array}{c|cc} n & A & B \\ K & j & l \end{array}\right)) \cdot \overset{P}{e}_z$$

$$= \left(\begin{array}{cc|c} v & m & R \\ F & K & p \end{array}\right) \left(\begin{array}{c|cc} p & K & P \\ R & n & z \end{array}\right) \left(\begin{array}{cc|c} i & k & K \\ A & B & m \end{array}\right) \left(\begin{array}{c|cc} n & A & B \\ K & j & l \end{array}\right) \Lambda^{KR}_{FP}\, \overset{F}{e}_v$$

We can derive a more compact expression for these constraints. Indeed we have, for the first part of the equality

$$\sum_{DCE} \Lambda^{AC}_{FE}\Lambda^{BD}_{EP}$$

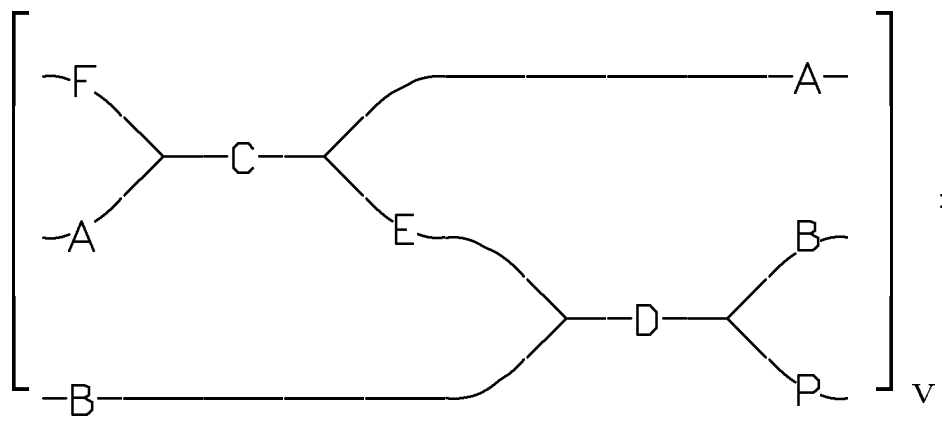

$$=$$

$$= \sum_{DCEU} \Lambda^{AC}_{FE}\Lambda^{BD}_{EP} \left\{\begin{array}{cc|c} A & E & C \\ B & U & D \end{array}\right\}$$

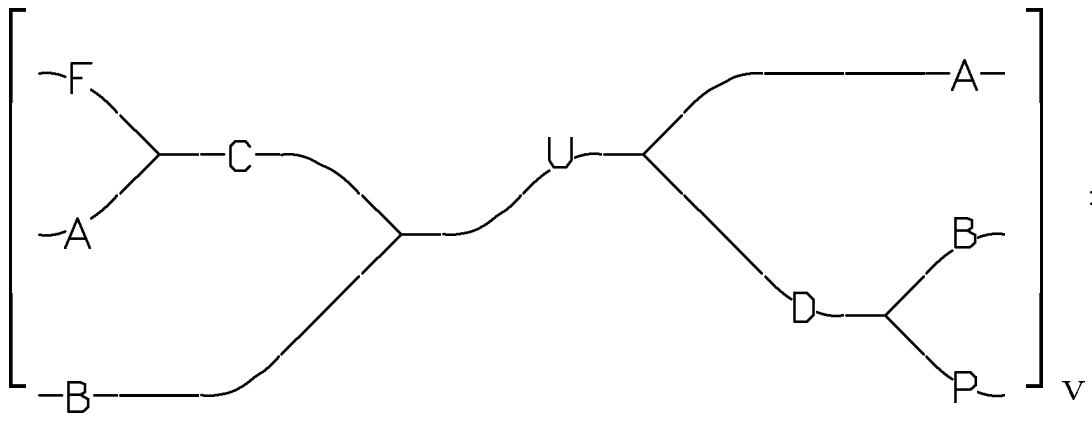

$$=$$

using a similar treatment for the other part, we obtain easily

$$\sum_{KT} \Lambda^{KT}_{FP}$$

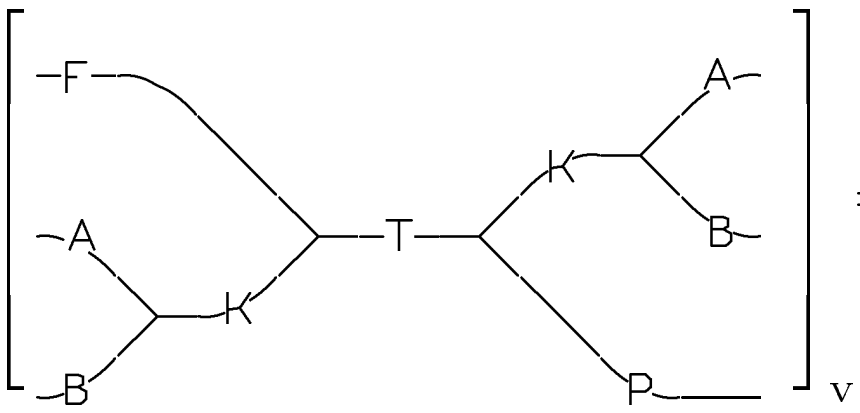

$$=$$

$$= \sum_{KUDC} \Lambda^{KU}_{FP} \left\{\begin{array}{cc|c} F & A & C \\ B & U & K \end{array}\right\} \left\{\begin{array}{cc|c} A & B & K \\ P & U & D \end{array}\right\}$$

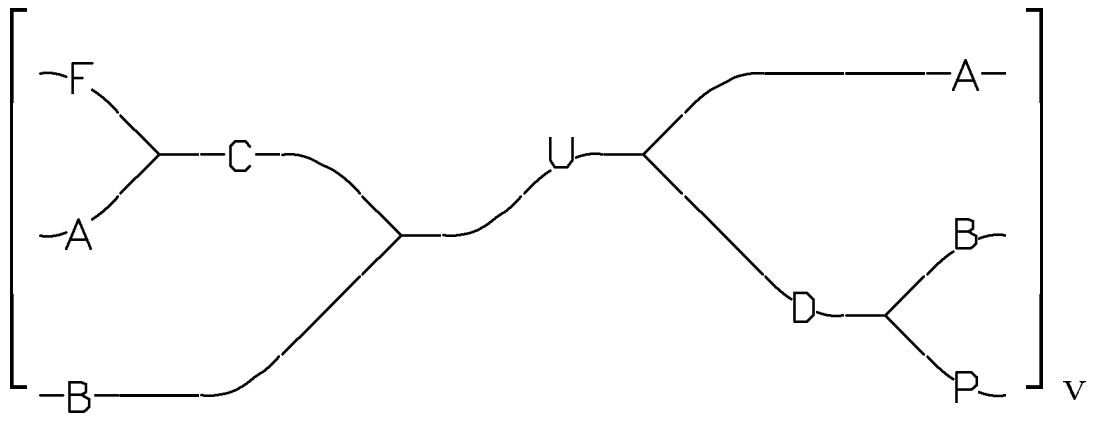

which gives the first announced constraint (110) on coefficients $\Lambda^{CD}_{EP}$.
The second constraint is derived from the action of Casimir elements. Let us write this action explicitely. Using the expressions (106)(107) and the expression of the Casimir elements we obtain

$$\Omega^{\pm} \cdot \overset{B}{e}_j = \sum_D \overset{A\frac{1}{2}}{R}{}^{(\pm)ib}_{lc}\overset{\frac{1}{2}}{\mu}{}^{-1a}_{b} \left(\begin{array}{cc|c} l & c & D \\ A & \frac{1}{2} & k \end{array}\right) \left(\begin{array}{c|cc} k & \frac{1}{2} & B \\ D & a & j \end{array}\right) \Lambda^{\frac{1}{2}D}_{AB}\, \overset{A}{e}_i =$$

$$= \sum_D (\frac{v_A v_{\frac{1}{2}}}{v_D})^{\pm\frac{1}{2}} \overset{\frac{1}{2}}{\mu}{}^{-1a}_{b} \left(\begin{array}{cc|c} b & i & D \\ \frac{1}{2} & A & k \end{array}\right) \left(\begin{array}{c|cc} k & \frac{1}{2} & B \\ D & a & j \end{array}\right) \Lambda^{\frac{1}{2}D}_{AB}\, \overset{A}{e}_i =$$

$$= \sum_D (\frac{v_A v_{\frac{1}{2}}}{v_D})^{\pm\frac{1}{2}} \Lambda_{AB}^{\frac{1}{2}D} e^{2i\pi\frac{1}{2}}$$

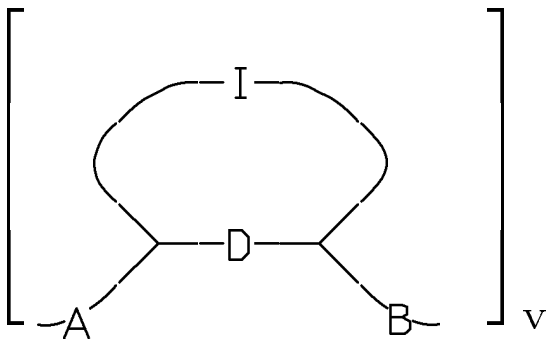

$$\overset{A}{e}_i =$$

$$= \sum_D (\frac{v_A v_{\frac{1}{2}}}{v_D})^{\pm\frac{1}{2}} \Lambda_{AB}^{\frac{1}{2}D} \frac{[2D+1]}{[2A+1]}$$

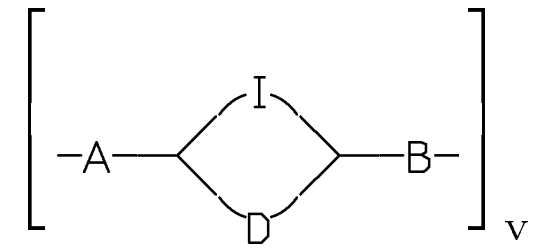

$$\overset{A}{e}_i =$$

$$= \delta^i_j \delta_{A,B} \sum_D (\frac{v_A v_{\frac{1}{2}}}{v_D})^{\pm\frac{1}{2}} \Lambda_{AB}^{\frac{1}{2}D} \frac{[2D+1]}{[2A+1]} \overset{A}{e}_i$$

which are the equations (111). As a first conclusion, we obtain that unitary irreducible representations of $\mathfrak{U}_q(sl(2,\mathbb{C})_{\mathbb{R}})$ give solutions of the relations (110,111) on the set of complex numbers $\Lambda_{AD}^{BC}$, and inversely solutions of these equations give a representation of $\mathfrak{U}_q(sl(2,\mathbb{C})_{\mathbb{R}})$ (the irreducibility and the unitarity has to be shown).
Let $\Lambda_{AD}^{BC}$ be given by the expression (108) which can be represented as:

$$\Lambda_{AD}^{BC} = \sum_{\aleph_2} \frac{v_C^{\frac{1}{2}}}{v_B^{\frac{1}{2}} v_A^{\frac{1}{4}} v_D^{\frac{1}{4}}}$$

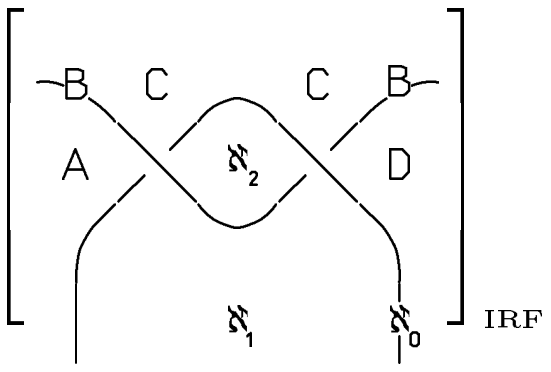

Let us first show the equality (109). The square of the expressions (108) and (109) are, from the expression (176) rational functions in $q^{i\frac{\rho}{2}}$. In order to show that these expressions are equal it is therefore sufficient to show that they are equal for an infinite set of values of $z = i\rho \in \frac{1}{2}\mathbb{Z}^+$. But this is the case because we have (92), and the complex continuation $\left\{\begin{matrix} B & C \\ \aleph_0 & \aleph_1 \end{matrix}\middle|\begin{matrix} D \\ \aleph_2 \end{matrix}\right\}$ is equal to the ordinary $6j$ coefficient when $z$ is an integer sufficiently large.
We want now to prove that the relation (110) is satisfied. This relation is also equivalent to the following form, using unitarity relations on $6j$ coefficients:

$$\Lambda_{FG}^{AC}\Lambda_{GH}^{BD} = \sum_{KU} e^{i\pi(G+U-C-D)} (\frac{[2U+1][2G+1]}{[2C+1][2D+1]})^{\frac{1}{2}} \left\{\begin{matrix} B & C \\ A & D \end{matrix}\middle|\begin{matrix} U \\ G \end{matrix}\right\} \left\{\begin{matrix} F & A \\ B & U \end{matrix}\middle|\begin{matrix} C \\ K \end{matrix}\right\} \left\{\begin{matrix} A & B \\ H & U \end{matrix}\middle|\begin{matrix} K \\ D \end{matrix}\right\} \Lambda_{FH}^{KU}. \tag{112}$$

This can also be written as :

$$e^{i\pi(2G-C-D)} \frac{[2C+1]^{\frac{1}{2}}[2D+1]^{\frac{1}{2}}}{[2G+1]} \Lambda_{FG}^{AC}\Lambda_{GH}^{BD} =$$
$$= \sum_{KU} e^{i\pi(G-U)} (\frac{[2U+1]}{[2G+1]})^{\frac{1}{2}} \left\{\begin{matrix} B & C \\ A & D \end{matrix}\middle|\begin{matrix} U \\ G \end{matrix}\right\} \left\{\begin{matrix} F & A \\ B & U \end{matrix}\middle|\begin{matrix} C \\ K \end{matrix}\right\} \left\{\begin{matrix} A & B \\ H & U \end{matrix}\middle|\begin{matrix} K \\ D \end{matrix}\right\} \Lambda_{FH}^{KU}. \tag{113}$$

Using the pictorial representation, it is easy to show that this relation is true. Indeed, the righthandside of the relation (113) can also be represented as

$$\sum_{KMU\aleph_2}\left(\frac{v_G v_F}{v_A v_B v_C v_D}\right)^{\frac{1}{2}}$$

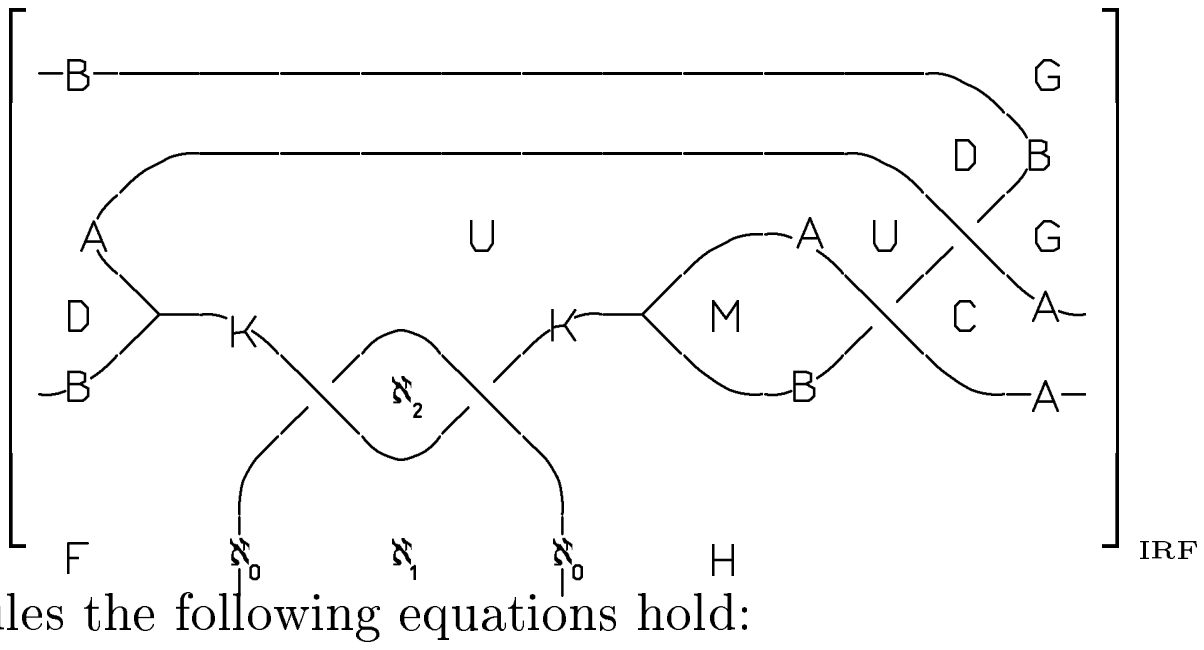

and, using the usual graphical rules the following equations hold:

$$\sum_{KMU\aleph_2}$$

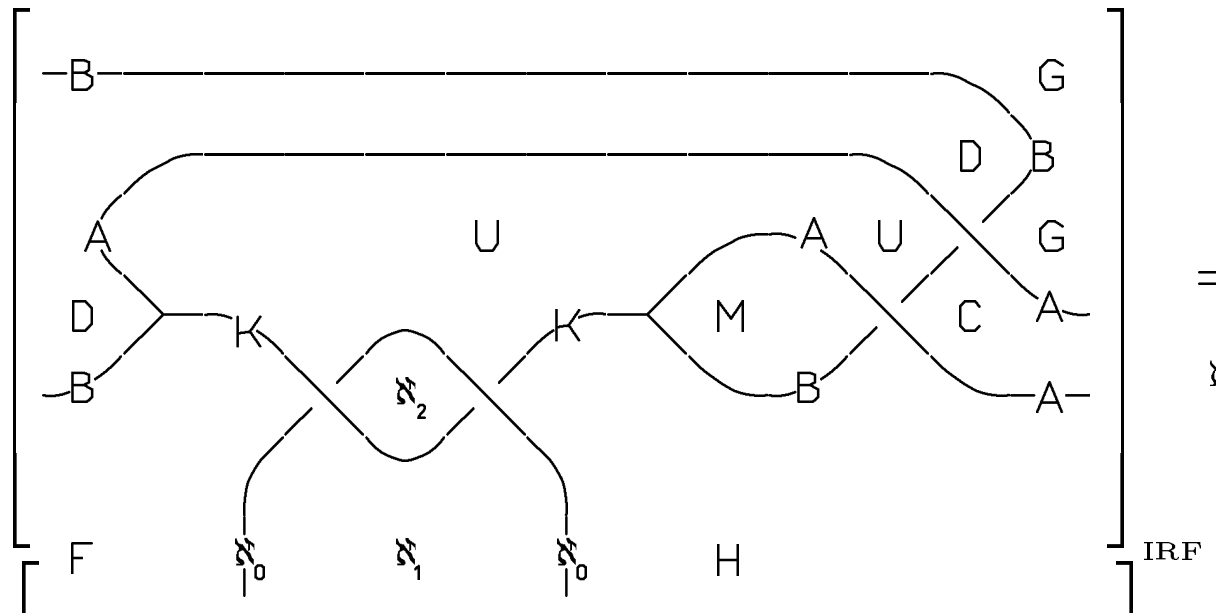

$$=\sum_{\substack{KMU\\ \aleph_2\aleph_3\aleph_4}}$$

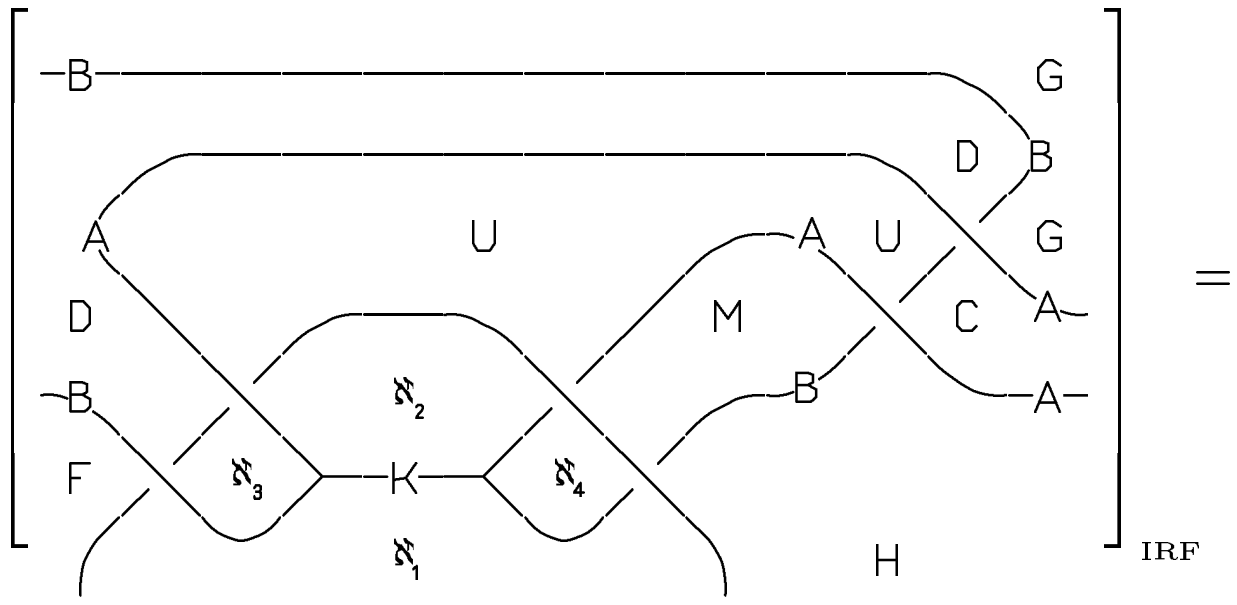

$$=$$

$$\sum_{MU\aleph_2\aleph_3}$$

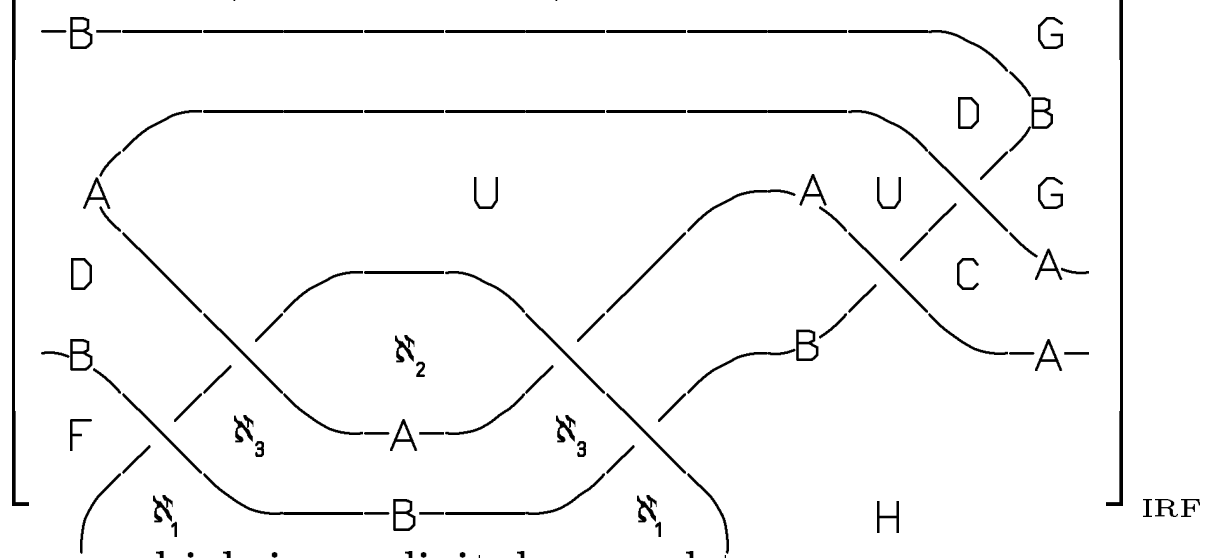

$$=\sum_{\aleph_2\aleph_3}$$

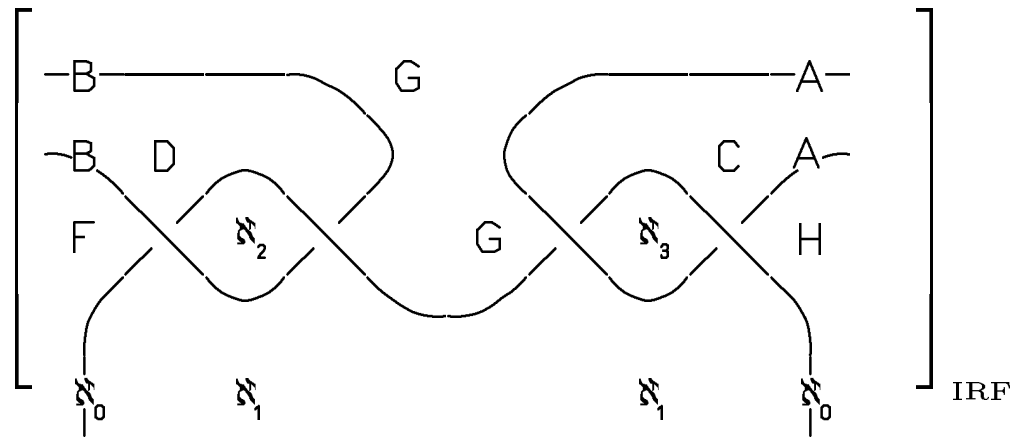

which is explicitely equal to

$$\sum_{\aleph_2\aleph_3} e^{i\pi(2G-C-D)}\frac{[2C+1]^{\frac{1}{2}}[2D+1]^{\frac{1}{2}}}{[2G+1]}\left\{\begin{array}{cc|c} B & D & F\\ \aleph_0 & \aleph_1 & \aleph_2\end{array}\right\}\left\{\begin{array}{cc|c} B & D & G\\ \aleph_0 & \aleph_1 & \aleph_2\end{array}\right\}\left\{\begin{array}{cc|c} A & C & G\\ \aleph_0 & \aleph_1 & \aleph_3\end{array}\right\}\left\{\begin{array}{cc|c} A & C & H\\ \aleph_0 & \aleph_1 & \aleph_3\end{array}\right\}\frac{v_{\aleph_2}v_{\aleph_3}}{v_{\aleph_1}^2}\left(\frac{v_F^{\frac{1}{2}}v_H^{\frac{1}{2}}v_G}{v_C v_D}\right).$$

which conclude the proof of the relation (113).

In order to prove the irreducibility and the unitarity of this representation we need an explicit formula for the action of $\overset{\frac{1}{2}}{g}{}^i_j$, i.e we need to compute $\Lambda^{\frac{1}{2}C}_{AD}$. This is trivial to do using the explicit formula for the $6j$ symbols when one of the spin is fixed to $\frac{1}{2}$, which are given in the appendix. We obtain the following formulas:

$$\begin{aligned}
\Lambda^{\frac{1}{2}C}_{C+\frac{1}{2}C+\frac{1}{2}}(\aleph_0,\aleph_1) &= -\frac{q^{2\aleph_1+1+\frac{m}{2}}[C+\frac{1}{2}-\frac{m}{2}]+q^{-2\aleph_1-1-\frac{m}{2}}[C+\frac{1}{2}+\frac{m}{2}]}{[2C+1]}\\
\Lambda^{\frac{1}{2}C}_{C-\frac{1}{2}C-\frac{1}{2}}(\aleph_0,\aleph_1) &= \frac{q^{2\aleph_1+1+\frac{m}{2}}[C+\frac{1}{2}+\frac{m}{2}]+q^{-2\aleph_1-1-\frac{m}{2}}[C+\frac{1}{2}-\frac{m}{2}]}{[2C+1]}\\
\Lambda^{\frac{1}{2}C}_{C-\frac{1}{2}C+\frac{1}{2}}(\aleph_0,\aleph_1) &= -i\frac{q-q^{-1}}{[2C+1]}\sqrt{[C+\frac{1}{2}-\frac{m}{2}][C+\frac{1}{2}+\frac{m}{2}][2\aleph_1+C+\frac{m}{2}+\frac{3}{2}][2\aleph_1-C+\frac{m}{2}+\frac{1}{2}]}\\
&= \Lambda^{\frac{1}{2}C}_{C+\frac{1}{2}C-\frac{1}{2}}(\aleph_0,\aleph_1).
\end{aligned}$$

The irreducibility of the representation can be easily proved. Indeed, the action of $\mathfrak{U}_q(su(2))$, on each vector space generated by $\{\overset{A}{e}_i, i = 1\cdots dim\overset{A}{V}\}$ is irreducible, any operator commuting with the action of the whole algebra must be scalar in each of these vector spaces . Let us denote by $\lambda_A$ its value in the vector space labelled by $A$. Then by analysing the action of $\overset{\frac{1}{2}}{g}{}^i_j$ we can easily see that, if this operator commutes with the generators $\overset{\frac{1}{2}}{g}{}^i_j$, then we have $(\lambda_A - \lambda_D)\Lambda^{\frac{1}{2}C}_{AD} = 0$. From the computation of $\Lambda^{\frac{1}{2}C}_{AD}$, we deduce that $\lambda_A = \lambda_{A+1}$, which allows us to conclude that this representation is irreducible. It remains to show that there exists an Hermitian form such that $\pi$ is unitary. Due to the fact that $\overset{A}{e}_i$ is an orthonormal basis for the representation $\overset{A}{\pi}$, it is sufficient to show that elements of $\mathfrak{U}_q(an(2))$ are represented by unitary actions. We have

$$(\overset{B}{e}_j \mid (\overset{C}{g}{}^k_l)^\star \mid \overset{A}{e}_i) = (\overset{B}{e}_j \mid S^{-1}(\overset{C}{g}{}^l_k) \mid \overset{A}{e}_i)$$
$$= \sum_D \overset{C}{w}_{km}\overset{C}{w}{}^{-1\,nl} \left(\begin{array}{cc|c} j & m & D \\ B & C & p \end{array}\right) \left(\begin{array}{c|cc} p & C & A \\ D & n & i \end{array}\right) \Lambda^{CD}_{AB}(\aleph_0,\aleph_1)$$
$$= \sum_D \left(\begin{array}{c|cc} j & D & C \\ B & p & k \end{array}\right) \frac{e^{i\pi(B-D)}[2D+1]^{\frac{1}{2}}}{[2B+1]^{\frac{1}{2}}} \left(\begin{array}{cc|c} l & p & A \\ C & D & i \end{array}\right) \frac{e^{i\pi(A-D)}[2D+1]^{\frac{1}{2}}}{[2A+1]^{\frac{1}{2}}} \Lambda^{CD}_{AB}(\aleph_0,\aleph_1)$$
$$= \sum_{DM} \left(\begin{array}{cc|c} l & j & M \\ C & B & q \end{array}\right) \left(\begin{array}{c|cc} q & A & C \\ M & i & k \end{array}\right) \frac{e^{i\pi(B+A-2D)}[2D+1]}{[2B+1]^{\frac{1}{2}}[2A+1]^{\frac{1}{2}}} \left\{\begin{array}{cc|c} C & D & B \\ C & M & A \end{array}\right\} \Lambda^{CD}_{AB}(\aleph_0,\aleph_1).$$

But we also have

$$\sum_D \frac{e^{i\pi(B+A-2D)}[2D+1]}{[2B+1]^{\frac{1}{2}}[2A+1]^{\frac{1}{2}}} \left\{\begin{array}{cc|c} C & D & B \\ C & M & A \end{array}\right\} \Lambda^{CD}_{AB}(\aleph_0,\aleph_1) = \sum_{D\aleph_2} \left(\frac{v_A^{\frac{1}{2}} v_B^{\frac{1}{2}}}{v_C v_M}\right)^{\frac{1}{2}}$$

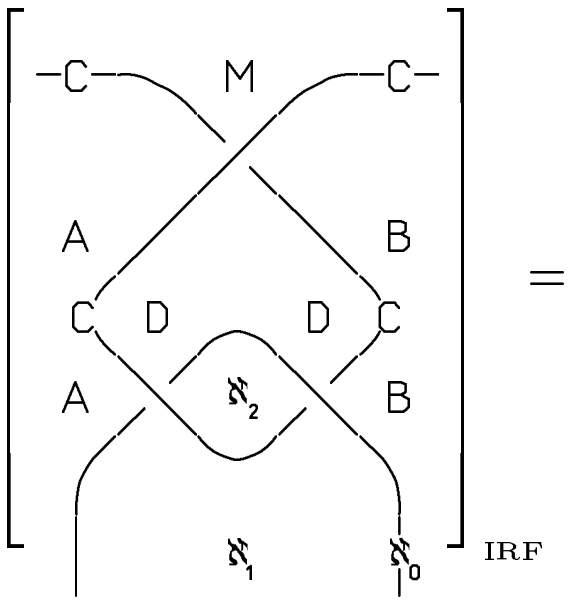

$$=$$

$$= \sum_{\aleph_2\aleph_3}$$

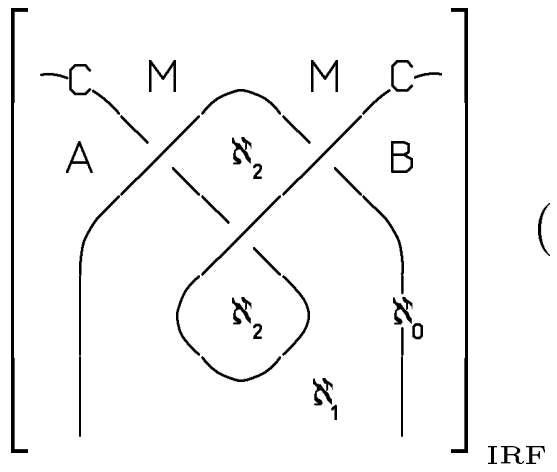

$$\left(\frac{v_A^{\frac{1}{2}} v_B^{\frac{1}{2}}}{v_C v_M}\right)^{\frac{1}{2}} = \sum_{\aleph_2}$$

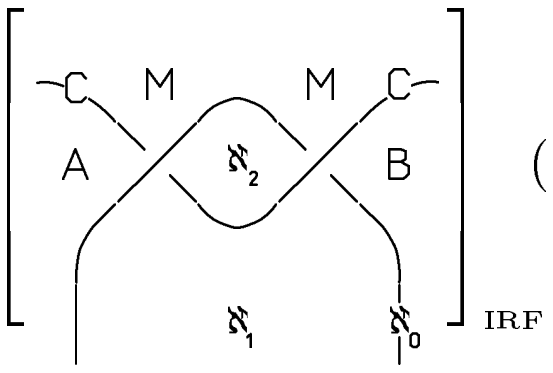

$$\left(\frac{v_A^{\frac{1}{2}} v_B^{\frac{1}{2}} v_C}{v_M}\right)^{\frac{1}{2}} = \Lambda^{CM}_{AB}(\aleph_1,\aleph_0)$$

where the last equality comes from equation (109). As a first conclusion, we then have

$$(\overset{B}{e}_j \mid (\overset{C}{g}{}^k_l)^\star \mid \overset{A}{e}_i) = \sum_M \left(\begin{array}{cc|c} l & j & M \\ C & B & q \end{array}\right) \left(\begin{array}{c|cc} q & A & C \\ M & i & k \end{array}\right) \Lambda^{CM}_{AB}(\aleph_1,\aleph_0).$$

On another hand we trivially get:

$$\overline{(\overset{A}{e}_i \mid \overset{C}{g}{}^k_l \mid \overset{B}{e}_j)} = \sum_M \left(\begin{array}{cc|c} l & j & M \\ C & B & q \end{array}\right) \left(\begin{array}{c|cc} q & A & C \\ M & i & k \end{array}\right) \overline{\Lambda^{CM}_{AB}(\aleph_0,\aleph_1)}$$

We thus have shown that unitary of the representation is equivalent to the relation: $\Lambda^{BC}_{AD}(\aleph_1,\aleph_0) = \overline{\Lambda^{BC}_{AD}(\aleph_0,\aleph_1)}$. Because $\overset{\frac{1}{2}}{g}{}^k_l$ generate $\mathfrak{U}_q(an(2))$, we only have to show that this last relation is satisfied for $B = \frac{1}{2}$. This is indeed true from the explicit form of $\Lambda^{BC}_{AD}(\aleph_0,\aleph_1)$ with $B = \frac{1}{2}$, computed above. In particular the value of $\Lambda^{\frac{1}{2}C}_{C-\frac{1}{2}C+\frac{1}{2}}(\aleph_0,\aleph_1)$ uses a square root and implies the constraint $\rho\in]0,1[$ for the two complementary series. □

From the definition of the complex continuation of the $6j$ symbols we would expect that $\Lambda^{BC}_{AA}$ is a rational function in $x$ where $x=q^{2\aleph_1+1}$. We will show that in fact $\Lambda^{BC}_{AA}$ as well as $(\Lambda^{BC}_{AD})^2$ are Laurent polynomials in the variable $x$ (this is clear for $B=\frac{1}{2}$ from the explicit expression). The reason for such a property is explained in the following section.

## 4.2 Connections between the coefficients $\Lambda^{BC}_{AD}$ and the universal shifted cocycle $F(x)$

In the article [7], O.Babelon introduces the universal matrix $F(x)\in\mathfrak{U}_q(sl(2))^{\otimes 2}$ in the context of quantization of Liouville theory on a lattice:

$$F_{12}(x)=\sum_{k=0}^{+\infty}(q-q^{-1})^k\frac{(-1)^k}{[k]!}\frac{x^k}{\prod_{\nu=k}^{2k-1}(xq^\nu q^{H_2}-x^{-1}q^{-\nu}q^{-H_2})}q^{k(H_1+H_2)/2}J_+^k\otimes J_-^k, \tag{114}$$

where we have denoted $H=2J_z$.

It is shown that

$$F_{12}(x)^{-1}=\sum_{k=0}^{+\infty}(q-q^{-1})^k\frac{1}{[k]!}\frac{x^k}{\prod_{\nu=1}^{k}(xq^\nu q^{H_2}-x^{-1}q^{-\nu}q^{-H_2})}q^{k(H_1+H_2)/2}J_+^k\otimes J_-^k \tag{115}$$

and that it satisfies the shifted cocycle identity:

$$(\mathrm{id}\otimes\Delta)(F(x))F_{23}(x)=(\Delta\otimes\mathrm{id})(F(x))F_{12}(xq^{H_3}). \tag{116}$$

In [8] a quasi-Hopf algebra interpretation of this construction is given as well as a connection between complex continuation of 6j coefficients and $F(x)$, which is encoded by the following formula:

$$\overset{AB}{F}(x)^{\sigma_1\sigma_2}_{\sigma_1'\sigma_2'}=\sum_C\frac{\mathcal{N}^{(C)}(x,\sigma_1+\sigma_2)\delta_{\sigma_1+\sigma_2,\sigma_1'+\sigma_2'}}{\mathcal{N}^{(A)}(xq^{2\sigma_2'},\sigma_1')\mathcal{N}^{(B)}(x,\sigma_2')}\left(\begin{matrix}\sigma_1 & \sigma_2\\ A & B\end{matrix}\middle|\begin{matrix}C\\ \sigma_1+\sigma_2\end{matrix}\right)\left\{\begin{matrix}A & B\\ \aleph(x) & \aleph(x)+\sigma_1'+\sigma_2'\end{matrix}\middle|\begin{matrix}C\\ \aleph(x)+\sigma_2'\end{matrix}\right\}, \tag{117}$$

with $x=q^{2\aleph(x)+1}$ and $\mathcal{N}^{(A)}(x,\sigma_1)$ are normalisation factors, the exact expression for them are derived in the appendix of [8]. From this result we easily obtain that:

$$(\overset{AB}{F}(x)^{-1})^{\sigma_1'\sigma_2'}_{\sigma_1\sigma_2}=\sum_C\frac{\mathcal{N}^{(A)}(xq^{2\sigma_2'},\sigma_1')\mathcal{N}^{(B)}(x,\sigma_2')\delta_{\sigma_1+\sigma_2,\sigma_1'+\sigma_2'}}{\mathcal{N}^{(C)}(x,\sigma_1+\sigma_2)}\left(\begin{matrix}\sigma_1+\sigma_2\\ C\end{matrix}\middle|\begin{matrix}A & B\\ \sigma_1 & \sigma_2\end{matrix}\right)\left\{\begin{matrix}A & B\\ \aleph(x) & \aleph(x)+\sigma_1'+\sigma_2'\end{matrix}\middle|\begin{matrix}C\\ \aleph(x)+\sigma_2'\end{matrix}\right\}. \tag{118}$$

We first have the following simple result which relates $\Lambda^{BC}_{AD}(\aleph_0,\aleph_1)$ to $F(x)$ :

Lemma 4:

*The coefficients $\Lambda^{BC}_{AD}(\aleph_0,\aleph_1)$ satisfy the relation:*

$$\Lambda^{BC}_{AD}(\aleph_0,\aleph_1)=\left(\begin{matrix}\frac{m}{2}\\ A\end{matrix}\middle|\begin{matrix}C & B\\ \sigma_2 & \sigma_1\end{matrix}\right)\left(\frac{v_A}{v_D}\right)^{\frac{1}{4}}\frac{\mathcal{N}^{(D)}(x,\frac{m}{2})}{\mathcal{N}^{(A)}(x,\frac{m}{2})}(\Lambda(x))^{\sigma_1\sigma_2}_{\sigma_1'\sigma_2'}\left(\begin{matrix}\sigma_1' & \sigma_2'\\ B & C\end{matrix}\middle|\begin{matrix}D\\ \frac{m}{2}\end{matrix}\right) \tag{119}$$

*where as usual $\aleph_0-\aleph_1=\frac{m}{2}$, and $x=q^{2\aleph_1+1}$, and we have denoted $\Lambda(x)=RF(x)G(x)F(x)^{-1}$, with $G(x)\in\mathfrak{U}_q(sl(2))^{\otimes 2}$ the element $G(x)=q^{-\frac{(H_1+H_2)^2}{2}}x^{-H_1}q^{\frac{H_2^2}{2}}$.*

Proof: Using the relations connecting $F(x)$ to the complex 6j coefficients, we can rewrite equation (109) as:

$$\Lambda^{BC}_{AD}(\aleph_0,\aleph_1)=(\frac{v_A}{v_D})^{\frac{1}{4}}\frac{\mathcal{N}^{(D)}(x,\frac{m}{2})}{\mathcal{N}^{(A)}(x,\frac{m}{2})}\times$$
$$\times\sum(\frac{v_Bv_C}{v_A})^{\frac{1}{2}}\left(\begin{matrix}\frac{m}{2}\\ A\end{matrix}\middle|\begin{matrix}B & C\\ \sigma_1 & \sigma_2\end{matrix}\right)\overset{BC}{F}(x)^{\sigma_1\sigma_2}_{\sigma_1'\sigma_2'}\frac{v_{\aleph_1+\frac{m}{2}}}{v_{\aleph_1+\sigma_2'}}(\overset{BC}{F}(x)^{-1})^{\sigma_1'\sigma_2'}_{\sigma_1''\sigma_2''}\left(\begin{matrix}\sigma_1'' & \sigma_2''\\ B & C\end{matrix}\middle|\begin{matrix}D\\ \frac{m}{2}\end{matrix}\right). \quad (120)$$

In the last sum we necessarily have $\sigma_1+\sigma_2=\frac{m}{2}$, it is then natural to introduce the matrix $G^{\sigma_1,\sigma_2}_{\sigma_1',\sigma_2'}=\frac{v_{\aleph_1+\sigma_1+\sigma_2}}{v_{\aleph_1+\sigma_2'}}\delta_{\sigma_1,\sigma_1'}\delta_{\sigma_2,\sigma_2'}$, which is the representation in the tensor product $\overset{B}{\pi}\otimes\overset{C}{\pi}$ of the element $G(x)$ introduced in the proposition. Using the twist relation on the Clebsch-Gordan coefficients we finally obtain the relation of the proposition. □

Let us denote by $B(x)$ the element $B(x)=q^{\frac{H_2^2}{2}}x^{H_2}$, we have the following result:

**Lemma 5:**

*For each $p\in\mathbb{N}$,*

$$F(x)B(x)^pF(x)^{-1}=B(x)^p\sum_{n=0}^{+\infty}\frac{(q-q^{-1})^n}{[n]!}q^{\frac{n}{2}(H_1+H_2)}x^n g_{p,n}(xq^{H_2})J_+^n\otimes J_-^n, \quad (121)$$

*where $g_{p,n}(x)\in\mathbb{C}[x]$, and in particular $g_{0,n}(x)=\delta_{0,n}$ and $g_{1,n}(x)=(-1)^nq^{\frac{n(n+1)}{2}}x^n$. From this result we obtain the following identity on $F(x)$:*

$$F_{12}(x)q^{\frac{H_2^2}{2}}x^{H_2}F_{12}(x)^{-1}=R^{-1}q^{\frac{H\otimes H}{2}}q^{\frac{H_2^2}{2}}x^{H_2} \quad (122)$$

Proof: By applying the proof to the case where $p=0$ we will obtain a simpler proof of the fact that $F(x)^{-1}$ is given by equation (115). Let us introduce the rational functions $f_k(x),\tilde{f}_k(x)$ defined by:

$$f_k(x)=\frac{1}{\prod_{\nu=k}^{2k-1}(xq^\nu-x^{-1}q^{-\nu})},\tilde{f}_k(x)=\frac{1}{\prod_{\nu=1}^{k}(xq^\nu-x^{-1}q^{-\nu})},$$

from the relation $J_\pm^k\psi(H)=\psi(H\mp 2k)J_\pm^k$ we obtain that relation (121) holds true with

$$g_p=\sum_{k=0}^{n}(-1)^k\binom{n}{k}_q f_k(x)\tilde{f}_{n-k}(xq^{2k})q^{2pk^2}x^{2pk},$$

(we will omit the variable $n$ in the sequel). Let us prove that $g_p(x)$ is a polynomial in $x$. Let $R_{k,p}(x)=f_k(x)\tilde{f}_{n-k}(xq^{2k})x^{2pk}$, this rational function is equal to

$$R_{k,p}(x)=\frac{(x^2q^{2k}-1)x^{n+2pk}}{\prod_{\nu=k}^{k+n}(x^2q^\nu-q^{-\nu})},$$

and we can decompose this rational function in simple poles:

$$R_{k,p}(x) = P_{k,p}(x) + \sum_{l=0}^{n} \frac{A_{l,k,p}x + B_{l,k,p}}{x^2 q^{k+l} - q^{-k-l}}$$

with $P_{k,p}(x)$ polynomial of degree less or equal to $2pk - n$. Because $R_{k,p}(x)$ has the same parity as $n$, we necessarely have $A_{l,k,p} = 0$ if $n$ is even, and $B_{l,k,p} = 0$ if $n$ is odd. When $n$ is odd we get $A_{l,k,p} = (-1)^l (q - q^{-1})^{1-n} [n]! \binom{n}{l}_q [k-l] q^{-2pk(k+l)}$, therefore

$$g_p(x) = \sum_{k=0}^{n} P_{k,p}(x) + \sum_{k,l=0}^{n} \frac{(-1)^{k+l}(q-q^{-1})^{1-n}}{xq^{k+l} - x^{-1}q^{-k-l}} \binom{n}{k}_q \binom{n}{l}_q q^{-2pkl}[k-l].$$

The second term of the right handside is zero due to the antisymmetry in the exchange of $k$ and $l$. As a result $g_p(x)$ is a polynomial of degree less or equal to $(2p-1)n$. From $R_{k,p}(x) = x^{2pk+n} Q_{k,p}(x)$, with $Q_{k,p}(x)$ rational function not having 0 as pole, we get by differentiating that $R_{k,p}(x)^{(m)}(0) = 0$, for $0 \leq m < n$, and $R_{0,p}(x)^{(n)} = n!(-1)^n q^{\frac{n(n+1)}{2}}$. The proof is similar if $n$ is even. As a result we obtain that $g_{0,n}(x) = \delta_{0,n}$ and $g_{1,n}(x) = (-1)^n q^{\frac{n(n+1)}{2}} x^n$. An immediate application of this result shows that

$$F(x)B(x)F(x)^{-1} = \{\sum_{n=0}^{+\infty} \frac{(q-q^{-1})^n}{[n]!} (-1)^n q^{\frac{n}{2}(H_1 - H_2)} q^{\frac{n(n+1)}{2} - 2n^2} J_+^n \otimes J_-^n\} B(x) = R^{-1} q^{\frac{H \otimes H}{2}} B(x).$$

□

We summarize the previous results in the following proposition, which gives a new expression for the coefficients $\Lambda_{AD}^{BC}(\aleph_0, \aleph_1)$.

**Proposition 1:**

$$\Lambda_{AD}^{BC}(\aleph_0, \aleph_1) = (\frac{v_A}{v_D})^{\frac{1}{4}} \frac{\mathcal{N}^{(D)}(x, \frac{m}{2})}{\mathcal{N}^{(A)}(x, \frac{m}{2})} \left( \begin{array}{c|cc} \frac{m}{2} & C & B \\ A & \sigma_2 & \sigma_1 \end{array} \right) q^{-2\sigma_1 z} \left( \begin{array}{cc|c} \sigma_1 & \sigma_2 & D \\ B & C & \frac{m}{2} \end{array} \right) \tag{123}$$

*where again* $2\aleph_0 + 1 = \frac{m}{2} + z, 2\aleph_1 + 1 = -\frac{m}{2} + z$. *As a result* $\Lambda_{AA}^{BC}$ *as well as* $(\Lambda_{AD}^{BC})^2$ *are Laurent polynomials in the variable* $x$.

Proof: The universal matrix $\Lambda(x)$ can be simplified using the previous Lemma, and we have:

$$\Lambda(x) = q^{-\frac{H \otimes H}{2}} q^{-\frac{H_1^2}{2}} x^{-H_1}.$$

Using this simple expression, it is then trivial to obtain the announced formula. □

Remark: The linear relation $F_{12}(x) q^{\frac{H_2^2}{2}} x^{H_2} = R^{-1} q^{\frac{H \otimes H}{2}} q^{\frac{H_2^2}{2}} x^{H_2} F_{12}(x)$ is a very important relation. We have computed [3], see also [24], using a generalisation of this identity, the exact form of the shifted cocycles $F : \mathfrak{H}^\star \to \mathfrak{U}_q(\mathfrak{G})$, where $\mathfrak{H}$ is a Cartan subalgebra of $\mathfrak{G}$.

## 4.3 Characters

$Fun_{cc}(SL_q(2,\mathbb{C})_{\mathbb{R}}) = Pol(SU_q'(2) \otimes Fun_c(AN_q(2))$ is a multiplier Hopf algebra, let us define on the vector space $\tilde{\mathcal{D}} = Pol(SU_q(2))^* \otimes Pol(SU_q(2))$, where $Pol(SU_q(2))^*$ is the restricted dual of $Pol(SU_q(2))$ and $\otimes$ is the algebraic tensor product, a structure of multiplier Hopf algebra in duality with $Fun_{cc}(SL_q(2,\mathbb{C})_{\mathbb{R}})$. Let $\overset{A}{X}{}^i_j \in Pol(SU_q(2))^*$ be the dual basis of $\overset{A}{g}{}^i_j$, i.e $<\overset{A}{X}{}^i_j, \overset{B}{g}{}^k_l>= \delta_{AB}\delta^i_l\delta^k_j$, the structure of algebra on $Pol(SU_q(2))^* \otimes Pol(SU_q(2))$, is defined by the quantum double, i.e

$$\overset{I}{g}_1\overset{J}{g}_2 = \sum_M \phi^{IJ}_M \overset{M}{g}\, \psi^M_{IJ}, \tag{124}$$

$$\overset{A}{X}{}^i_j\overset{B}{X}{}^k_l = \delta_{AB}\delta^i_l \overset{A}{X}{}^k_j \tag{125}$$

and a braided relation between $\overset{A}{X}{}^i_j$ and $\overset{B}{g}{}^k_l$ that can be computed from (16). By duality the coproduct on $\overset{A}{X}{}^i_j$ is defined by

$$\Delta(\overset{A}{X}{}^i_j) = \sum_{BC} \begin{pmatrix} m & r \\ B & C \end{pmatrix}\!\!\left|\begin{matrix} A \\ j\end{matrix}\right)\begin{pmatrix} i \\ A \end{pmatrix}\!\!\left|\begin{matrix} B & C \\ n & s\end{matrix}\right) \overset{B}{X}{}^n_m \otimes \overset{C}{X}{}^s_r. \tag{126}$$

The element $\overset{A}{L}{}^{\pm} \in \mathcal{D}$ are multipliers of $\tilde{\mathcal{D}}$ and we have the trivial relation $\overset{A}{L}{}^{(\pm)}{}^i_j = \sum^{\oplus}_B \overset{AB}{R}{}^{\pm ik}_{jl} \overset{B}{X}{}^l_k$.

From the expression of an irreducible unitary representation $\pi$ of $\mathcal{D} = \mathfrak{U}_q(sl(2,\mathbb{C})_{\mathbb{R}})$, associated to $\aleph_0, \aleph_1$ we can construct a representation $\tilde{\pi}$ of $\tilde{\mathcal{D}}$ as follows:

$$\begin{aligned} \tilde{\pi}(\overset{B}{X}{}^i_j)\overset{C}{e}_m &= \delta_{BC}\delta^i_m \overset{C}{e}_j \\ \tilde{\pi}(\overset{B}{g}{}^i_j)\overset{C}{e}_m &= \sum_{DE} \overset{E}{e}_p \begin{pmatrix} p & i \\ E & B \end{pmatrix}\!\!\left|\begin{matrix} D \\ x\end{matrix}\right) \begin{pmatrix} x \\ D \end{pmatrix}\!\!\left|\begin{matrix} B & C \\ j & m\end{matrix}\right) \Lambda^{BD}_{EC}(\aleph_0,\aleph_1). \end{aligned}$$

For each $\phi \in Fun_{cc}(SL_q(2,\mathbb{C})_{\mathbb{R}})$, let us define $\pi(\phi)$ to be the element of $End(V)$ defined by:

$$\pi(\phi) = \sum_{AB} \tilde{\pi}(\overset{A}{X}{}^i_j \otimes \overset{B}{g}{}^k_l) h((\overset{A}{k}{}^j_i \otimes \overset{B}{E}{}^l_k)\phi), \tag{127}$$

this sum involves just a finite number of non zero terms. From the structure of $\tilde{\pi}$ we immediately see that the matrix of $\pi(\phi)$ in the Hilbert basis $\{\overset{A}{e}_m, A \in \frac{1}{2}\mathbb{Z}^+, m = -A, ..., A\}$ has only a finite number of non zero matrix elements. As a result we can define $\chi_\pi(\phi) = tr_V(\pi(\mu^{-1})\pi(\phi))$, which defines the character associated to $\pi$.

We have a natural inclusion $Fun(SL_q(2,\mathbb{C})_{\mathbb{R}}) \overset{\iota}{\hookrightarrow} (Fun_{cc}(SL_q(2,\mathbb{C})_{\mathbb{R}}))^*$ given as usual by: if $f \in Fun(SL_q(2,\mathbb{C})_{\mathbb{R}})$ we have $\iota(f)(a) = h(fa), \forall a \in Fun_{cc}(SL_q(2,\mathbb{C})_{\mathbb{R}})$.

Let $\psi \in (Fun_{cc}(SL_q(2,\mathbb{C})_{\mathbb{R}}))^*$ where the $*$ denotes the algebraic dual, we can endow this vector space with a structure of topological space with the weak * topology. Let $\psi(\overset{I}{k}{}^m_n \otimes \overset{J}{E}{}^r_s) = \overset{IJ}{\psi}{}^{ns}_{mr}$, we always have $\psi = \sum_{IJ} \overset{IJ}{\psi}{}^{ns}_{mr} \iota(\overset{I}{k}{}^m_n \otimes \overset{J}{E}{}^r_s)$ where the convergence always holds.

As a result we have $\chi_\pi = \sum_{A,B,C} \begin{pmatrix} j & l \\ A & B \end{pmatrix}\!\!\left|\begin{matrix} C \\ n\end{matrix}\right) \begin{pmatrix} n \\ C \end{pmatrix}\!\!\left|\begin{matrix} B & A \\ m & t\end{matrix}\right) \Lambda^{BC}_{AA}\ \overset{A}{\mu}{}^{-1t}_{\ \ i}\ \iota(\overset{A}{k}{}^i_j \otimes \overset{B}{E}{}^m_l)$. We will use the notation $\chi_{AB}{}^{jl}_{im} = \sum_C \begin{pmatrix} j & l \\ A & B \end{pmatrix}\!\!\left|\begin{matrix} C \\ n\end{matrix}\right) \begin{pmatrix} n \\ C \end{pmatrix}\!\!\left|\begin{matrix} B & A \\ m & t\end{matrix}\right) \Lambda^{BC}_{AA}\ \overset{A}{\mu}{}^{-1t}_{\ \ i}$ and $\Lambda^{BC}_A = \Lambda^{BC}_{AA}$.

We have already defined $\overset{R}{\triangleright}$ and $\overset{L}{\triangleright}$ as right and left action of $\mathfrak{U}_q(sl(2,\mathbb{C})_{\mathbb{R}})$ on $Fun_{cc}(SL_q(2,\mathbb{C})_{\mathbb{R}})$, let us define a right and left action of $\mathfrak{U}_q(sl(2,\mathbb{C})_{\mathbb{R}})$ on $(Fun_{cc}(SL_q(2,\mathbb{C})_{\mathbb{R}}))^*$ : let $\psi \in (Fun_{cc}(SL_q(2,\mathbb{C})_{\mathbb{R}}))^*$ and $x \in \mathfrak{U}_q(sl(2,\mathbb{C})_{\mathbb{R}})$ we denote $x \overset{R}{\triangleright} \psi, (resp. x \overset{L}{\triangleright} \psi)$ to be the elements of $(Fun_{cc}(SL_q(2,\mathbb{C})_{\mathbb{R}}))^*$ defined as follow:

$$x \overset{L}{\triangleright} \psi(a) = \psi(S^{-1}(x) \overset{L}{\triangleright} a), \quad x \overset{R}{\triangleright} \psi(a) = \psi(S(x) \overset{R}{\triangleright} a), \quad \forall a \in Fun_{cc}(SL_q(2,\mathbb{C})_{\mathbb{R}}). \tag{128}$$

This definition is chosen in order that if $f \in Fun(SL_q(2,\mathbb{C})_{\mathbb{R}})$, $x \overset{R,L}{\triangleright} \iota(f) = \iota(x \overset{R,L}{\triangleright} f)$.

From these definitions we obtain the following proposition:

Proposition 2:

$\forall x \in \mathfrak{U}_q(sl(2,\mathbb{C})_{\mathbb{R}}$

$\pi(x \overset{L}{\triangleright} \phi) = \pi(S(x))\pi(\phi), \pi(x \overset{R}{\triangleright} \phi) = \pi(\phi)\pi(S^{-1}(x)), \forall \phi \in Fun_{cc}(SL_q(2,\mathbb{C})_{\mathbb{R}})$

$x \overset{R}{\triangleright} \chi_\pi = S^2(x) \overset{L}{\triangleright} \chi_\pi$, *($ad^+$ invariance)*

$\Omega^{\pm} \overset{L}{\triangleright} \chi_\pi = \omega^{\pm}\chi_\pi$, *where $\omega^{\pm}$ are the eigenvalues of the Casimir operator.*

Proof: Easy application of the definitions. □

We can define on $Fun_{cc}(SL_q(2,\mathbb{C})_{\mathbb{R}})$ a convolution product as follows: let $\phi, \psi \in Fun_{cc}(SL_q(2,\mathbb{C})_{\mathbb{R}})$, we define $\phi * \psi = \sum_{(\phi)} \phi_{(1)} h(S(\phi_{(2)})\psi)$.

It is easy to show that the following proposition is verified:

Proposition 3:

$$\forall \phi, \psi \in Fun_{cc}(SL_q(2,\mathbb{C})_{\mathbb{R}}), \qquad \pi(\phi \star \psi) = \pi(\phi)\pi(\psi) \tag{129}$$

$$\pi(S(\phi^\star)) = \pi(\phi)^\dagger. \tag{130}$$

Proof: This is a consequence of the fact that $\tilde{\pi}$ is a unitary representation of $\tilde{\mathcal{D}}$. □

$\chi_\pi$ being an $ad^+$ invariant element and an eigenvector of the action of $\Omega^{\pm}$, this implies linear relations satisfied by the elements $\Lambda^{BC}_{AA}$. We will need this linear system in the proof of Plancherel formula.

Proposition 4: Invariance under $\mathfrak{U}_q(an(2))$

*The element $\chi_\pi$ of $(Fun_{cc}(SL_q(2,\mathbb{C})_{\mathbb{R}}))^*$ is invariant under the adjoint action of $\mathfrak{U}_q(an(2))$, which implies that $\Lambda^{AB}_{CC}$ satisfy the linear equations:*

$$\sum_{AB} \Lambda^{BD}_{AA} \left\{ \begin{matrix} M & C \\ A & D \end{matrix} \middle| \begin{matrix} B \\ R \end{matrix} \right\} \left\{ \begin{matrix} C & A \\ M & K \end{matrix} \middle| \begin{matrix} R \\ T \end{matrix} \right\} \left\{ \begin{matrix} A & M \\ C & D \end{matrix} \middle| \begin{matrix} T \\ B \end{matrix} \right\} =$$

$$\sum_{AB} \Lambda^{BK}_{AA} \left\{ \begin{matrix} A & C \\ M & K \end{matrix} \middle| \begin{matrix} R \\ B \end{matrix} \right\} \left\{ \begin{matrix} C & M \\ A & K \end{matrix} \middle| \begin{matrix} B \\ T \end{matrix} \right\} \left\{ \begin{matrix} M & A \\ C & D \end{matrix} \middle| \begin{matrix} T \\ R \end{matrix} \right\} \tag{131}$$

*for all $K, T, C, M, D, R \in \frac{1}{2}\mathbb{Z}^+$.*

Proof: The right action of $\mathfrak{U}_q(an(2))$ on $\chi_\pi$ is given by the formula:

$$\overset{C}{g}{}^p_n \overset{R}{\triangleright} (\sum_{A,B} \chi_{AB}{}^{jl}_{im} \; \iota(\overset{A}{k}{}^i_j \otimes \overset{B}{E}{}^m_l)) = \sum_{ABMC} \chi_{AB}{}^{jl}_{im} \begin{pmatrix} u & d \\ C & M \end{pmatrix} \begin{matrix} B \\ l \end{pmatrix} \begin{pmatrix} m \\ B \end{matrix} \begin{pmatrix} C & M \\ v & b \end{pmatrix} \iota(S^{-1}(\overset{A}{k})^v_n \overset{C}{k}{}^i_j \overset{A}{k}{}^p_u \otimes \overset{M}{E}{}^b_d)$$

whereas the left action is (132)

$$\overset{C}{g}{}^p_n \overset{L}{\triangleright} (\sum_{A,B} \chi_{AB}{}^{jl}_{im} \; \iota(\overset{A}{k}{}^i_j \otimes \overset{B}{E}{}^m_l)) = \sum_{ABMC} \chi_{AB}{}^{jl}_{im} \left( \begin{matrix} c & p \\ M & C \end{matrix} \middle| \begin{matrix} B \\ l \end{matrix} \right) \left( \begin{matrix} m \\ B \end{matrix} \middle| \begin{matrix} M & C \\ a & n \end{matrix} \right) \iota(\overset{A}{k}{}^i_j \otimes \overset{M}{E}{}^a_c). \tag{133}$$

The element $\chi_\pi$ is therefore $ad^+$ invariant if the following relation holds:

$$\sum_{ABD} \left(\begin{array}{cc|c} j & l & D \\ A & B & x \end{array}\right) \left(\begin{array}{c|cc} x & B & A \\ D & m & i \end{array}\right) \left(\begin{array}{cc|c} c & p & B \\ M & C & l \end{array}\right) \left(\begin{array}{c|cc} m & M & C \\ B & a & n \end{array}\right) \left(\begin{array}{cc|c} n & i & R \\ C & A & e \end{array}\right) \left(\begin{array}{c|cc} f & C & A \\ R & t & j \end{array}\right) \Lambda^{BD}_{AA} =$$

$$= \sum_{ABD} \left(\begin{array}{cc|c} u & c & B \\ C & M & l \end{array}\right) \left(\begin{array}{c|cc} m & C & M \\ B & t & a \end{array}\right) \left(\begin{array}{cc|c} i & p & R \\ A & C & e \end{array}\right) \left(\begin{array}{c|cc} f & A & C \\ R & j & u \end{array}\right) \left(\begin{array}{cc|c} j & l & D \\ A & B & x \end{array}\right) \left(\begin{array}{c|cc} x & B & A \\ D & m & i \end{array}\right) \Lambda^{BD}_{AA}.$$

The left handside of this equation is nicely described by the following picture

$$\sum_{ABD} \Lambda^{BD}_{A}$$ 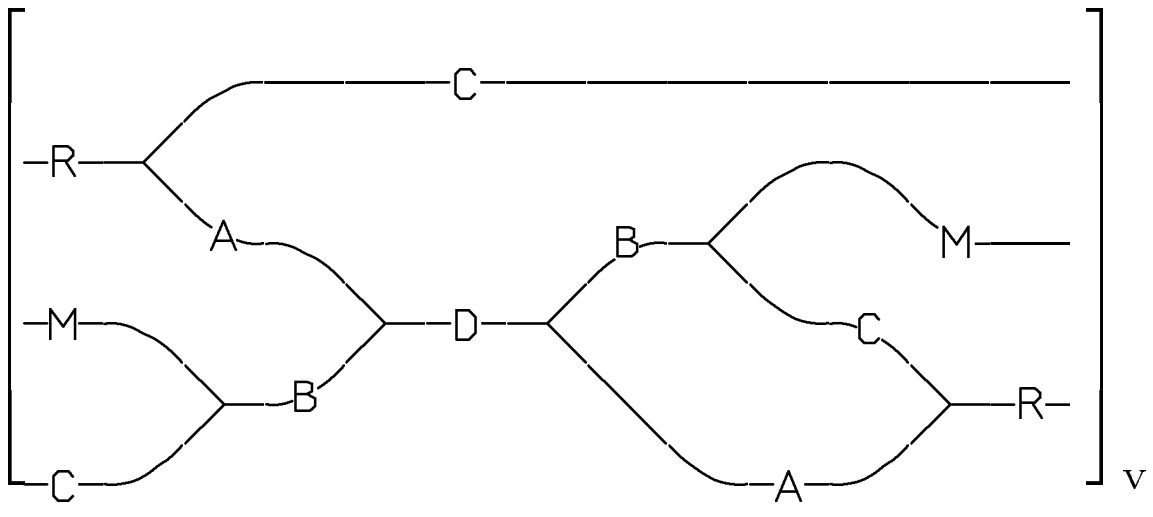 ,

which can be written using the definitions of $6j$ symbols and the unitarity relations on Clebsch-Gordan coefficients as

$$\sum_{ABDKT} \Lambda^{BD}_{A} \left\{\begin{array}{cc|c} M & C & B \\ A & D & R \end{array}\right\} \left\{\begin{array}{cc|c} C & A & R \\ M & K & T \end{array}\right\} \left\{\begin{array}{cc|c} A & M & T \\ C & D & B \end{array}\right\}$$ 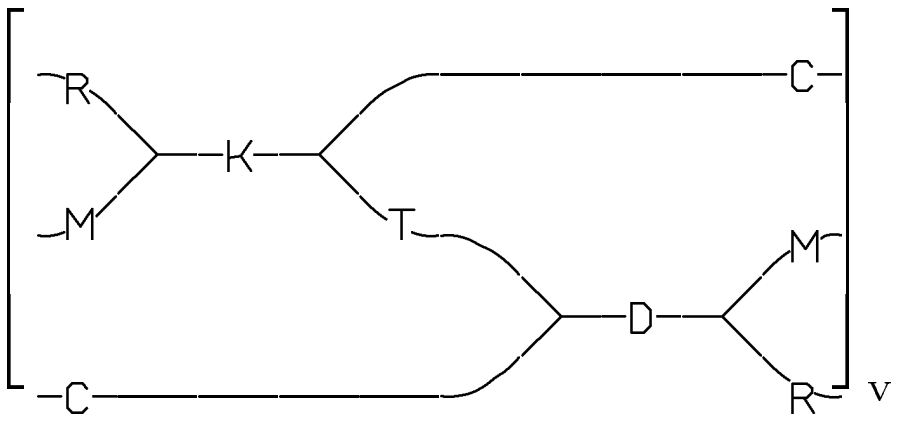

The righthandside is represented by the picture:

$$\sum_{ABD} \Lambda^{BD}_{A}$$ 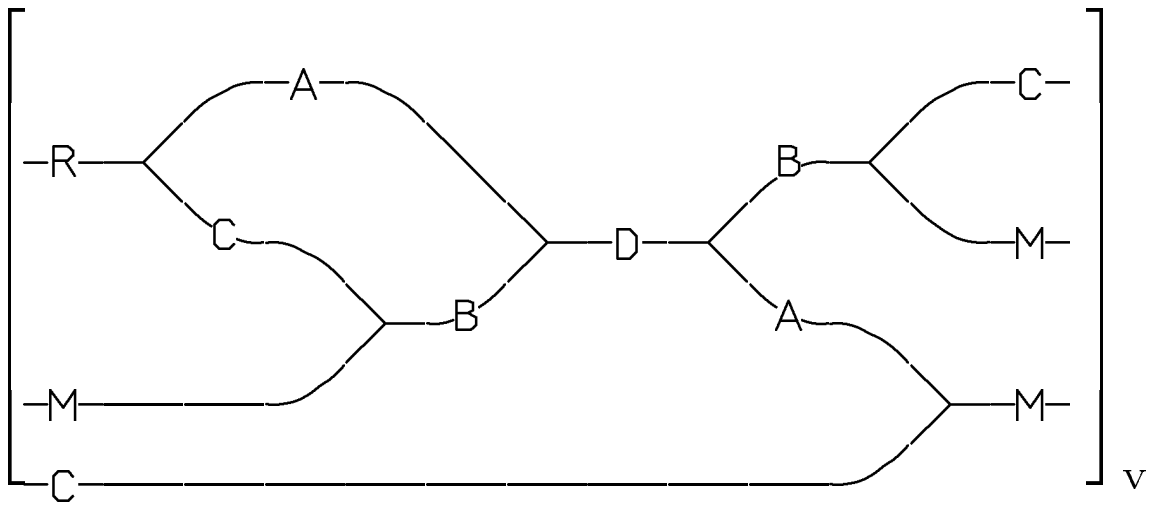 ,

which is equal after similar manipulations to

$$\sum_{ABDKT} \Lambda^{BK}_{A} \left\{\begin{array}{cc|c} A & C & R \\ M & K & B \end{array}\right\} \left\{\begin{array}{cc|c} C & M & B \\ A & K & T \end{array}\right\} \left\{\begin{array}{cc|c} M & A & T \\ C & D & R \end{array}\right\}$$ 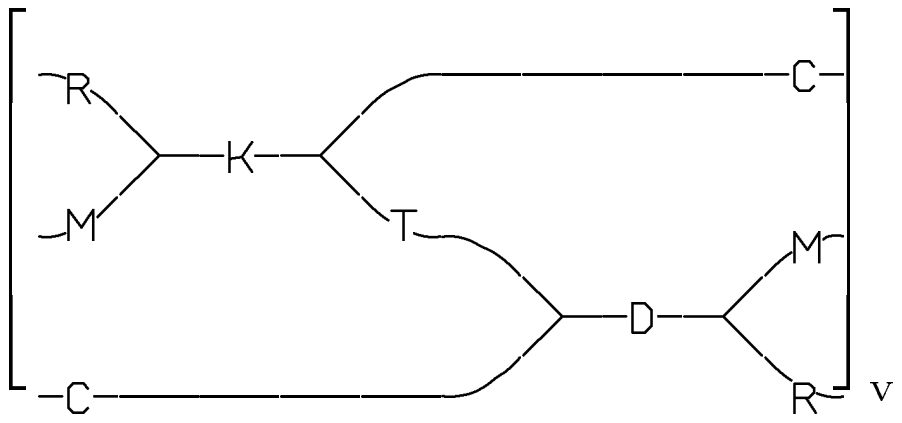

As a result we obtain the announced equation (131).

□

We will need the explicit form of this system of equation for $C = \frac{1}{2}$. Using the formulas for the 6-j symbols with one spin fixed to $\frac{1}{2}$ (175) the system is written as: $\forall (R, M, T) \in \frac{1}{2}\mathbb{Z}^+ \times \frac{1}{2}\mathbb{Z}^+ \times \frac{1}{2}\mathbb{Z}^{+*}$,

$$
\begin{aligned}
&Z(M,T,R)\left\{\frac{[R{+}M{+}T{+}\frac{5}{2}]}{[2R{+}2]}\Lambda^{M+\frac{1}{2}T+\frac{1}{2}}_{R+\frac{1}{2}}+\frac{[R{+}T{-}M{+}\frac{3}{2}]}{[2R{+}2]}\Lambda^{M-\frac{1}{2}T+\frac{1}{2}}_{R+\frac{1}{2}}-\frac{[T{+}M{-}R-\frac{1}{2}]}{[2R{+}2]}\Lambda^{M-\frac{1}{2}T-\frac{1}{2}}_{R+\frac{1}{2}}+\frac{[M{+}R{-}T{+}\frac{3}{2}]}{[2R{+}2]}\Lambda^{M+\frac{1}{2}T-\frac{1}{2}}_{R+\frac{1}{2}}\right\}=\\
&=Z(M,T,R)\left\{-\frac{[T{+}M{-}R{+}\frac{3}{2}]}{[2R]}\Lambda^{M+\frac{1}{2}T+\frac{1}{2}}_{R-\frac{1}{2}}+\frac{[R{+}M{-}T-\frac{1}{2}]}{[2R]}\Lambda^{M-\frac{1}{2}T+\frac{1}{2}}_{R-\frac{1}{2}}+\frac{[R{+}T{-}M-\frac{1}{2}]}{[2R]}\Lambda^{M+\frac{1}{2}T-\frac{1}{2}}_{R-\frac{1}{2}}+\frac{[R{+}M{+}T{+}\frac{1}{2}]}{[2R]}\Lambda^{M-\frac{1}{2}T-\frac{1}{2}}_{R-\frac{1}{2}}\right\}
\end{aligned}
\tag{134}
$$

where $\quad Z(M,R,T) = \dfrac{[R{+}T{+}M{+}\frac{3}{2}]^{\frac{1}{2}}[R{+}T{-}M{+}\frac{1}{2}]^{\frac{1}{2}}[M{+}R{-}T{+}\frac{1}{2}]^{\frac{1}{2}}[M{+}T{-}R{+}\frac{1}{2}]^{\frac{1}{2}}}{[2T]^{\frac{1}{2}}}.$

The coefficients $\Lambda^{BC}_{AD}$ satisfy additional linear relations, which will be used in the proof of Plancherel Theorem.

Lemma 6: A useful identity

$$S_{B,\aleph_0} = \sum_{C,\aleph_2}$$

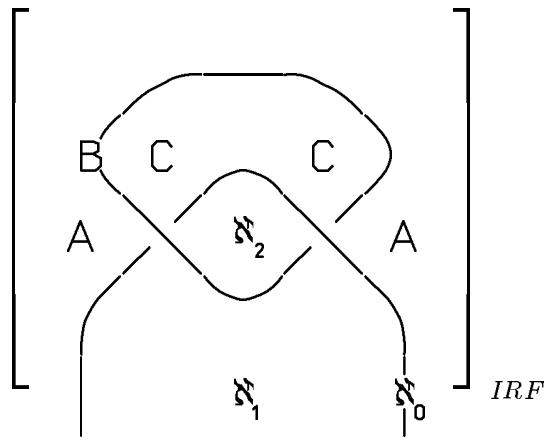

$$= (-1)^{2B}\frac{[(2B+1)(2\aleph_0+1)]}{[2\aleph_0+1]}Y(A,\aleph_0,\aleph_1),$$

*and the same formula with the undercrossing exchanged with the overcrossing.*

Proof: The proof is obtained by an easy induction on $B$. It is clear that this equation is true for $B = 0$. Furthermore, using the definition (176) it is an easy check to verify that it is also true for $B = \frac{1}{2}$. Now, we can establish a recursion relation on the half-integer $B$. Indeed, we have

$$S_{B,\aleph_0}S_{C,\aleph_0} = \sum_{DE\aleph_2\aleph_3}$$

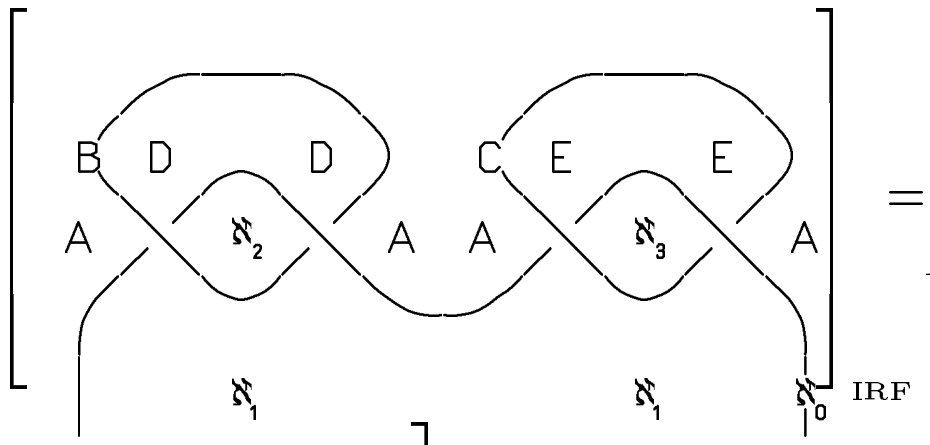

$$= \sum_{DE\aleph_2\aleph_3}$$

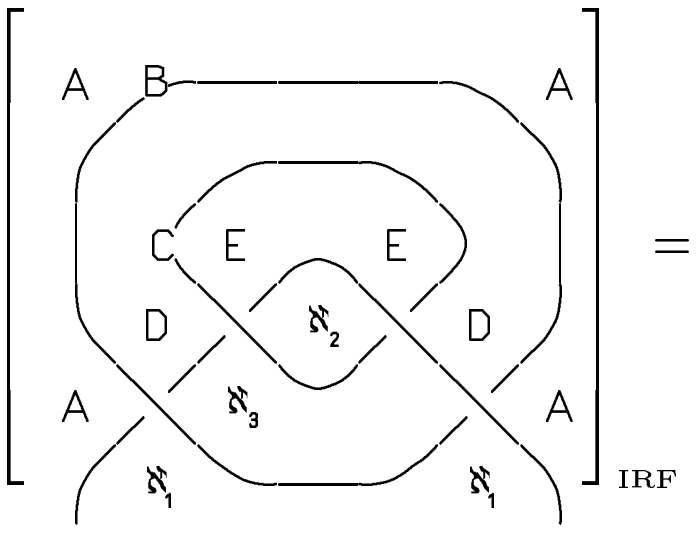

$$= \sum_{DEF\aleph_2\aleph_3}$$

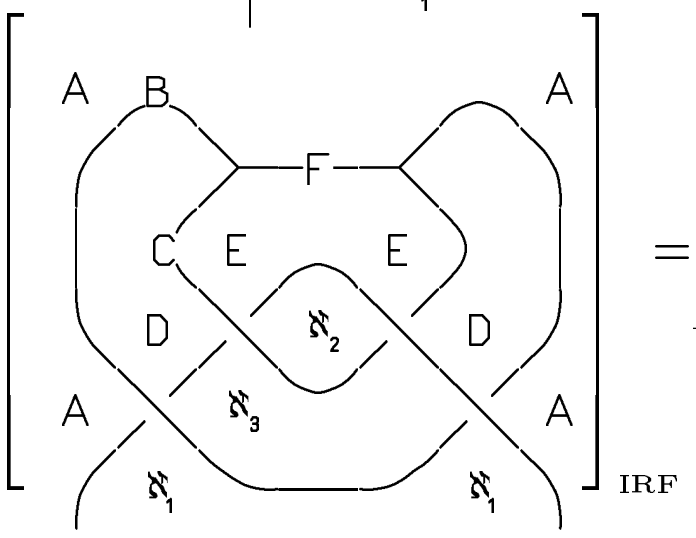

$$= \sum_{FC\aleph_2}$$

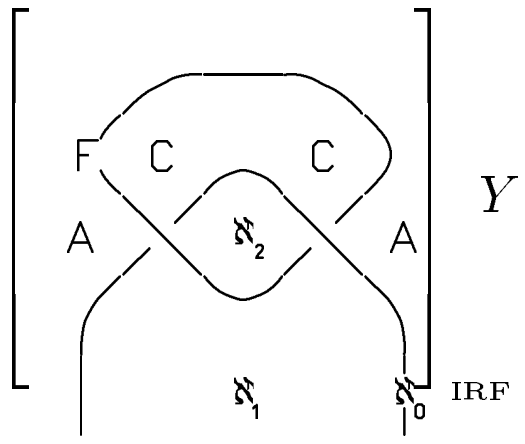

$$Y(B,C,F) = \sum_{F} Y(B,C,F)S_{F\aleph_0}$$

In particular by taking $C = \frac{1}{2}$, we have $-S_{B,\aleph_0}(q^{(2\aleph_0+1)} + q^{-(2\aleph_0+1)}) = S_{B+\frac{1}{2},\aleph_0} + S_{B-\frac{1}{2},\aleph_0}$. As a result, because $(-1)^{2B}\frac{[(2B+1)(2\aleph_0+1)]}{[2\aleph_0+1]}Y(A,\aleph_0,\aleph_1)$ satisfies the same functional equation, it is equal to $S_{B,\aleph_0}$. □

**Proposition 5**: Action of Casimir elements

*The coefficients* $\Lambda_{AD}^{BC}(\aleph_0, \aleph_1)$ *satisfy the following linear relations:*

$$\sum_{MN} \Lambda_{AD}^{MN} \begin{Bmatrix} B & I & M \\ N & A & C \end{Bmatrix} \begin{Bmatrix} B & I & M \\ N & D & C \end{Bmatrix} [d_N] (\frac{v_M}{v_N})^{\pm\frac{1}{2}} = \omega^{\pm}(I)[d_C] \Lambda_{AD}^{BC} (\frac{v_B}{v_C})^{\pm\frac{1}{2}} \tag{135}$$

*for all* $A, B, C, D, I \in \frac{1}{2}\mathbb{Z}^+$, *with the definition* $\omega^+(I) = \frac{[(2I+1)(2\aleph_0+1)]}{[2\aleph_0+1]}, \omega^-(I) = \frac{[(2I+1)(2\aleph_1+1)]}{[2\aleph_1+1]}$. *This system of linear equation, when* $I = \frac{1}{2}$ *is equivalent to the fact that* $\chi_\pi$ *is an eigenvector of the right action of* $\Omega^{(\pm)}$ *with eigenvalues* $\omega_\pm$.

Proof: From the relations satisfied by the generalized $6j$ symbols, we have the following chain of equations:

$\sum_{MN\aleph_2}$ 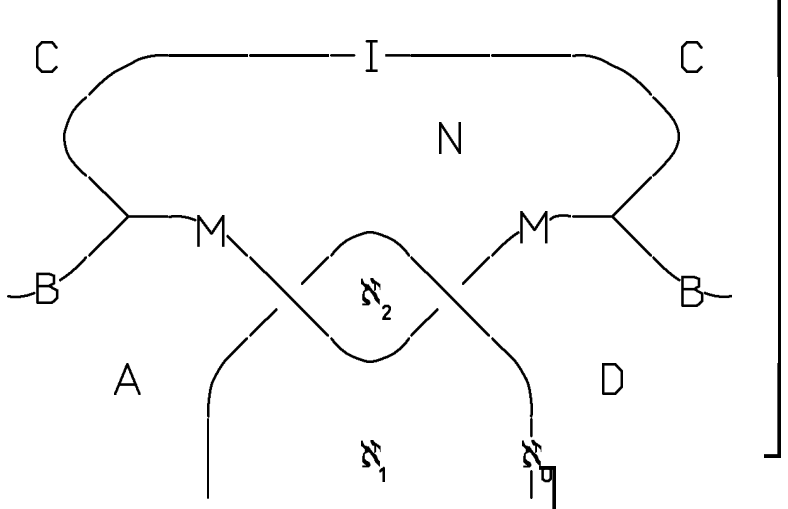 $_{\mathrm{IRF}}$ $= \sum_{MN\aleph_2}$ 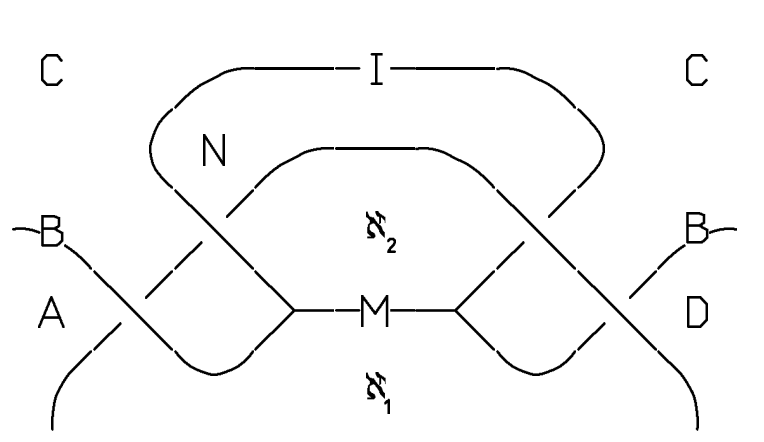 $_{\mathrm{IRF}}$ $=$

$= \sum_{N\aleph_2}$ 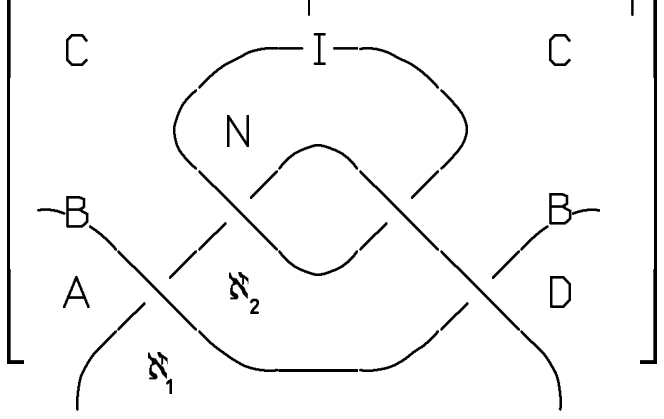 $_{\mathrm{IRF}}$ $= (-1)^{(2I)} \frac{[(2I+1)(2\aleph_0+1)]}{[(2\aleph_0+1)]}$ 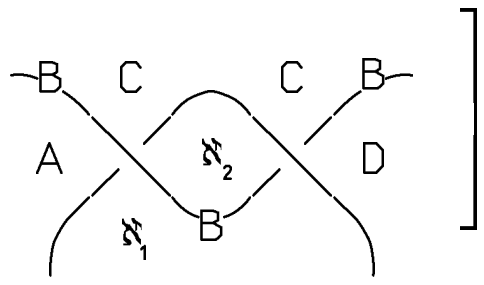 $_{\mathrm{IRF}}$

where we have used the previous lemma to obtain the last equation. If we replace the graphics by their explicit expressions in terms of $6j$ symbols, and complex continuation of $6j$ symbols, the left handside of this chain of equations is the opposite of lefthandside of (135) and similarly for the righthandside. The equation for $\omega^-(I)$ is obtained similarly, using the other equivalent definition of $\Lambda_{AD}^{BC}$ (109). It is an easy computation to show that indeed $\chi_\pi$ is an eigenvector of the right action of $\Omega^{(\pm)}$ with eigenvalues $\omega_\pm$ if and only if this linear relations are satisfied with $A = D$ and $I = \frac{1}{2}$. □

We will need the explicit form of these linear equations when $I = \frac{1}{2}$ the system is equivalent to:

$$\begin{aligned} &[B{+}K{-}A{+}1][A{+}B{+}K{+}2]q^{\pm(K-B)}\Lambda_A^{B+\frac{1}{2}K+\frac{1}{2}} - [A{+}K{-}B][A{+}B{-}K{+}1]q^{\mp(K+B+1)}\Lambda_A^{B+\frac{1}{2}K-\frac{1}{2}} + \\ &+[A{+}K{+}B{+}1][B{+}K{-}A]q^{\mp(K-B)}\Lambda_A^{B-\frac{1}{2}K-\frac{1}{2}} - [A{+}B{-}K][A{+}K{-}B{+}1]q^{\pm(K+B+1)}\Lambda_A^{B-\frac{1}{2}K+\frac{1}{2}} = \\ &= \omega_\pm[2K{+}1][2B+1]\Lambda_A^{BK} \end{aligned} \tag{136}$$

Let $A \in \frac{1}{2}\mathbb{Z}^+$ we will define $\Gamma_A = \{(B, C) \in (\frac{1}{2}\mathbb{Z}^+)^{\times 2}, Y(A, B, C) = 1\}$, and we will denote by $\partial\Gamma_A^{d+}, \partial\Gamma_A^{d-}, \partial\Gamma_A^{a}$ the "boundaries" of $\Gamma_A$ ($d$ is for diagonal, and $a$ is for antidiagonal) defined by $\partial\Gamma_A^{d\pm} = \{(B, C) \in \Gamma_A, B - C = \pm A\}$, $\partial\Gamma_A^{a} = \{(B, C) \in \Gamma_A, B + C = A\}$.

An inspection of the system of linear equations (136) shows that if we fix $A \in \frac{1}{2}\mathbb{Z}^+$, and $\Lambda_A^{0A}$ the solution $\Lambda_A^{BC}$ for $(B, C) \in \Gamma_A$ of (136), if it exists, is unique. The idea of the proof is very simple, it uses the fact that we can always eliminate one of the indeterminate $\Lambda_A^{B+\frac{\sigma}{2}K+\frac{\epsilon}{2}}$ in the system (136) by taking linear combination of the equation containing $\omega_+$ and the equation containing $\omega_-$.

As a result we can compute $\Lambda_A^{BC}$ on the antidiagonal $\partial\Gamma_A^a$ by eliminating the point $\Lambda_A^{B+\frac{1}{2}C+\frac{1}{2}}$. The point $\Lambda_A^{B-\frac{1}{2}C-\frac{1}{2}}$ does not contribute because its coefficient vanishes (we also assumed that $\Lambda_A^{BC}=0$ for $Y(A,B,C)=0$ which is then the case). As a result we obtain a linear recurrent equation which is written as:

$$[A+C-B][A+B-C+1][2C+1]\Lambda_A^{B+\frac{1}{2}\,C-\frac{1}{2}}-[A+B-C][A+C-B+1][2B+1]\Lambda_A^{B-\frac{1}{2}\,C+\frac{1}{2}}$$
$$=([2C+1][2B+1](\omega_+q^{-C+B}-\omega_-q^{C-B})\Lambda_A^{BC}$$

which then defines uniquely the coefficients $\Lambda_A^{BC}=0$ on $\partial\Gamma_A^a$ in term of $\omega_+,\omega_-$ and $\Lambda_A^{0A}$ because the coefficient of $\Lambda_A^{-\frac{1}{2}\,A-\frac{1}{2}}$ vanishes.

If we eliminate $\Lambda_A^{B-\frac{1}{2}\,C-\frac{1}{2}}$, we obtain a linear system which is triangular with respect to the gradation defined by $B+C$. Because $\Lambda_A^{BC}$ are known on the boundary $\Gamma_A^a$, this triangular system defines uniquely the $\Lambda_A^{BC}$. This triangular system is an efficient tool to compute $\Lambda_A^{BC}$.

# 5 Plancherel Theorem

In this section we give a proof of Plancherel formula for $SL_q(2,\mathbb{C})_{\mathbb{R}}$. Let us denote by $\pi(m,\rho)$ the element of $\hat{\mathcal{D}}_p$ associated to $m\in\mathbb{N},\rho\in[-\frac{2\pi}{\hbar},\frac{2\pi}{\hbar}]=I_\hbar$ , in this case we have $2\aleph_0+1=\frac{1}{2}(m+i\rho)$ and $2\aleph_1+1=\frac{1}{2}(-m+i\rho)$. Note that we have doubled the interval for $m=0$, and that in this case we have $\pi(0,\rho)=\pi(0,-\rho)$ for $\rho\in I_\hbar$. Here we will often denote abusively $\Lambda_A^{BC}(m,\rho)$ the coefficients defining this principal serie of representation, i.e. $\Lambda_A^{BC}(m,\rho)=\Lambda_A^{BC}(\aleph_0,\aleph_1)$ and $\chi(m,\rho)=\chi_{\pi(m,\rho)}$.
Let us define the Plancherel density

$$\begin{aligned}P(m,\rho)&=\frac{\hbar}{2\pi}(1-\frac{1}{2}\delta_{m,0})(\cosh(m\hbar)-\cos(\hbar\rho))\\&=-\frac{\hbar}{4\pi}(q-q^{-1})^2(1-\frac{1}{2}\delta_{m,0})[2\aleph_0+1][2\aleph_1+1],\forall m\in\mathbb{N},\rho\in I_\hbar.\end{aligned}$$

We can equivalently see this Plancherel density as a measure on the space $\hat{\mathcal{D}}_p$, denoted $dP(\pi)$. We will prove the following theorem:

Theorem 7: 1.Plancherel Formula for $SL_q(2,\mathbb{C})_{\mathbb{R}}$

*Let $\phi$ be any element of $Fun_{cc}(SL_q(2,\mathbb{C})_{\mathbb{R}})$, we have :*

$$\sum_{m=0}^{+\infty}\int_{I_\hbar}\chi(m,\rho)(\phi)P(m,\rho)d\rho=\epsilon(\phi). \tag{137}$$

The proof of this theorem consists in three lemmas which are interesting in themselves.

Lemma 8:

*Assume that* $(f(m,.))_{m\in\mathbb{N}}$ *is a family of* $C^1$ *even functions on* $I_\hbar$ *such that Plancherel formula holds, i.e*

$$\sum_{m=0}^{+\infty}\int_{I_\hbar}\chi(m,\rho)(\phi)f(m,\rho)d\rho=\epsilon(\phi)\quad\forall\phi\in Fun_{cc}(SL_q(2,\mathbb{C})_{\mathbb{R}}) \tag{138}$$

*then* $f(m,\rho)=P(m,\rho)$.

Proof: We first compute $h(\ \chi(m,\rho)\ S^{-1}(\overset{A}{k}{}^t_n)\otimes\overset{B}{E}{}^r_s)$, this is equal to:

$$h(\sum_{CDF}\overset{C}{k}{}^u_j S^{-1}(\overset{A}{k}{}^t_n)\otimes\overset{D}{E}{}^p_l\overset{B}{E}{}^r_s\begin{pmatrix} j & l & F\\ C & D & x\end{pmatrix}\begin{pmatrix} x & D & C\\ F & p & i\end{pmatrix}\overset{C}{\mu}{}^{-1i}{}_u\Lambda^{DF}_C(m,\rho))$$

$$=\sum_{CDF}\frac{\delta_{AC}\overset{A}{\mu}{}^u_n}{[d_A]}\delta^t_j\delta_{BD}\overset{B}{\mu}{}^{-1p}{}_s[d_B]\delta^r_l\begin{pmatrix} j & l & F\\ C & D & x\end{pmatrix}\begin{pmatrix} x & D & C\\ F & p & i\end{pmatrix}\overset{C}{\mu}{}^{-1i}{}_u\Lambda^{DF}_C(m,\rho)$$

$$=\frac{[d_B]}{[d_A]}\sum_F\begin{pmatrix} t & r & F\\ A & B & x\end{pmatrix}\begin{pmatrix} x & B & A\\ F & p & n\end{pmatrix}\overset{B}{\mu}{}^{-1p}{}_s\Lambda^{BF}_A(m,\rho).$$

Assume that Plancherel formula holds with the densities $f(m,.)$, the previous computation implies that Plancherel formula is equivalent to:

$$\frac{[d_B]}{[d_A]}\sum_{m=0}^{+\infty}\int_{I_\hbar}\sum_F\begin{pmatrix} t & r & F\\ A & B & x\end{pmatrix}\overset{F}{\mu}{}^{-1x}{}_y\begin{pmatrix} y & B & A\\ F & s & n\end{pmatrix}\Lambda^{BF}_A(m,\rho)f(m,\rho)d\rho=\overset{A}{\mu}{}^{-1t}{}_n\delta_{B,0}$$

Therefore if we multiply both sides of this equation by $P_{12}\overset{AB}{R}{}^{\pm}$ and if we take the trace on $\overset{A}{V}\otimes\overset{B}{V}$ we obtain that:

$$\sum_F\frac{[d_B][d_F]}{[d_A]}\sum_{m=0}^{+\infty}\int_{I_\hbar}\Lambda^{BF}_A(m,\rho)f(m,\rho)(\frac{v_Av_B}{v_F})^{\pm\frac{1}{2}}d\rho=[d_A]\delta_{B,0}. \tag{139}$$

Using the relation

$$\sum_F\frac{[d_F]}{[d_A]}\Lambda^{BF}_A(m,\rho)(\frac{v_Av_B}{v_F})^{\pm\frac{1}{2}}=\tilde{Y}(A,\frac{m}{2})\frac{[(2B+1)\frac{(\pm m+i\rho)}{2}]}{[\frac{(\pm m+i\rho)}{2}]}$$

proved in the previous lemma, the system of equations (139) is equivalent to:

$$\sum_{m|A\pm\frac{m}{2}\in\mathbb{N}}\int_{I_\hbar}\frac{[(2B+1)\frac{(\pm m+i\rho)}{2}]}{[\frac{(\pm m+i\rho)}{2}]}\ f(m,\rho)\ =[d_A]\delta_{B,0}. \tag{140}$$

We can develop $f(m,.)$ in Fourier series and write: $f(m,\rho)=\sum_{p\in\mathbb{Z}}a^m_p q^{\frac{ip\rho}{2}}$. Using the relation:

$$\frac{[(2B+1)\frac{(\pm m+i\rho)}{2}]}{[\frac{(\pm m+i\rho)}{2}]}=\sum_{k=0}^{2B}q^{(k-B)(\pm m+i\rho)},$$

the system (140) is equivalent to a linear system on $a_m^p$ which is trivial to solve and which unique non zero solution under the assumption that $f(0,\rho)$ is an even function gives $f(m,\rho) = P(m,\rho)$. In particular, only the Fourier modes $(0,-1,1)$ of $f(m,.)$ are non zero. □

It remains to show that the system of equations

$$\frac{[d_B]}{[d_A]} \sum_{m=0}^{+\infty} \int_{I_\hbar} \sum_F \left( \begin{array}{cc|c} t & r & F \\ A & B & x \end{array} \right) {}^{F}\mu^{-1}{}^{x}_{y} \left( \begin{array}{c|cc} y & B & A \\ F & s & n \end{array} \right) \Lambda_A^{BF}(m,\rho) P(m,\rho) d\rho = {}^{A}\mu^{-1}{}^{t}_{n} \delta_{B,0} \qquad (141)$$

is satisfied. This is trivially the case for $B = 0$, because if we put $B = 0$ in the system (139) , we get:

$$\sum_{m=0}^{+\infty} \int_{I_\hbar} \Lambda_A^{0A}(m,\rho) f(m,\rho) d\rho = [d_A]$$

which is trivially equivalent to (141) with $B = 0$.

Therefore, the proof of Plancherel formula is reduced to the proof of (141) when $B \neq 0$, which is equivalent to:

$$\sum_{m=0}^{+\infty} \int_{I_\hbar} \Lambda_C^{BA}(m,\rho) f(m,\rho) d\rho = 0, \quad B \neq 0. \qquad (142)$$

Lemma 9:

*In order to show that equations (142) are true, it is sufficient to show that they hold for any $(B,C) \in \partial\Gamma_A^{d\pm}, \partial\Gamma_A^a$ where $A$ is any half-integer and $B \neq 0$.*

Proof: The idea is to use the linear system coming from the invariance under $\mathfrak{U}_q(an(2))$. Let us define $\nu_A^{BC} = \sum_{m=0}^{+\infty} \int_{I_\hbar} \Lambda_A^{BC}(m,\rho) f(m,\rho) d\rho$, and assume that $\nu_A^{BC} = 0, \forall A \in \frac{1}{2}\mathbb{Z}^+, (B,C) \in \partial\Gamma_A^{d\pm}, \partial\Gamma_A^a, B \neq 0$. By linearity, the system (134) is trivially satisfied by $\nu_A^{BC}$. It is triangular with respect to the gradation $A+B+C$, therefore if we assume that that $\nu_A^{BC} = 0$ for all $(B,C) \in \partial\Gamma_A^{d\pm}, \partial\Gamma_A^a$, $B \neq 0$ we deduce, by an easy induction, that $\nu_{A+1}^{BC} = 0, \forall B, C \in \frac{1}{2}\mathbb{Z}^+$, $B \neq 0$ except for the special case $B = 0, C = A+1$. Indeed in this case the linear system express $\nu_{A+1}^{1A+1}$ in terms of linear combination of $\nu_A^{BC}, B+C < A+3, \nu_{A+1}^{BC}, B+C < A+2$ containing $\nu_{A+1}^{0A+1}$ and $\nu_A^{0A}$ which are non zero. The linear system (134) applied to this case gives, using the induction assumption,

$$\frac{[2A+4]}{[2A+3]} \nu_{A+1}^{1A+1} = \frac{[2A+2]}{[2A+1]} \nu_A^{0A} - \frac{[2A+2]}{[2A+3]} \nu_{A+1}^{0A+1} = \frac{[2A+2]}{[2A+1]} [2A+1] - \frac{[2A+2]}{[2A+3]} [2A+3] = 0. \quad (143)$$

Which concludes the proof of the lemma □

Lemma 10:

*The explicit expressions for $\Lambda_A^{BC}$ on the boundaries of $\Gamma_A$ are:*

$$\begin{aligned}
\Lambda_A^{BC}(m,\rho) &= \tilde{Y}(A,\frac{m}{2})\frac{(-1)^{2B}}{\binom{2A}{2B}_q}\sum_{k=-B}^{B} q^{-ik\rho}\binom{A+\frac{m}{2}}{B+k}_q\binom{A-\frac{m}{2}}{B-k}_q, \quad (\partial\Gamma^a)\\
\Lambda_A^{BC}(m,\rho) &= \frac{\tilde{Y}(A,\frac{m}{2})}{\binom{2C+1}{2B}_q}\sum_{k=-B}^{B} q^{-ik\rho}\binom{A+\frac{m}{2}+B-k}{B-k}_q\binom{A-\frac{m}{2}+B+k}{B+k}_q, \quad (\partial\Gamma^{d-})\\
\Lambda_A^{BC}(m,\rho) &= \tilde{Y}(A,\frac{m}{2})\frac{(-1)^{2A}q^{-i\frac{\rho m}{2}}}{\binom{2B+1}{2C}_q}\sum_{k=-C}^{C} q^{-ik\rho}\binom{A-\frac{m}{2}+C-k}{C-k}_q\binom{A+\frac{m}{2}+C+k}{C+k}_q, \quad (\partial\Gamma^{d+}).
\end{aligned} \tag{144}$$

*The equations (141) are satisfied on $\partial\Gamma_A^{d\pm},\partial\Gamma_A^a$.*

Proof: In order to show that $\Lambda_A^{BC}$ are given by the expressions (144) we have to integrate the linear equations (136). This is trivial, because we know that these systems admit a unique solution and it is then sufficient to show that the right handside of (144) satisfy these linear equations, which is easy to show. We will first show that the equations (141) are satisfied on $\partial\Gamma_A^a$.

From the expression of $\Lambda_A^{BC}$ on $\partial\Gamma_A^a$, we see that if $2B$ is an odd integer then $\Lambda_A^{BC}$ has just Fourier modes in $\mathbb{Z}+\frac{1}{2}$. As a result the integration $\int_{I_\hbar}\Lambda_A^{BC}(m,\rho)P(m,\rho)d\rho = 0$, and Plancherel formula is trivially satisfied. The only non trivial check is therefore for $2B$ even. From the explicit expression of $\Lambda_A^{BC}(m,\rho)$ we obtain that

$$\begin{aligned}
&(-1)^{2B}\binom{2A}{2B}_q\int_{I_\hbar}\Lambda_A^{BC}(m,\rho)P(m,\rho)d\rho =\\
&= \tilde{Y}(A,\frac{m}{2})(\binom{A+\frac{m}{2}}{B}_q\binom{A-\frac{m}{2}}{B}_q(q^m+q^{-m}) - \binom{A+\frac{m}{2}}{B+1}_q\binom{A-\frac{m}{2}}{B-1}_q - \binom{A+\frac{m}{2}}{B-1}_q\binom{A-\frac{m}{2}}{B+1}_q).
\end{aligned}$$

We have to consider two cases, the first one is when $A$ is an integer, the constraint $\tilde{Y}(A,\frac{m}{2})=1$ imposes that we put $m=2s$, the Plancherel formula (142), is then equivalent to the following identity:

$$\begin{aligned}
&\sum_{s=0}^{A}\{\binom{A+s}{B}_q\binom{A-s}{B}_q(q^{2s}+q^{-2s}) - \binom{A+s}{B+1}_q\binom{A-s}{B-1}_q - \binom{A+s}{B-1}_q\binom{A-s}{B+1}_q)\} =\\
&= \binom{A}{B}_q\binom{A}{B}_q - \binom{A}{B+1}_q\binom{A}{B-1}_q.
\end{aligned} \tag{145}$$

The second case is when $A=l+\frac{1}{2}$, with $l\in\mathbb{Z}$, the constraint $\tilde{Y}(A,\frac{m}{2})=1$ imposes that we put $m=2s+1$ with $0\leq s\leq l$, the Plancherel formula (142), is then equivalent to the following identity:

$$\sum_{s=0}^{l}\{\binom{l+s+1}{B}_q\binom{l-s}{B}_q(q^{2s+1}+q^{-2s-1}) - \binom{l+s+1}{B+1}_q\binom{l-s}{B-1}_q - \binom{l+s+1}{B-1}_q\binom{A-s}{B+1}_q)\} = 0$$

We will now show that these identities on $q$-binomials hold, and will only give an explicit proof of the first identity, the second case uses the same method of proof. We will show the two identities, for $B \neq 0$

$$\sum_{s=0}^{A}\{\binom{A+s}{B}_q\binom{A-s}{B}_q q^{-2s}-\binom{A+s}{B+1}_q\binom{A-s}{B-1}_q\}=-q^{A-B+1}\binom{A}{B+1}_q\binom{A+1}{B}_q \quad (146)$$

$$\sum_{s=0}^{A}\{\binom{A+s}{B}_q\binom{A-s}{B}_q q^{2s}-\binom{A-s}{B+1}_q\binom{A+s}{B-1}_q\}=q^{A-B}\binom{A}{B}_q\binom{A+1}{B+1}_q. \quad (147)$$

If we take the sum of these two equations, the left handside is the lefthandside of (145) and the right handside is, after the use of two $q$-analogues of Pascal relations, equal to the righthandside of (145). Let us first recall that if $x, y$ are two elements satisfying the relation $yx = q^2xy$, then we have the identity:

$$(x+y)^n=\sum_{p=0}^{n}\binom{n}{p}_q q^{p(n-p)}x^p y^{n-p} \quad (148)$$

Let consider the $\mathbb{C}$ algebra generated by $x, y, a, b$ and satisfying the relations $yx = q^2xy, ba = q^2ab, ax = xa, ay = ya, bx = xb, by = yb$. A basis of this algebra is $(a^m y^n a^r b^s)_{m,n,r,s\in\mathbb{N}}$. If $P = \sum_{m,n,r,s\in\mathbb{N}} P_{mnrs}a^m y^n a^r b^s$ we will denote $C_{m,r}(P) = \sum_{n,s\in\mathbb{N}} P_{mnrs}$. From the defining relations, we have:

$$\binom{A+s}{B}_q\binom{A-s}{B}_q=C_{B,B}(q^{-2(A-B)B}(x+y)^{A+s}(a+b)^{A-s}),$$

from which we deduce that:

$$\sum_{s=0}^{A}\binom{A+s}{B}_q\binom{A-s}{B}_q q^{-2s}=q^{-2(A-B)B}C_{B,B}((x+y)^A(a+b)^A\sum_{s=0}^{A}(\frac{x+y}{a+b}q^{-2})^s)$$

$$=q^{-2(A-B)B}C_{B,B}((x+y)^A\frac{(a+b)^{A+1}-(q^{-2}(x+y))^{A+1}}{a+b-q^{-2}(x+y)})$$

where the reader is urged to check that because $a + b$ commutes with $x + y$, these expressions have a perfect definition. A similar reasonning implies that

$$\binom{A+s}{B+1}_q\binom{A-s}{B-1}_q=q^{-2(A-B)B+2}C_{B,B}(x^{-1}a(x+y)^A\frac{(a+b)^{A+1}-(q^{-2}(x+y))^{A+1}}{a+b-q^{-2}(x+y)}).$$

(Note that we have localized the algebra in $x$, with a trivial extension of the definition of $C_{m,r}$.) Therefore we have

$$\sum_{s=0}^{A}\{\binom{A+s}{B}_q\binom{A-s}{B}_q q^{-2s}-\binom{A+s}{B+1}_q\binom{A-s}{B-1}_q\}=$$

$$=-q^{-2(A-B)B}C_{B,B}(x^{-1}q^2(a-q^{-2}x)(x+y)^A\frac{(a+b)^{A+1}-(q^{-2}(x+y))^{A+1}}{a+b-q^{-2}(x+y)}).$$

We now use the following trick

$$C_{B,B}(x^{-1}(b-q^{-2}y)a^m y^n a^r b^s)=0$$

which is trivial to show, and obtain

$$\sum_{s=0}^{A}\{\begin{pmatrix}A+s\\ B\end{pmatrix}_q\begin{pmatrix}A-s\\ B\end{pmatrix}_q q^{-2s}-\begin{pmatrix}A+s\\ B+1\end{pmatrix}_q\begin{pmatrix}A-s\\ B-1\end{pmatrix}_q\}=$$
$$=-q^{-2(A-B)B}C_{B,B}(x^{-1}q^2(a-q^{-2}x+b-q^{-2}y)(x+y)^A\frac{(a+b)^{A+1}-(q^{-2}(x+y))^{A+1}}{a+b-q^{-2}(x+y)})$$
$$=-q^{-2(A-B)B}C_{B,B}(x^{-1}q^2(x+y)^A((a+b)^{A+1}-(q^{-2}(x+y))^{A+1}))$$
$$=-q^{A-B+1}\begin{pmatrix}A\\ B+1\end{pmatrix}_q\begin{pmatrix}A+1\\ B\end{pmatrix}_q$$

which ends the proof of the first relation (146). The second relation (147) is proved with exactly the same method. The proofs of (141) on $\partial\Gamma^{d\pm}_A$ are similar and we leave them to the reader.

□

Let $\phi\in Fun_{cc}(SL_q(2,\mathbb{C})_{\mathbb{R}})$, it is trivial to show that $\epsilon(\phi * S(\phi^\star))=h(\phi^\star\phi)$. As a result we obtain using $\pi(\phi * S(\phi^\star))=\pi(\phi)\pi(\phi)^\dagger$ and the Plancherel formula that:

Theorem 11: Plancherel Formula .2

*Let $\phi$ be any element of $Fun_{cc}(SL_q(2,\mathbb{C})_{\mathbb{R}})$, we have :*

$$||\phi||^2_{L^2}=\int_{\hat{\mathcal{D}}_p} tr(\pi(\mu^{-1})\pi(\phi)\pi(\phi)^\dagger)dP(\pi). \tag{149}$$

$V_\pi$ admits a Hilbert basis $\mathcal{E}=\{\overset{A}{e}_m, A\in\frac{1}{2}\mathbb{Z}^+, m=-A,...,A\}$, and denote by $\mathcal{L}_{qHS}(V_\pi)$ the Hilbert space $\mathcal{L}_{qHS}(V_\pi)=\{\xi\in\mathbb{C}^{\mathcal{E}\times\mathcal{E}},\sum_{x,y\in\mathcal{E}}\mu^{-1}_{xx}|\xi(x,y)|^2<+\infty\}$, with Hermitian product: $<\xi,\eta>=\sum_{x,y\in\mathcal{E}}\mu^{-1}_{xx}\eta(x,y)\overline{\xi}(x,y),\forall\xi,\eta\in\mathcal{L}_{qHS}(V_\pi)$.

When $q=1$, $\mu=1$ and $\mathcal{L}_{qHS}(V_\pi)$ is the Hilbert space of Hilbert-Schmidt bounded operator of $V_\pi$. If $\phi\in Fun_{cc}(SL_q(2,\mathbb{C})_{\mathbb{R}})$, we can associate to $\pi(\phi)$ an element $\xi_\pi(\phi)\in\mathcal{L}_{qHS}$ defined by $\xi_\pi(\phi)(x,y)=<x|\pi(\phi)(y)>,\ \ x,y\in\mathcal{E}$. We have $<\xi_\pi(\phi),\xi_\pi(\phi)>=tr(\pi(\mu^{-1})\pi(\phi)\pi(\phi)^\dagger)$.

From this last proposition we have defined an isometric mapping $T$ from $Fun_{cc}(SL_q(2,\mathbb{C})_{\mathbb{R}})$ with the norm $||.||_{L^2}$ to the Hilbert space $\mathcal{H}$ defined as the direct integral $\mathcal{H}=\int^{\oplus}_{\hat{\mathcal{D}}_p}\mathcal{L}_{qHS}(V_\pi)dP(\pi)$, with the following definition: $T(\phi)_\pi=\xi_\pi(\phi)$.

As a result, $T$ being uniformly continuous we can uniquely extend $T$ to an isometric map $T^\#$ from $L^2(SL_q(2,\mathbb{C})_{\mathbb{R}})$ to $\int^{\oplus}_{\hat{\mathcal{D}}_p}\mathcal{L}_{qHS}(V_\pi)dP(\pi)$.

Theorem 12: Plancherel Theorem

*The spherical transform $T^\#:L^2(SL_q(2,\mathbb{C})_{\mathbb{R}})\to\int^{\oplus}_{\hat{\mathcal{D}}_p}\mathcal{L}_{qHS}(V_\pi)dP(\pi)$ is a unitary operator.*

Proof: Because $T^\#$ is an isometry, it is an injection. The only non trivial point is to show that $T^\#$ is a surjection. Because $T^\#$ is an extension of $T$, it is sufficient to show that $T(Fun_{cc}(SL_q(2,\mathbb{C})_{\mathbb{R}}))$ is dense in $\mathcal{H}$, which is equivalent to show that $T(Fun_{cc}(SL_q(2,\mathbb{C})_{\mathbb{R}}))^\perp=\{0\}$. Let $\eta\in\mathcal{H}$ and assume that $\int_{\hat{\mathcal{D}}_p}<\eta_\pi,\xi_\pi(\phi)>dP(\pi)=0$ for any $\phi\in Fun_{cc}(SL_q(2,\mathbb{C})_{\mathbb{R}})$. From the expression of $\tilde{\pi}(\overset{A}{X}{}^i_j\otimes\overset{B}{g}{}^k_l)$ in the orthonormal basis $\mathcal{E}$, this orthogonality condition is equivalent to the condition

$$\forall A,B,C,D\in\frac{1}{2}\mathbb{Z}^+,\int_{\hat{\mathcal{D}}_p}\overline{\eta_\pi}(\overset{A}{e}_i,\overset{D}{e}_j)\Lambda^{BC}_{AD}(\pi)dP(\pi)=0. \tag{150}$$

As a result it will be sufficient to show that if $(f_m)_{m\leq inf(A,D)} \in L^2(I_\hbar)$, with $(\tilde{Y})(A, \frac{m}{2}) = (\tilde{Y})(D, \frac{m}{2}) = 1$, is such that

$$\sum_m \int_{I_\hbar} f_m(\rho)\Lambda^{BC}_{AD}(m,\rho)dP(m,\rho) = 0, \forall B, C \in \frac{1}{2}\mathbb{Z}^+ \tag{151}$$

then $f_m = 0, \forall m$. Let us prove this result. If we assume that $\sum_m \int_{I_\hbar} f_m(\rho)\Lambda^{BC}_{AD}(m,\rho)dP(m,\rho) = 0, \forall B, C$, we obtain using the linear equations (135) that

$$\sum_m \int_{I_\hbar} \omega^+(I)\omega^-(J) f_m(\rho)\Lambda^{BC}_{AD}(m,\rho)dP(m,\rho) = 0, \forall B, C, I, J.$$

This can be written as

$$\sum_m \int_{I_\hbar} [\frac{r}{2}(-m+i\rho)][\frac{s}{2}(+m+i\rho)] f_m(\rho)\Lambda^{BC}_{AD}(m,\rho)d\rho = 0, \forall r, s \in \mathbb{Z}^+. \tag{152}$$

We will first show that this implies that $\int_{I_\hbar} f_m(\rho)\Lambda^{BC}_{AD}(m,\rho)q^{i\frac{n}{2}\rho}d\rho = 0, \forall n \in \mathbb{Z}$. This is done by an induction on $m_0$ the largest integer appearing in the finite sum (151). Indeed from the equation (152), we obtain the equality for every integer $p, r, s$

$$\int_{I_\hbar} f_{m_0}(\rho)\Lambda^{BC}_{AD}(m_0,\rho)q^{i\frac{n}{2}\rho}d\rho = \tag{153}$$

$$= \int_{I_\hbar} f_{m_0}(\rho)\Lambda^{BC}_{AD}(m_0,\rho)(q^{\frac{pm_0}{2}}[\frac{r}{2}(-m_0+i\rho)][\frac{s}{2}(+m_0+i\rho)] + q^{i\frac{n}{2}\rho})d\rho + \tag{154}$$

$$\sum_{m<m_0} \int_{I_\hbar} f_m(\rho)\Lambda^{BC}_{AD}(m,\rho)q^{\frac{pm_0}{2}}[\frac{r}{2}(-m+i\rho)][\frac{s}{2}(+m+i\rho)]d\rho. \tag{155}$$

As a result we get:

$$\begin{aligned} &|\int_{I_\hbar} f_{m_0}(\rho)\Lambda^{BC}_{AD}(m_0,\rho)q^{i\frac{n}{2}\rho}d\rho| \\ &\leq ||\Lambda^{BC}_{AD}(m_0,.)f_{m_0}||_{L^2}||g_{r,s}||_{L^2} + \sum_{m<m_0} ||\Lambda^{BC}_{AD}(m,.)f_m||_{L^2}||h_{r,s,m}||_{L^2} \end{aligned} \tag{156}$$

where $g_{r,s}, h_{r,s,m}$ are the functions: $g_{r,s}(\rho) = q^{\frac{pm_0}{2}}[\frac{r}{2}(-m_0+i\rho)][\frac{s}{2}(+m_0+i\rho)] + q^{i\frac{n}{2}\rho}$ and $h_{r,s,m}(\rho) = q^{\frac{pm_0}{2}}[\frac{r}{2}(-m+i\rho)][\frac{s}{2}(+m+i\rho)]$. Let us define a sequence $r_p, s_p$ by $r_p - s_p = n, r_p + s_p = p$, (with $p$ odd or even depending on the parity of $n$), then it is easily shown that $g_{r_p,s_p}, h_{r_p,s_p,m}$ are uniformally convergent toward 0 when $p$ goes to infinity. As a result, by fixing $r = r_p, s = s_p$ in (156) and taking the limit $p \to +\infty$, we obtain that $\int_{I_\hbar} f_{m_0}(\rho)\Lambda^{BC}_{AD}(m_0,\rho)q^{i\frac{n}{2}\rho}d\rho = 0$. This implies that $f_{m_0}\Lambda^{BC}_{AD}(m_0,.) = 0$ as a $L^2$ function. There always exist $B, C$ such that $\Lambda^{BC}_{AD}(m_0,\rho)$ is a non constant function, and because $\Lambda^{BC}_{AD}(m_0,\rho)^2$ is a trigonometric polynomial, it has only a finite number of zeros. As a result we obtain that $f_{m_0} = 0$ as a $L^2$ function. A trivial induction on $m_0$ shows that $f_m = 0$ for all $m \neq 0$. The induction ends up when $m = 0$, in this case (152) implies that $\int_{I_\hbar}[i\frac{r}{2}\rho][i\frac{s}{2}\rho]f_0(\rho)\Lambda^{BC}_{AD}(0,\rho)d\rho = 0$, from which we deduce that $f_0(\rho) = 0$ because we assumed $f_0(\rho) = f_0(-\rho)$. This concludes the proof. □

# 6 Conclusion

The main theme of this article was the recognition of the central role played by generalized $6j$ symbols in the harmonic analysis of $SL_q(2,\mathbb{C})_{\mathbb{R}}$. There are many ways in which our work can be pursued.

From a purely mathematical point of view, it would be nice to generalize our work to the quantization of other complex Lie algebras. This should be quite a hard and technical problem, involving the precise study of generalized $6j$ symbols to compact quantum groups, but the main structures that we develop in our article should not suffer a lot of modifications, apart from the multiplicities, because many of the proofs were graphical.

It would be interesting to study the precise image of different classes of functions on $SL_q(2,\mathbb{C})_{\mathbb{R}}$ by $T^{\#}$. In particular there should exists a whole class of theorems, q-analogues of Paley-Wiener theorems. We also do not precisely know how to recover the standard characters of the unitary serie of $SL(2,\mathbb{C})_{\mathbb{R}}$ from the classical limit of the characters of $SL_q(2,\mathbb{C})_{\mathbb{R}}$.

Our work certainly has implication for the physics of low dimensional field theory where the non compact groups enter in an essential way. This is the case for Chern-Simons theory with the group $SL(2,\mathbb{C})_{\mathbb{R}}$ which is related to $2+1$ and $3+1$ quantum gravity with positive cosmological constant in the canonical quantization program. The program of combinatorial quantization can certainly be extended to this case, and representations of the moduli algebra [2] obtained by acting on wave packets of unitary representations. This will be a difficult task because it will necessarely amount to study the tensor product of unitary representations of $SL_q(2,\mathbb{C})_{\mathbb{R}}$ and its decomposition in terms of integral of unitary representations, which will be studied in a second paper [10].

Because complex continuation of $6j$ symbols played an important role in the study of special points in the strong coupling regime of Liouville theory, it is natural to suspect that harmonic analysis of $SL_q(2,\mathbb{C})_{\mathbb{R}}$ is hidden in this problem and that the study of the tensor product of unitary representations may be an important tool in the study of the algebra of fields of Liouville theory [15]

**Acknowlegments** We warmfully thank O.Babelon, J.L.Gervais, E. Ragoucy, J.Schnittger, J.Teschner for discussions on this work and on related topics.

# 7 Appendix

## 7.1 Conventions and commutation relations

Let $x$ be a complex number, we will denote by $[x]_q = [x]$ the q-number associated to $x$ defined by $[x] = \frac{q^x - q^{-x}}{q - q^{-1}}$. The $q-$ factorial is defined by $[n]! = \prod_{k=1}^{n} [k] \quad n \in \mathbb{Z}^{+*}$ with the convention that $[0]! = 1$. The $q-$ binomial coefficients are defined by:

$$\binom{n}{p}_q = \frac{[n]!}{[p]![n-p]!}, \quad n \geq 0, 0 \leq p \leq n. \tag{157}$$

They satisfy the two relations:

$$\binom{n+1}{p}_q = q^p \binom{n}{p}_q + q^{-n+p-1} \binom{n}{p-1}_q \tag{158}$$

$$\binom{n+1}{p}_q = q^{n+1-p} \binom{n}{p}_q + q^{-p} \binom{n}{p-1}_q \tag{159}$$

$\mathfrak{U}_q(su(2))$ is the star Hopf algebra defined by the relations (41,42,43). This Hopf algebra is quasitriangular and the expression of the $R$ matrix is given by:

$$R = q^{2J_z \otimes J_z} e_{q^{-1}}^{(q-q^{-1})\ (q^{J_z} J_+ \otimes J_- q^{-J_z})} \quad \text{with} \quad e_\alpha^z = \sum_{k=0}^{+\infty} \alpha^{-\frac{k(k-1)}{2}} \frac{z^k}{[k]_\alpha!} \tag{160}$$

The action of the generators of $A$ on an orthonormal basis for an irreducible representation of spin $J$ of $Irr(A)$ is given by the following expressions:

$$q^{J_z} \overset{I}{e}_m = q^m \overset{I}{e}_m \tag{161}$$

$$J_\pm \overset{I}{e}_m = q^{\mp\frac{1}{2}} \sqrt{[I \pm m + 1][I \mp m]} \overset{I}{e}_{m\pm1} \tag{162}$$

(the unusual presence of $q^{\mp\frac{1}{2}}$ in the left handside comes from our definition of the $\star$ on $J_\pm$).

The element $\mu$ is given by $\mu = q^{2J_z}$. It is easy to compute $v_K^2 = q^{-4K(K+1)}$, and we will choose

$$v_K^{\frac{1}{2}} = \exp(i\pi K) q^{-K(K+1)}.$$

With this choice of signs the relations of the Clebsch Gordan coefficients written in 5.2 are satisfied. The reader is advised that there exists a large number of conventions largely unexplained in the litterature which consequences are that there are lots of misprints and phase ambiguities in certain papers on this subject.

$Pol(SU_q(2))$ is the star Hopf algebra generated by the elements of

$$\overset{\frac{1}{2}}{u} = \begin{bmatrix} a & b \\ c & d \end{bmatrix}$$

satisfying the relations:

$$\begin{gathered} qab = ba \quad qac = ca \quad qbd = db \quad qcd = dc \quad bc = cb \quad ad - da = (q^{-1} - q)bc \quad ad - q^{-1}bc = 1 \\ \Delta(a) = a \otimes a + b \otimes c \quad \Delta(b) = b \otimes d + a \otimes b \quad \Delta(c) = c \otimes a + d \otimes c \quad \Delta(d) = d \otimes d + c \otimes b \\ a^\star = d \qquad d^\star = a \qquad b^\star = -q^{-1}c \qquad c^\star = -qb \end{gathered} \tag{163}$$

The relations between generators of $\mathfrak{U}_q(su(2))$ and matrix elements of $\overset{\frac{1}{2}}{L}{}^{(\pm)}$

$$\overset{\frac{1}{2}}{L}{}^{(+)} = \begin{bmatrix} q^{J_z} & (1-q^{-2})J_- \\ 0 & q^{-J_z} \end{bmatrix} \quad \overset{\frac{1}{2}}{L}{}^{(-)} = \begin{bmatrix} q^{-J_z} & 0 \\ (1-q^{2})J_+ & q^{J_z} \end{bmatrix}. \tag{164}$$

The mixed relations of the quantum double are

$$\begin{array}{llll} q^{J_z}c = qcq^{J_z} & q^{J_z}b = q^{-1}bq^{J_z} & [q^{J_z},a]=0 & [q^{J_z},d]=0 \\ [J_+,c]=0 & [J_-,b]=0 & [J_+,b]=q^{-1}(q^{J_z}a - q^{-J_z}d) & [J_-,c]=q(q^{J_z}d-q^{-J_z}a) \\ J_-a = q^{-1}aJ_- + bq^{J_z} & aJ_+ = qJ_+a + q^{-J_z}c & dJ_- = q^{-1}J_-d + q^{-J_z}b & J_+d = qdJ_+ + cq^{J_z}. \end{array} \tag{165}$$

The application of the factorization theorem gives an explicit isomorphism of complex algebras between $\mathcal{D}$ and $\mathfrak{U}_q(sl(2,\mathbb{C}))\otimes\mathfrak{U}_q(sl(2,\mathbb{C}))$. If we denote

$$\overset{\frac{1}{2}}{M}{}^{(i+)} = \begin{bmatrix} q^{J_z^{(i)}} & (1-q^{-2})J_-^{(i)} \\ 0 & q^{-J_z^{(i)}} \end{bmatrix} \quad \overset{\frac{1}{2}}{M}{}^{(i-)} = \begin{bmatrix} q^{-J_z^{(i)}} & 0 \\ (1-q^{2})J_+^{(i)} & q^{J_z^{(i)}} \end{bmatrix}, \tag{166}$$

then (59,60) gives the explicit change of generators:

$$J_z = J_z^{(l)} + J_z^{(r)} \qquad J_- = q^{J_z^{(r)}}J_-^{(l)} + J_-^{(r)}q^{-J_z^{(l)}} \qquad J_+ = q^{J_z^{(r)}}J_+^{(l)} + J_+^{(r)}q^{-J_z^{(l)}} \tag{167}$$

$$a = q^{J_z^{(l)}-J_z^{(r)}} \quad b=(1-q^{-2})q^{-J_z^{(r)}}J_-^{(l)} \quad c=(1-q^{2})q^{J_z^{(l)}}J_+^{(r)} \quad d = q^{-J_z^{(l)}+J_z^{(r)}} - (q-q^{-1})^2 J_+^{(r)}J_-^{(l)}$$

The expression of the Casimir elements $\Omega^{\pm}$ is given by:

$$\Omega^+ = q(q-q^{-1})J_+b + q^{-1}q^{J_z}a + qq^{-J_z}d \tag{168}$$

$$\Omega^- = -q^{-1}(q-q^{-1})J_-c + q^{-1}q^{-J_z}a + qq^{J_z}d \tag{169}$$

From these commutation relations it is easy to show that the precise relation between the generators $\{J^{\pm}, q^{J_z}, a, b, c, d\}$ and the generators $\{A, N, N^{\star}, \gamma, \gamma^{\star}, \alpha, \alpha^{\star}\}$ used in [36][35] is:

$$A = q^{J_z}, N = J^+, N^{\star} = q^{-1}J^-, \gamma = b, \alpha = d, \gamma^* = -q^{-1}c, \alpha^{\star} = a. \tag{170}$$

As a result the relation between the Casimir element $\Omega^+$ and the Casimir element $X$ of [35] is $-q\Omega^+ = \sqrt{1+q^2}X$.

## 7.2 Clebsch-Gordan coefficients, and Vertex pictures

### Definitions

In the sequel we shall frequently use the following notations:

$$[P]_k = \prod_{i=1}^{k} [P+i-1]$$

The $R$ matrix elements and the Clebsh-Gordan coefficients can be computed and are given by (see [22] for a precise derivation)

$$\overset{AB}{R}{}^{n_1\; n_2}_{n_1-n\; n_2-n} = \frac{(q-q^{-1})^n}{[n]!} q^{2n_1 n_2 - n(n_2-n_1+\frac{n+1}{2})} \sqrt{\frac{[A-n_1]![A+n_1+n]![B+n_2]![B-n_2+n]!}{[A-n_1-n]![A+n_1]![B+n_2-n]![B-n_2]!}}$$

$$\left(\begin{matrix} m & n \\ I & J \end{matrix} \middle| \begin{matrix} K \\ p \end{matrix}\right) = \delta_{m+n,p} q^{m(p+1)+\frac{1}{2}(J(J+1)-I(I+1)-K(K+1))} e^{i\pi(I-m)} \sqrt{\frac{[2K+1][I+J-K]![I-m]![J-n]![K-p]![K+p]!}{[K+J-I]![I+K-J]![I+J+K+1]![I+m]![J+n]!}} \times$$
$$\times \sum_{V=0}^{K-p} \frac{q^{V(K+p+1)} e^{i\pi V} Y(I,J,K)[I+m+V]![J+K-m-V]!}{[V]!\,[K-p-V]!\,[I-m-V]!\,[J-K+m+V]!}.$$

### Pictorial representation

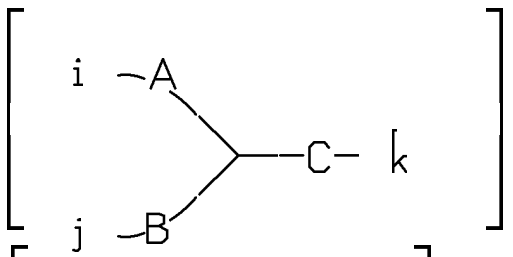

$= \left(\begin{matrix} i & j \\ A & B \end{matrix} \middle| \begin{matrix} C \\ k \end{matrix}\right)$

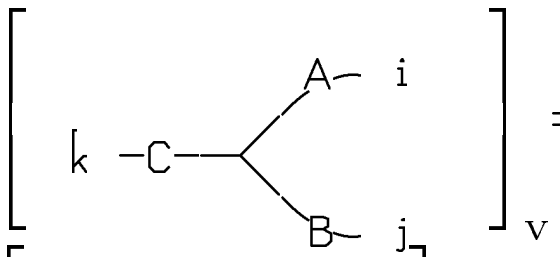

$= \left(\begin{matrix} k \\ C \end{matrix} \middle| \begin{matrix} A & B \\ i & j \end{matrix}\right)$

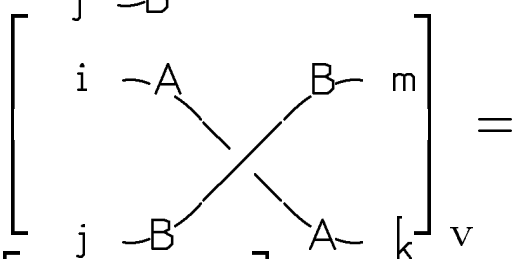

$= \overset{AB}{R}{}^{-1\, ij}_{\quad km}$

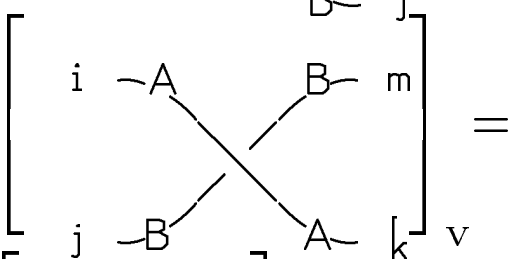

$= \overset{AB}{R'}{}^{ij}_{km} = \overset{BA}{R}{}^{ji}_{mk}$

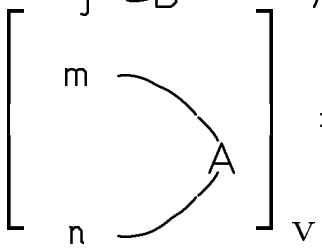

$= \delta_{m,-n} e^{-i\pi m} q^m = v_A^{\frac{1}{2}} \overset{A}{w}{}^{mn}$

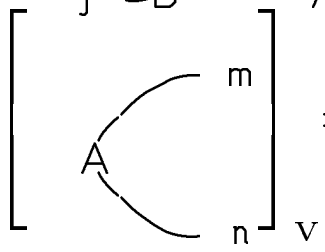

$= \delta_{m,-n} e^{-i\pi m} q^m = v_A^{\frac{1}{2}} \overset{A}{w}{}_{mn}$

### Relations

Unitarity:

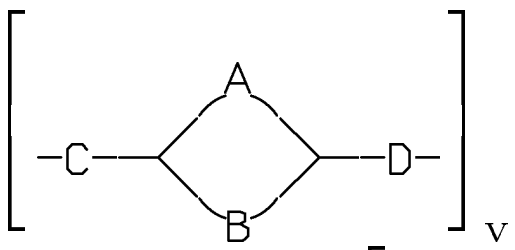

$= \delta_{C,D} Y(A,B,C)$

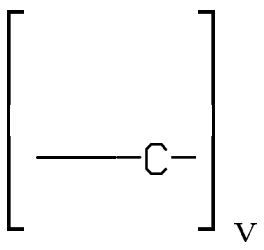

, $\sum_C$

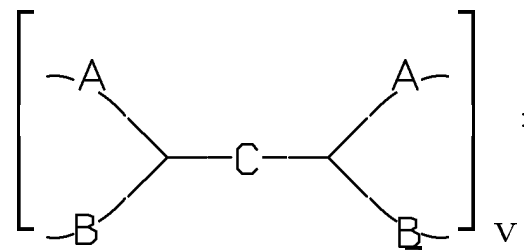

$=$

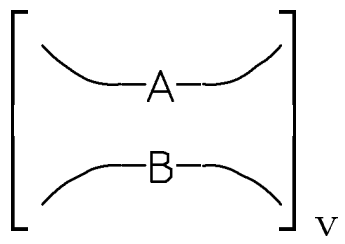

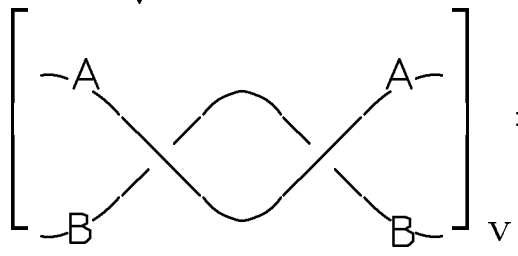

$=$

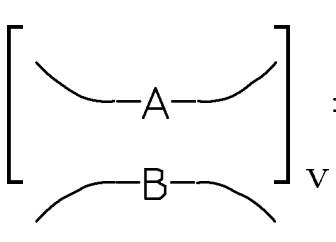

$=$

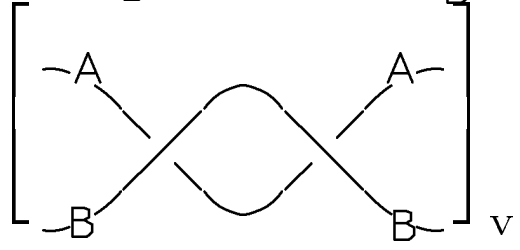

Twist:

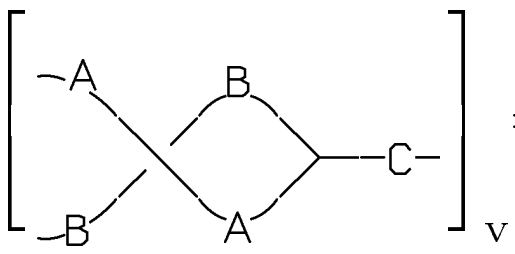

$= (\frac{v_A v_B}{v_C})^{\frac{1}{2}}$

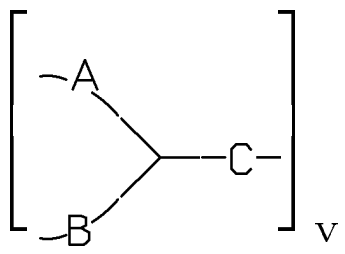

,

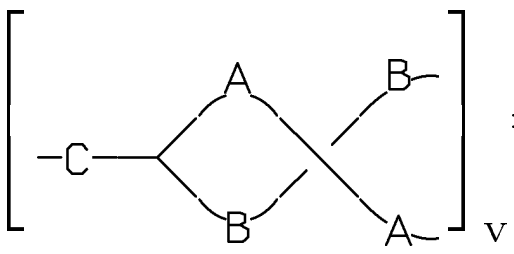

$= (\frac{v_A v_B}{v_C})^{\frac{1}{2}}$

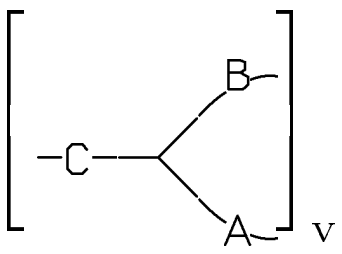

Contragredient:

$$\Big[$$ 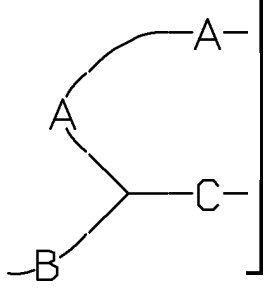 $$\Big]_{V} = e^{i\pi(-C+B)}\Big(\frac{[d_C]}{[d_B]}\Big)^{\frac{1}{2}} \Big[$$ 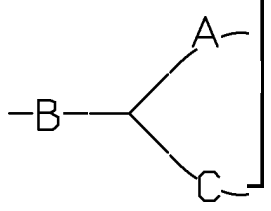 $$\Big]_{V} ,\qquad \Big[$$ 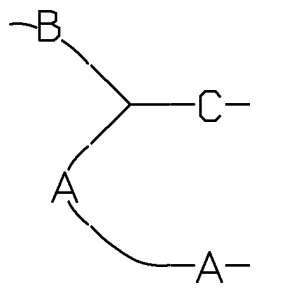 $$\Big]_{V} = e^{i\pi(-C+B)}\Big(\frac{[d_C]}{[d_B]}\Big)^{\frac{1}{2}} \Big[$$ 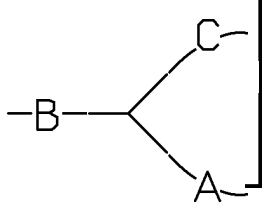 $$\Big]_{V}$$

Normalization:

$$\Big[$$ 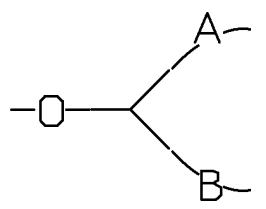 $$\Big]_{V} = \frac{e^{i\pi B}\delta_{B,A}}{\sqrt{[d_B]}} \Big[$$ 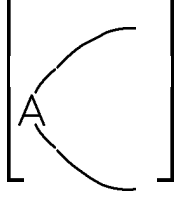 $$\Big]_{V} \qquad \Big[$$ 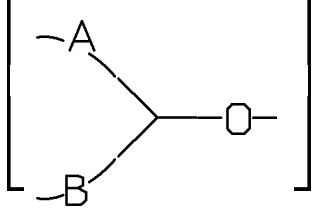 $$\Big]_{V} = \frac{e^{i\pi B}\delta_{B,A}}{\sqrt{[d_B]}} \Big[$$ 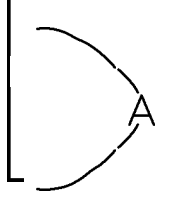 $$\Big]_{V}$$

One has to add to these relations those coming from elementary topological moves on ribbon graphs, which come from Yang-Baxter equation and quasitriangularity equations.

## Definitions of $6j$ symbols, relations and formulas

$$\Big[$$ 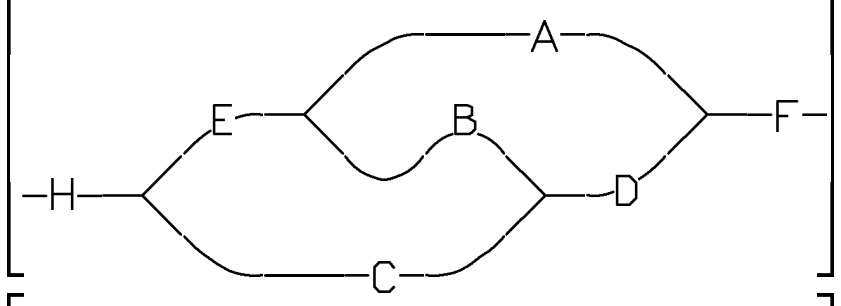 $$\Big]_{V} = \left\{\begin{matrix} A & B \\ C & F \end{matrix}\middle|\begin{matrix} E \\ D \end{matrix}\right\} \delta_{F,H} \Big[$$ 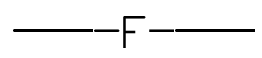 $$\Big]_{V}$$

$$\Big[$$ 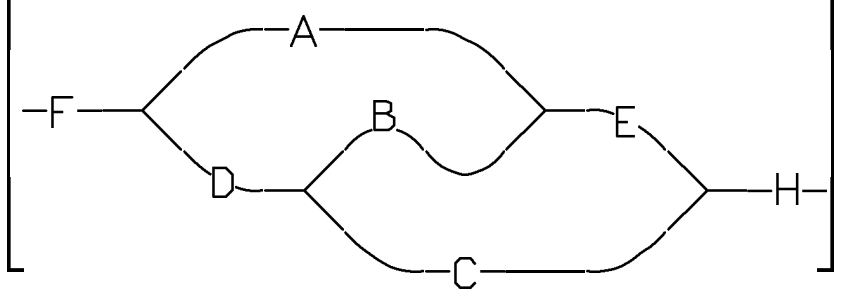 $$\Big]_{V} = \left\{\begin{matrix} C & B \\ A & F \end{matrix}\middle|\begin{matrix} D \\ E \end{matrix}\right\} \delta_{F,H} \Big[$$ 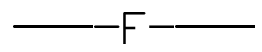 $$\Big]_{V}$$

$$\left\{\begin{matrix} A & B \\ D & C \end{matrix}\middle|\begin{matrix} E \\ F \end{matrix}\right\} = e^{-i\pi(C+D+2E-A-B)}([2E+1][2F+1])^{\frac{1}{2}}\Delta(A,B,E)\Delta(A,C,F)\Delta(C,E,D)\Delta(D,B,F)\times$$

$$\times\sum_{\substack{U\in\mathbb{Z} \\ \left.\begin{matrix}A+C+F\\A+B+E\\B+D+F\\D+C+E\end{matrix}\right\}\le U\le\left\{\begin{matrix}A+B+C+D\\A+D+E+F\\B+C+E+F\end{matrix}\right.}} \frac{e^{i\pi U}\,Y(A,B,E)Y(A,C,F)Y(C,E,D)Y(D,B,F)\ [U+1]!}{[U-A-C-F]!\,[U-A-B-E]!\,[U-B-D-F]!\,[U-D-C-E]!\,[A+B+C+D-U]!\,[A+D+E+F-U]!\,[B+C+E+F-U]!}$$

$$\text{with}\qquad \Delta(I,J,K)=\sqrt{\frac{[J+K-I]!\,[J+I-K]!\ [I+K-J]!}{[I+J+K+1]!}}\quad \forall I,J,K\in\frac{1}{2}\mathbb{Z}^{+}. \tag{171}$$

The $6j$ symbols satisfy the symmetry relations:

$$\left\{\begin{matrix} A & B \\ C & D \end{matrix}\middle|\begin{matrix} E \\ F \end{matrix}\right\} = \left\{\begin{matrix} B & A \\ D & C \end{matrix}\middle|\begin{matrix} E \\ F \end{matrix}\right\} = \left\{\begin{matrix} C & D \\ A & B \end{matrix}\middle|\begin{matrix} E \\ F \end{matrix}\right\} = \left\{\begin{matrix} A & D \\ C & B \end{matrix}\middle|\begin{matrix} F \\ E \end{matrix}\right\} \tag{172}$$

They also satisfy the following relation, known as Racah-Wigner relation:

$$\left\{\begin{matrix} A & C \\ B & E \end{matrix}\middle|\begin{matrix} F \\ D \end{matrix}\right\} = (-1)^{C-F+E-D}\left\{\begin{matrix} A & F \\ B & D \end{matrix}\middle|\begin{matrix} C \\ E \end{matrix}\right\}\Big(\frac{[d_F][d_D]}{[d_C][d_E]}\Big)^{\frac{1}{2}} \tag{173}$$

which is nicely pictured in the next subsection.
We give the value of 6j symbols for one spin fixed to the value 0 or $\frac{1}{2}$ :

$$\left\{\begin{matrix} 0 & B \\ A & C \end{matrix}\middle|\begin{matrix} B \\ C \end{matrix}\right\} = Y(A,B,C) \qquad \left\{\begin{matrix} A & A \\ B & B \end{matrix}\middle|\begin{matrix} 0 \\ C \end{matrix}\right\} = Y(A,B,C)e^{i\pi(A+B-C)}\sqrt{\frac{[d_C]}{[d_A][d_B]}} \quad (174)$$

$$\begin{aligned}
\left\{\begin{matrix} A & B \\ \frac{1}{2} & E+\frac{1}{2} \end{matrix}\middle|\begin{matrix} E \\ B+\frac{1}{2} \end{matrix}\right\} &= (\frac{[B+E-A+1][A+B+E+2]}{[2E+2][2B+1]})^{\frac{1}{2}}Y(A,B,E) \\
\left\{\begin{matrix} A & B \\ \frac{1}{2} & E+\frac{1}{2} \end{matrix}\middle|\begin{matrix} E \\ B-\frac{1}{2} \end{matrix}\right\} &= -(\frac{[A+B-E][A+E-B+1]}{[2E+2][2B+1]})^{\frac{1}{2}}Y(A,B,E) \\
\left\{\begin{matrix} A & B \\ \frac{1}{2} & E-\frac{1}{2} \end{matrix}\middle|\begin{matrix} E \\ B+\frac{1}{2} \end{matrix}\right\} &= (\frac{[A+E-B][A+B-E+1]}{[2E][2B+1]})^{\frac{1}{2}}Y(A,B,E) \\
\left\{\begin{matrix} A & B \\ \frac{1}{2} & E-\frac{1}{2} \end{matrix}\middle|\begin{matrix} E \\ B-\frac{1}{2} \end{matrix}\right\} &= (\frac{[A+B+E+1][E+B-A]}{[2E][2B+1]})^{\frac{1}{2}}Y(A,B,E) \qquad (175)
\end{aligned}$$

## 7.3 $6-j$ coefficients and IRF pictures

## Pictorial representation

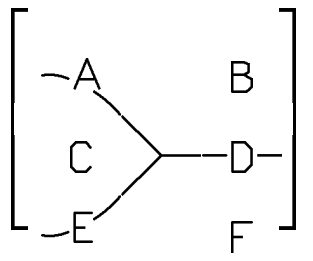

$= \left\{\begin{matrix} E & A \\ B & F \end{matrix}\middle|\begin{matrix} D \\ C \end{matrix}\right\}$,

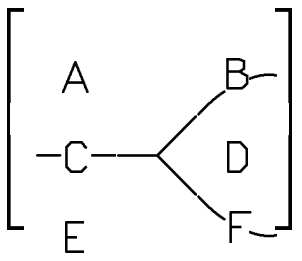

$= \left\{\begin{matrix} F & B \\ A & E \end{matrix}\middle|\begin{matrix} C \\ D \end{matrix}\right\}$,

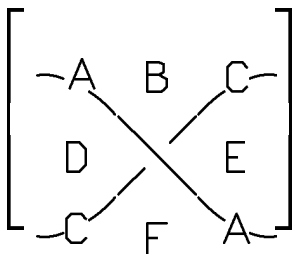

$= (\frac{v_D v_E}{v_B v_F})^{\frac{1}{2}} \left\{\begin{matrix} A & B \\ C & F \end{matrix}\middle|\begin{matrix} D \\ E \end{matrix}\right\}$,

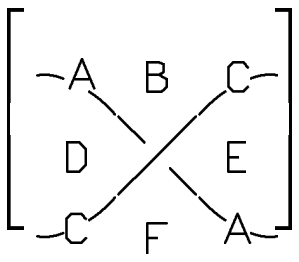

$= (\frac{v_B v_F}{v_D v_E})^{\frac{1}{2}} \left\{\begin{matrix} A & B \\ C & F \end{matrix}\middle|\begin{matrix} D \\ E \end{matrix}\right\}$,

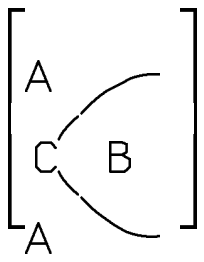

$= e^{i\pi(A-B)}(\frac{[d_B]}{[d_A]})^{\frac{1}{2}}Y(A,B,C)$

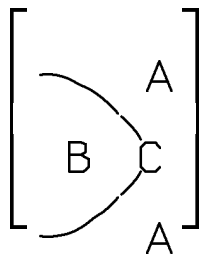

$= e^{i\pi(A-B)}(\frac{[d_B]}{[d_A]})^{\frac{1}{2}}Y(A,B,C)$

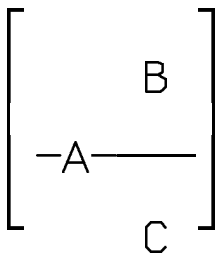

$= Y(A,B,C)$

### Relations

Unitarity:

$\sum_F$ 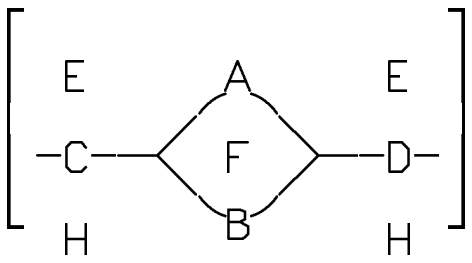

$= \delta_{C,D}\, Y(A,B,C)$ 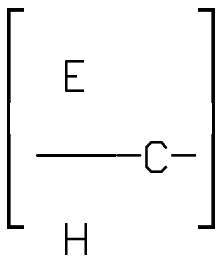

, $\sum_C$ 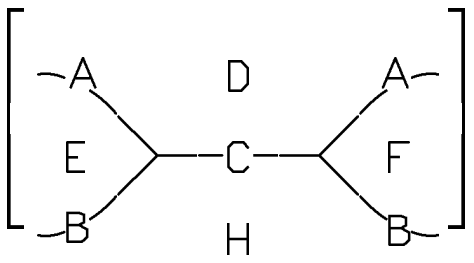

$= \delta_{E,F}$ 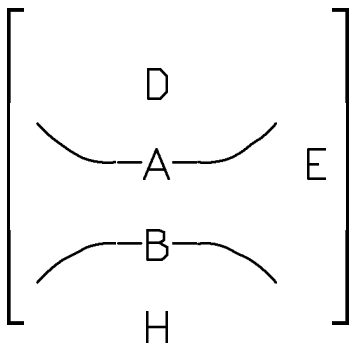

$\sum_F$ 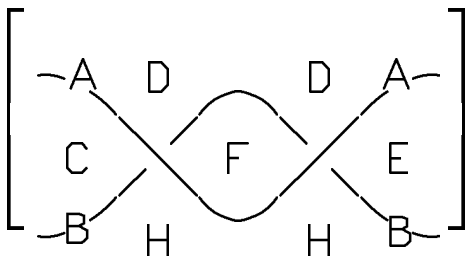

$= \delta_{C,E}$ 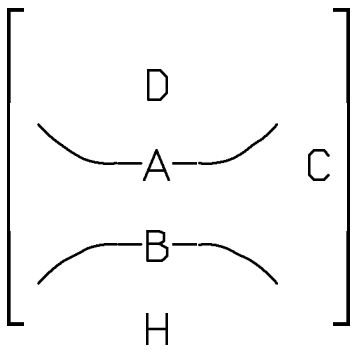

$= \sum_F$ 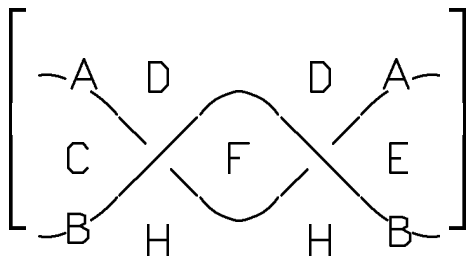

Racah relation:

$$\sum_{H}$$ 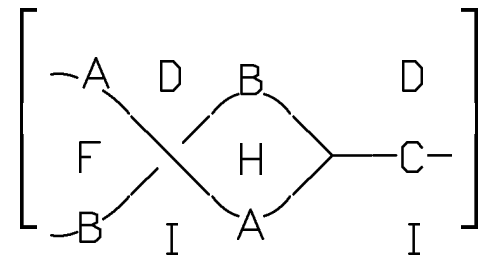 $$_{\text{IRF}} = \left(\frac{v_A v_B}{v_C}\right)^{\frac{1}{2}}$$ 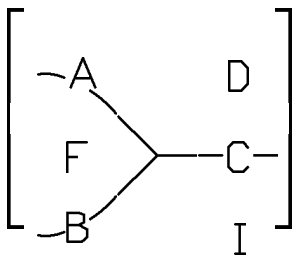 $$_{\text{IRF}}\ ,\qquad \sum_{E}$$ 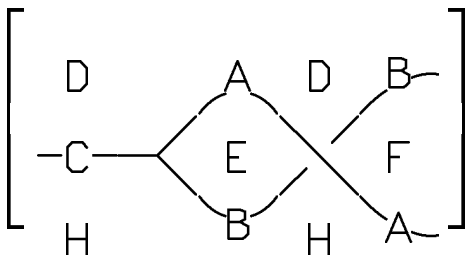 $$_{\text{IRF}} = \left(\frac{v_A v_B}{v_C}\right)^{\frac{1}{2}}$$ 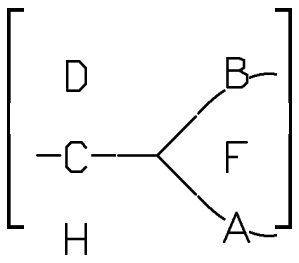 $$_{\text{IRF}}$$

Racah-Wigner Symmetry:

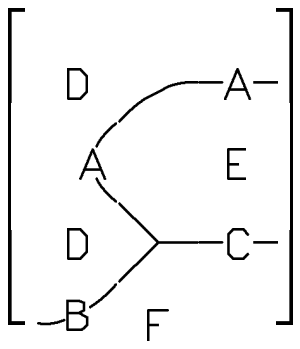 $$_{\text{IRF}} = e^{i\pi(-\text{C}+\text{B})}\left(\frac{[d_C]}{[d_B]}\right)^{\frac{1}{2}}$$ 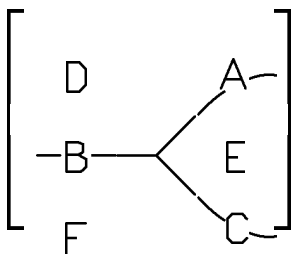 $$_{\text{IRF}}\ ,$$ 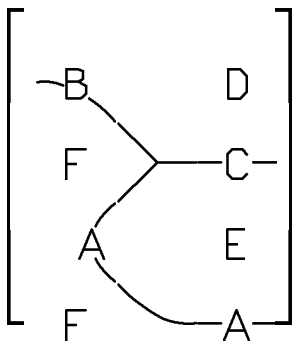 $$_{\text{V}} = e^{i\pi(-\text{C}+\text{B})}\left(\frac{[d_C]}{[d_B]}\right)^{\frac{1}{2}}$$ 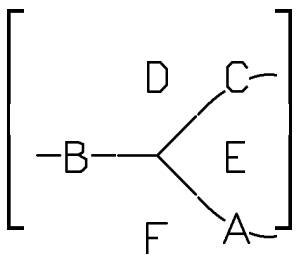 $$_{\text{IRF}}$$

Normalization:

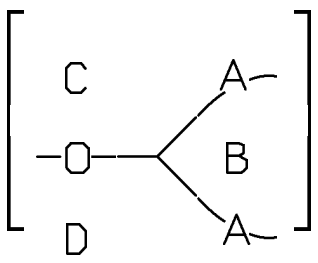 $$_{\text{V}} = \frac{e^{i\pi \text{A}}\delta_{D,C}}{\sqrt{[d_A]}}$$ 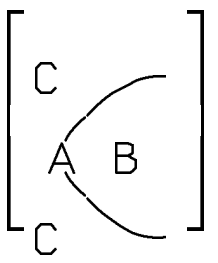 $$_{\text{V}}$$ 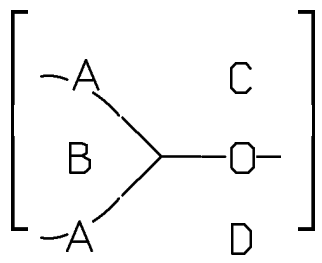 $$_{\text{V}} = \frac{e^{i\pi \text{A}}\delta_{D,C}}{\sqrt{[d_A]}}$$ 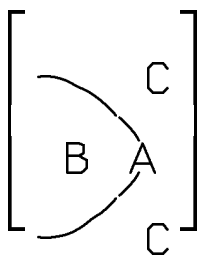 $$_{\text{V}}$$

One has to add to these relations the pentagonal relation, and the hexagonal one. It is important to remark that the ranges of summations in these relations are simply given by the $Y$ factors in the definitions of $6j$ coefficients entering in the formulas associated to the pictorial representation

## 7.4 Complex continuations of $6j$ symbols

The constructions described in this paper make an extensive use of complex continuation of quantum $6j$ symbols. The pionnering work on this subject is due to R.Askey and J.Wilson [5], and has been analyzed in great details in the works of E.Cremmer, J.L.Gervais, J.F.Roussel [14, 23]. In this appendix we just want to recall the basic ideas of their construction. Ordinary $6j$ coefficients can be rewritten using ${}_4F_3$ basic hypergeometric functions (see [23]), so it could seem natural to study the continuation from this expression. But the proof of polynomial equations usually verified by $6j$ coefficients is not so obvious in this formulation. In the work [14, 23] another method is developed, which consists in a step by step complex continuation of $6j$ coefficients from the situation where all spins are half integers, and where the polynomial equations are simply derived from features of the representation theory of the associated quantum group [27].

This continuation is done in different steps, depending on the number of combination of spins we want to let be half integers. For our purpose, the first step will be sufficient. The central idea is that the definition of usual $6j$ coefficients (171) requires as a basic condition, that some combinations of spins be integers, but does not require the spins to be half integers themselves. In order to give sense to this idea [14, 23] use the following definitions

$$\Theta_{J\ K}^{\ I} = \sqrt{\frac{[J+I-K+1]_{J+K-I}\,[K+I-J+1]_{J+K-I}\,[2I+2]_{J+K-I}}{[J+K-I]!}}$$

where $I, J, K$ can be complex numbers provided that $J+K-I \in \mathbb{Z}^+$.

$$Y(A,\aleph+N,\aleph+P)=\begin{cases} 1 & \text{if } A\pm(N-P)\in\mathbb{N} \\ 0 & \text{elsewhere}\end{cases}$$

In order to define the $6j$ coefficient of the form $\left\{\begin{matrix} B & C \\ \aleph+N & \aleph \end{matrix}\middle| \begin{matrix} A \\ \aleph+M \end{matrix}\right\}$ (where $A,B,C \in \frac{1}{2}\mathbb{Z}^+, M,N \in \mathbb{Z}, \aleph \in \mathbb{C} - \frac{1}{2}\mathbb{Z}$) we will proceed to the variable change

$$T = U - 2\aleph$$

in the summation in (171). In order to avoid the use of $q-\Gamma$ functions, we remark that the following simple equalities hold

$$\Theta_{I\,J}^{\ K} = \sqrt{\frac{[2I]!\,[2J]!}{[2K+1]!}}\frac{1}{\Delta(I,J,K)}, \qquad \frac{[U+1]!}{[U-I-J-K]!} = [U-I-J-K+1]_{I+J+K+1},$$

$$\Theta_{J\,\aleph+m}^{\ \aleph} = \sqrt{\frac{[2J]!}{[J+m]!\ [J-m]!}}\sqrt{[2\aleph+m-J+1]_{J+m}\,[2\aleph+2]_{J+m}} \qquad \Delta(I,J,J+K)^2 = \frac{[I+K]!\,[I-K]!}{[2J+K-I+1]_{2I+1}}.$$

We can therefore define the first type of one parameter complex continuation of $6j$ by

$$\left\{\begin{matrix} B & C \\ \aleph+N & \aleph \end{matrix}\middle| \begin{matrix} A \\ \aleph+M \end{matrix}\right\} = e^{i\pi(C+B-2A-N)}\frac{\Theta_{C\,\aleph+N}^{\ \aleph+M}\ \Theta_{B\,\aleph+M}^{\ \aleph}}{\Theta_{B\ C}^{\ A}\ \Theta_{A\,\aleph+N}^{\ \aleph}}[2\aleph+2M+1]\frac{[C+N-M]!\,[C+M-N]!}{[2\aleph+M+N-C+1]_{2C+1}}\frac{[B+M]!\,[B-M]!}{[2\aleph+M-B+1]_{2B+1}} \times$$
$$\times \sum_{\substack{T\in\mathbb{Z} \\ \left.\begin{matrix} B+M \\ A+N \\ C+M+N\end{matrix}\right\}\le T\le\left\{\begin{matrix} B+C+N \\ B+A+M+N \\ C+A+M\end{matrix}\right.}} \frac{Y(A,B,C)\,Y(B,\aleph,\aleph+M)\,Y(C,\aleph+N,\aleph+M)\,Y(A,\aleph+N,\aleph)\ e^{i\pi T}\ [T+2\aleph-B-C-A+1]_{A+B+C+1}}{[T-A-N]!\,[T-B-M]!\,[T-C-M-N]!\,[B+C+N-T]!\,[B+A+M+N-T]!\,[C+A+M-T]!} \tag{176}$$

and observe that this definition has the same form as (171), except for two features: the definition of $Y$ is changed (i.e. one of the three constraints is relaxed), and the summation forgets one of the constraints entering in the usual definition (i.e. $A+B+C \le U$). In fact, it was necessary to proceed to these changes because, if $\aleph$ is in $\mathbb{C} - \frac{1}{2}\mathbb{Z}^+$, these constraints do not have sense anymore. Nevertheless, if $\aleph$ is an half integer, sufficiently large, the forgotten constraints are automatically verified, once the other constraints are fullfilled. Hence, in this case, both definitions are equal. In this sense, we can speak of "continuation" to complex spins. Up to this point, we have defined continuation of $6j$ coefficients with three complex spins, depending just on one arbitary complex number, located in the lower row of the $6j$ symbol. The definition of the continuation $6j$ coefficient is the unique one which maintains the usual trivial symmetries (172).

In the same way, [14, 23] define another one parameter complex continuation of $6j$ coefficients with four complex entries, depending just on one arbitrary complex number :

$$\left\{\begin{matrix} A & \aleph+N \\ B & \aleph \end{matrix}\middle| \begin{matrix} \aleph+P \\ \aleph+M \end{matrix}\right\} = e^{i\pi(A-B+N-2P)}\frac{\Theta_{B\,\aleph+N}^{\ \aleph+M}\ \Theta_{A\,\aleph+M}^{\ \aleph}}{\Theta_{A\,\aleph+N}^{\ \aleph+P}\ \Theta_{B\,\aleph+P}^{\ \aleph}}[2\aleph+2M+1]\frac{[B+N-M]!\,[B+M-N]!}{[2\aleph+M+N-B+1]_{2B+1}}\frac{[A+M]!\,[A-M]!}{[2\aleph+M-A+1]_{2A+1}} \times$$
$$\times \sum_{\substack{T\in\mathbb{Z} \\ \left.\begin{matrix} A+M \\ A+N+P \\ B+P \\ B+M+N\end{matrix}\right\}\le T\le\left\{\begin{matrix} B+A+N \\ B+A+M+P\end{matrix}\right.}} \frac{Y(A,\aleph+N,\aleph+P)\,Y(A,\aleph,\aleph+M)\,Y(B,\aleph+N,\aleph+M)\ Y(B,\aleph,\aleph+P)\ e^{i\pi T}\ [2\aleph+N+M+P-T+1]_{2T-N-M-P}}{[T-A-M]!\,[T-A-N-P]!\,[T-B-M-N]!\,[T-B-P]!\,[B+A+N-T]!\,[B+A+M+P-T]!} \tag{177}$$

The previous discussion can be repeated in this case, and the reader is invited to read the details in [23]. As in the previous case, the continuation of $6j$ coefficients obtained from the last one by the usual trivial symmetries (172) can be defined in the same way and these symmetries are again verified.

At this point, we have just given definitions and it is checked that the usual properties of $6j$ coefficients (such as Unitarity, Racah , Pentagonal, or Hexagonal relations) are verified [14, 23].

It is easy to compute using the formula (176) the $6j$ symbols when one of the spins is equal to 0 or $\frac{1}{2}$. The values of these coefficients are obtained from (174,175) by replacing $A = \aleph + N, B = \aleph$ in the formulas with $N \in \frac{1}{2}\mathbb{Z}$ and leaving $C, E$ positive half integers.